\documentstyle[twoside,titlepage]{report}

\newtheorem{thm}{\indent{\sc Theorem}}[section]
\newtheorem{defn}[thm]{\indent{\sc Definition}} 
\newtheorem{lem}[thm]{\indent{\sc Lemma}} 
\newtheorem{por}[thm]{\indent{\sc Porism}}
\newtheorem{prop}[thm]{\indent{\sc Proposition}} 
\newtheorem{cor}[thm]{\indent{\sc Corollary}}
\newtheorem{ex}{\indent {\sc Example}}[section]
\newtheorem{prob}[thm]{\indent{\sc Open Problem}\index{Problem{,} Open}} 

\newcommand{\sect}{\cleardoublepage\section}
\newcommand{\ssect}{\subsection}
\newcommand{\sssect}{\subsubsection}

\newcommand{\eqref}[1]{equation~(\ref{#1})}
\newcommand{\eref}[1]{\eqref{#1}}

\newcommand{~}{\nolinebreak[3] }

\newcommand{\thin}[1]{\makebox[0pt]{#1}}

\newcommand{\th}{\raisebox{0.6ex}{th}}
\newcommand{\del}{\nabla}
\newcommand{\seper}{\nopagebreak \begin{center}
	\underline{\hspace{2in}}
	\end{center}}

\newcommand{\emp}[1]{{\em #1}\index{#1}}
\newcommand{\bold}[1]{{\em (#1)}\index{#1}}
\newcommand{\bld}[1]{{\em (#1)}}

\newcommand{\dc}[4]{{\samepage$$\left\{ \begin{array}{rrcl}
{\bf Discrete} &\displaystyle #1 &=& \displaystyle  #2 \\*[0.1in]
{\bf Continuous} &\displaystyle #3 &=& \displaystyle  #4 
\end{array}\right\}$$}}
\newcommand{\dcthree}[6]{{\samepage$$\left\{ \begin{array}{rrcccl}
{\bf Discrete} &\displaystyle #1 &=& \displaystyle  #2 &=& \displaystyle
#3\\*[0.1in] 
{\bf Continuous} &\displaystyle #4 &=& \displaystyle  #5 &=& \displaystyle #6
\end{array}\right\}$$}}
\newcommand{\dclab}[5]{{\samepage \begin{equation}\label{#1}
\left\{ \begin{array}{rrcl}
{\bf Discrete} &\displaystyle #2 &=& \displaystyle  #3 \\*[0.1in]
{\bf Continuous} &\displaystyle #4 &=& \displaystyle  #5
\end{array}\right\}
\end{equation}}}

 \newcommand{\proof}[1]{\proo{#1}$\Box $}
\newcommand{\proo}[1]{{\em Proof: }#1}
\renewcommand{\Box}{ \hspace{1cm}{\rule{1.2ex}{2ex}}}

\newcommand{\rn}[1]{\left\lfloor #1 \right\rceil}

\newcommand{\smallD}{\mbox{{\small D}}}
\newcommand{\smallI}{{\mbox{\small I}}}
\newcommand{\D}{{{\rm D}}} 
 
\newcommand{\I}{{\rm I}}

\newcommand{\proj}{\mbox {\rm proj}}

\newcommand{\adj}[1]{\mbox {\rm adj}(#1)}

\newcommand{\romcoeff}[2]{\rn{{#1 \atop #2}}} 
 
\newcommand{\lr}{\langle x \rangle}
\newcommand{\Res}[1]{\mbox{\rm Res}(#1)} 
\newcommand{\action}[2]{\left\langle #2 \right\rangle 
 	\!\!_{_{\scriptstyle #1}}}

\newcommand{\phil}{\nonumber \\*[-1mm]  &&\: }

\newcommand{\SB}[1]{\mbox{{\bf \scriptsize #1}}}
\newcommand{\seq}[3]{#1_{#2}^{#3}(x)}
\newcommand{\comp}[2]{#1\left({\bf #2}\right)}
\newcommand{\digress}[1]{ \begin{quotation}
	 #1 {\em End of Digression.} \end{quotation}}
\newcommand{\discrete}[1]{\begin{description}
	\item[Discrete] #1 \end{description}}
\newcommand{\cont}[1]{\begin{description}
	\item[Continuous] #1 \end{description}}

\renewcommand{\tilde}[1]{\widetilde{#1}}
  
\newcommand{\fig}[2]{\begin{table}[htbp]
	\caption{#1} 
	\begin{center}\fbox{ \scriptsize \begin{tabular}#2 \end{tabular}}
	\end{center}
\end{table}}

\newcommand{\figx}[2]{\begin{figure}[htbp] \caption{#1} 
	\begin{center} \scriptsize \mbox{}#2 \end{center} \end{figure}}

\begin{document}

\newcommand{\thmplus}[2]{\renewcommand{\thethm}{{\bf 
	\thesection{}.\arabic{thm}{\bf #1}}} #2 
	\renewcommand{\thethm}{{\bf \thesection.\arabic{thm}}}}
\renewcommand{\thethm}{{\bf \thesection.\arabic{thm}}}
\newcommand{\explus}[2]{\renewcommand{\theex}{{\bf 
	\thesection{}.\arabic{ex}{\bf #1}}} #2 
	\renewcommand{\theex}{{\bf \thesection.\arabic{ex}}}}
\renewcommand{\theex}{{\bf \thesection.\arabic{ex}}}

\newcommand{\cref}[1]{chapter~\ref{#1}}
\newcommand{\Cref}[1]{Chapter~\ref{#1}}

\addtocontents{toc}{Please note: Sections marked {\bf D} apply only to the
discrete 
theory. Sections marked {\bf C} apply only to the continuous theory. All other
sections apply equally to both theories. Sections marked {\bf A} are
appendicies, and are generally independent of all later material.}

\renewcommand{\sect}[1]{\asectplus{#1}{}{#1}}
\newcommand{\sectplus}[2]{\asectplus{#1}{#2}{#1}}
\newcommand{\sapp}[1]{\asectplus{#1}{A}{Appendix: #1}}
\newcommand{\asectplus}[3]{\renewcommand{\thechapter}{\arabic{chapter}{\bf #2}}
	\cleardoublepage \chapter[#3]{#1} \indent}
\renewcommand{\ssect}[1]{\assectplus{#1}{}{#1}}
\newcommand{\ssectplus}[2]{\assectplus{#1}{#2}{#1}}
\newcommand{\ssapp}[1]{\assectplus{#1}{A}{Appendix: #1}}
\newcommand{\assectplus}[3]{\renewcommand{\thesection}
	{\thechapter.\arabic{section}{\bf #2}}
	\nopagebreak \section[#3]{#1} \nopagebreak \indent}
\renewcommand{\sssect}[1]{\asssectplus{#1}{}{#1}}
\newcommand{\sssectplus}[2]{\asssectplus{#1}{#2}{#1}}
\newcommand{\sssapp}[1]{\asssectplus{#1}{A}{Appendix: #1}}
\newcommand{\asssectplus}[3]{\renewcommand{\thesubsection}
	{\thesection.\arabic{subsection}{\bf #2}}
	\nopagebreak \subsection[#3]{#1} \nopagebreak  \indent}
\renewcommand{\numberline}[1]{\makebox[0pt][l]{#1}\hspace{1.0in}}
\renewcommand{\tilde}[1]{\widetilde{#1}}

\begin{abstract}
We generalize the Umbral Calculus of G-C. Rota \cite{RR} by 
studying not only sequences of polynomials and inverse power series, or even
the logarithms studied in \cite{LR}, but instead we study sequences of formal
expressions involving the iterated logarithms and $x$ to an arbitrary real
power.

Using a theory of formal power series with real exponents, and a more general 
 definition of factorial, binomial
coefficient, and Stirling numbers to all the real numbers,
we define  the Iterated Logarithmic Algebra
${\cal I}$. Its elements are the formal representations of the asymptotic
expansions of a large class of real functions, and we define the harmonic logarithm
 basis of ${\cal I}$ which  will be interpreted as a generalization of
the powers $x^{n}$ since it behaves nicely with respect to the derivative.

We classify all operators over ${\cal I}$ which commute with the derivative
(classically these are known as shift-invariant operators), and formulate
several equivalent definitions of a sequence of binomial type. 
We then derive many formulas useful towards the calculation of these sequences
including the Recurrence Formula, the Transfer Formula, and the Lagrange
Inversion Formula. Finally, we study Sheffer sequences, and give many examples.
\seper
{\bf L'alg\`{e}bre des Logarithmes Iter\'{e}s}

On g\'{e}n\'{e}ralise ici le calcul ombral de G-C Rota en \'{e}tudiant
non seulement les suites de polyn\^{o}mes, de s\'{e}ries \`{a}
exposants entier n\'{e}gatif, ou encore de logarithmes \cite{LR}, mais
d'une fa\,con plus g\'en\'erale on s'interesse aux suites formelles
des logarithms it\'eres de $x$ et de $x$ \`a un exposant r\'eel
quelconque. 

On s'appuie sur  la th\'eorie des s\'eries de puissances formelles \'a
exposants r\'eels, sur une d\'efinition g\'en\'erale des factorielles,
des coefficients du bin\^omes, et des nombres de Stirling \`a tous les
nombres r\'eel pour d\'efinir l'Algebre de Logarithmes It\'er\'es
${\cal I}$.  Cette alg\`ebre a comme \'el\'ements les
r\'epresentations asymptotiques de beaucoup de fonctions r\'eelles par
rapport de l'echelle des mon\^omes dans les Logarithmes It\'er\'es. On
d\'efinit alors une base de ${\cal I}$ de logarithms harmoniques,
lesquels pourront \^etre consid\'er\'es comme g\'en\'eralisant les
puissances $x^n$, puisqu'ils se comportent comme $x^n$ lors de la
d\'erivation. 

On caracterise tous les op\'erations sur ${\cal I}$ qui commutent avec
l'op\'erateur de d\'erivation (ce sont usuellement les ``op\'erateurs
invariants par translation''), ce qui nous permet de formuler
plusieurs d\'efinition equivalentes des suite de series logarithmiques
de type binomial. On en d\'eduit de nombreuses formules utiles au
calcul de ces series, en particulier la formule de r\'ecurrence, la
formule de transfert, et la formule d'inversion de Lagrange. Enfin, on 
\'etudie    les suites du Sheffer, ainsi que de nombreux exemples.
\end{abstract}

\begin{center}
{\it Dedicated to\\
Professor Gian-Carlo Rota}
\end{center}

\tableofcontents

\sect{Introduction}

Over the years, mathematicians have studied many special sequences of functions
especially polynomials. These sequences were found to have many important
similarities; however, it was not until G-C. Rota's Umbral Calculus \cite{RR}
that this notion was formalized. Only could one study under a single theory
most of 
the important sequences of polynomials---for example, $x^{n}$, Abel, lower
factorial, upper factorial. These sequences
are all sequences of binomial type.

Nevertheless, yet more work needed to be done, for the theory only applied to
polynomials. For example, the theory of ``factor sequences'' was developed
specifically to handle inverse formal power series. However, this theory was
completely separate from the Umbral calculus of polynomials. 
An early version of this article was published under the title
{\em Formal Power Series of Logarithmic Type} \cite{LR}. Not only did this
paper allow us to simultaneously consider polynomials and inverse formal power
series,  but it also allowed us to consider expressions involving the
logarithm $\log x$. 

This article differs from \cite{LR} principally in that it considers two further
generalizations to the logarithmic algebra ${\cal L}$. 
By introducing
the Iterated logarithms $\log \log x$, $\log \log \log x$, and so on to the
Logarithmic algebra, we have the {\em discrete} theory\index{Discrete Theory}.
Then by   including $x$ to any real power, we derive  the {\em
continuous}\index{Continuous Theory}
theory. Thus, the logarithmic algebra ${\cal L}$ is the set of all complex
functions on the real numbers with asymptotic expansions in a neighborhood of
$+\infty $ with respect to the ladder of comparison
$x^{a}(\log x)^{b}\cdots $.
 With only a handful of obvious exceptions, all of the results of {\em
Formal Power Series of Logarithmic Type} carry over into these new contexts
with only minor changes.
For example, we conclude that
every sequence of polynomials binomial type can be uniquely extended to a 
pseudobasis of ${\cal L}$ called a Roman graded
sequence (after Prof. Steve 
Roman). In fact, this sequence is usually easier to calculate than the sequence
of binomial type itself using classical techniques. Thus, even when one is only
interested in polynomials it still pays to introduce logarithms.

\ssect{Discrete and Continuous}

In this paper, all results and sections
regarding the discrete iterated logarithmic algebra will be
denoted ``Discrete'' and the results regarding the more general
continuous iterated logarithmic algebra will be denoted
``Continuous.'' Readers interested in only one of these theories may safely
omit all material pertaining to the other.
Sections and results numbered with a {\bf D} and paragraphs starting with
{\bf Discrete}
are relevant only to the discrete theory whereas sections and results numbered
with a {\bf C} and paragraphs starting with {\bf Continuous} are relevant only
to the continuous theory. Otherwise, the remainder of this article may be
interpreted discretely by 
supposing that all variables $a,b,c,\ldots $ are integers, and that
$\alpha ,\beta $ are vectors of integers, or it may be
interpreted continuously by supposing that all variables
$a,b,c,\ldots $ are real numbers and that $\alpha ,\beta $
are vectors of real numbers.

Digressions and all sections marked with the word ``Appendix'' or
the letter ``{\bf{A}}'' are independent
of all later material, and are included for their own sake.

\include{log}
\sect{$\D $-Invariant Operators}\label{opsec}
\ssect{The Operator Topology}\label{OpTop}
\label{Euler-sec}

It is a classical result  that the algebra of
formal differential operators $\sum _{n\geq 0}a_{n}\D ^{n}$
acts on the vector space of polynomials. In
view of the fact that the derivative $\D $ is invertible in
${\cal I}^{\alpha}$ for $\alpha \neq (0)$
(Proposition~\ref{invert}), we can define on ${\cal I}^{+}$ the
action of a more general class of differential operators,
called Artinian operators. These are linear operators which
commute with all powers of the derivative, and act on each level in a similar
fashion; this notion will be made more precise later. 

The primary goal of this chapter is to classify all Artinian operators. We see
(Theorem~\ref{sio1}) that they are merely Artinian series in the 
derivative, and that the 
operator topology---which we are shortly to define---corresponds to the
Artinian topology. This implies (see \S\ref{diffEqSec}) that all linear
differential equations have 
a unique canonical solution. Finally, along the way we find the
opportunity (in \S\ref{RA}) to apply this theory to Real Analysis.

We begin by defining a topology on the ring of continuous
linear operators acting on formal power series of
logarithmic type. 

Let ${\cal J}$ be a  subspace of ${\cal I}$. We say
that a sequence $(\theta _{n})_{n\geq 0}$ of continuous
linear operators of ${\cal J}$ into itself converges in the
{\em operator topology}\index{Operator Topology} on ${\cal J}$ when for every
$p(x)\in {\cal J}$ the sequence $(\theta _{n}p(x))_{n\geq 0}$ converges in
${\cal J}$. 

\begin{prop}
Let ${\cal J}={\cal I}$, ${\cal I}^{+}$, or ${\cal
I}^{\alpha }$. Then the set of continuous linear operators
in the operator topology of ${\cal J}$ is a complete 
topological $K$-algebra whose operations are given by:
\begin{eqnarray*}
(\theta \phi )p(x)&=&\theta (\phi p(x))\\*
(a\theta )p(x)&=&a(\theta p(x))\\*
(\theta +\phi )p(x)&=&(\theta p(x))+(\phi p(x)).
\end{eqnarray*}
\end{prop}

\proof{Let $(\theta _{n})_{n\geq 0}$ and $(\phi
_{n})_{n\geq 0}$ be convergent
 sequences of continuous linear
operators on ${\cal J}$. Thus, for any $p(x)\in {\cal J}$,
the sequences 
$(\theta _{n}p(x))_{n\geq 0}\mbox{ and }(\phi
_{n}p(x))_{n\geq 0}$
are Cauchy. In particular, 
$(\theta _{n}\phi _{k}p(x))_{n\geq 0}$
 is Cauchy for any
$k\geq 0$,  and since
$\theta _{k}$ is continuous, 
$(\theta _{k}\phi
_{n}p(x))_{n\geq 0}$
 is also Cauchy for any $k\geq 0$.
Hence, 
$(\theta _{n}\phi _{n}p(x))_{n\geq 0}$ is Cauchy, and
$(\theta_{n}\phi _{n})_{n\geq 0}$ converges. 

$((\theta _{n}+\phi _{n})p(x))_{n\geq 0}$ converges since
${\cal J}$ is a topological space. Hence, $(\theta
_{n}+\phi _{n})_{n\geq 0}$ converges. 

The ring is complete since $\theta p(x)=
\lim_{n\rightarrow +\infty } \theta _{n}p(x)$ is the limit
of a Cauchy sequence.}

We shall write infinite series $ \sum _{k\geq d}\theta _{k} $
of operators, which are understood to denote the limits of
their partial sums.

We list below some notable operators and comment briefly on them:

\begin{ex} \bold{Derivative} The derivative $\D^{a}$ is defined via
Definition~\ref{derDef} on all of ${\cal I}$ when $a$ is a 
nonnegative integer. Otherwise, it is defined  on only the positive logarithmic
algebra ${\cal I}^{+}$ via 
\begin{description}
\item[Discrete] Proposition~\ref{invert}.
\item[Continuous] Definition~\ref{fracDer}. Note that the map $a\mapsto \D
^{a}$ is continuous in the operator topology. 
\end{description}
$$ \D ^{a}\lambda_{b}^{\alpha }(x)= \frac{\rn{b}!}{\rn{b-a}!}
\lambda_{b-a}^{\alpha}(x).$$ 
\end{ex}

\begin{ex}\bold{Shift Operator} For all complex
numbers $z$,  the {\em shift operator}, $E^{z}: 
{\cal I} \rightarrow {\cal I}$ is given by the sum
$$ E^{z}=\sum_{n\geq 0}z^{n}\D ^{n}/n!. $$
We see that this is a convergent sum, so $E^{z}$ is
well defined. In fact, we show that $E^{z}$ is a field
isomorphism. 
Note that this is the first use we make of the fact that the field of complex
numbers {\bf C} has characteristic zero.
\end{ex}

\begin{ex}\bold{Elementary {\bf D}-Invariant Operator} 
For any pair of vectors $\alpha $ and $\beta$
the {\em elementary $\D $-invariant operator} from $\alpha $ to
$\beta$---denoted $E_{\beta \alpha }$---is the linear map on the
logarithmic algebra defined for all $a$ and $\gamma $ by
$$ E_{\beta \alpha }\lambda _{a}^{\gamma }(x) =
 \left\{\begin{array}{ll}
\lambda _{a}^{\beta}(x)&\mbox{if $\alpha =\gamma $, and}\\[.1in]
0&\mbox{if $\alpha \neq \gamma $}
\end{array}\right. $$
where $a$ is not a negative integer if $\alpha = (0)$. In the discrete case,
the elementary $\D -$invariant operators are continuous.

Note that  $\proj_{\alpha }=E_{\alpha \alpha }$ is the projection map
${\cal I} \rightarrow {\cal I}^{\alpha}$. In other words, 
$$ \proj_{\beta}(\lambda_{a}^{\alpha}(x)) =
\left\{\begin{array}{ll}\lambda _{a}^{\alpha}(x)&\mbox{if
$\alpha =\beta$}\\[.1in]
0&\mbox{if $\alpha \neq \beta$} \end{array} \right.$$ 
These projections commute with $\D^{a}$. Note however that not
all continuous, linear projections
which commute with $\D ^{a}$
are  expressible in terms of these projections.
\end{ex}

\begin{description}
\item[Continuous] A subspace ${\cal J}$ of ${\cal I}^{+}$ is said to be
\emp{$\D $-invariant} if it is invariant under the fractional
derivative $\D ^{a}$ for all real numbers $a$.  An operator $\theta $ on a
$\D $-invariant subspace ${\cal J}$ is said to be
{\em $\D $-invariant} when $\D ^{a}\theta  =\theta \D ^{a}$ for
all real $a$. 
\item[Discrete] A subspace ${\cal J}$ of ${\cal I}$ is said to be
{\em $\D $-invariant} or \index{Shift-Invariant}{\em shift-invariant} if it is
invariant under the derivative $\D $, or equivalently if it is invariant under
$E^{z}$ for all complex numbers $z$.  An operator $\theta $ on a
$\D $-invariant subspace ${\cal J}$ is said to be
{\em $\D $-invariant} or {\em shift-invariant} when
$\D\theta  =\theta \D$, or equivalently if $\D E^{z}=E^{z}\D$ for all complex
numbers $z$. 
\end{description}

For example, all the operators mentioned above are $\D $-invariant.

A continuous linear operator $\theta $
 on ${\cal I}$ or ${\cal I}^{+}$ is said to be a
{\em regular}\index{Regular} operator if it  commutes 
with every elementary $\D $-invariant operator $E_{\alpha
\beta}$ (except possibly when $\beta=0$),
that is, such that 
$$ \theta E_{\alpha \beta}=E_{\alpha \beta}\theta  $$ 
for $\beta \neq 0$.
For example, $\D ^{a}$ and $E^{z}$ are regular operators.;

A regular $\D $-invariant operator on ${\cal I}^{+}$ is called a {\em Artinian
operator}\index{Artinian Operator}, and one on 
all of ${\cal I}$ is called a  
{\em differential operator}\index{Differential Operator}. The set of Artinian
operators is denoted by $ \Lambda^{+} $, and the set of differential
operators by $\Lambda$. Beware that this notation is inconsistent with the
notation for the corresponding concept in \cite{LR}. The notation used here was
chosen because it is more logical.

Clearly every differential operator restricts to an Artinian
operator. Thus, $\Lambda^{+}\subseteq \Lambda$

\begin{prop}
\begin{enumerate}
\item The set of Artinian operators $\Lambda^{+} $ is a complete
topological ring in the operator topology of ${\cal I}^{+}$. 
\item The set of differential operators $\Lambda$
 is a complete topological ring in the operator topology of
${\cal I}$. 
\end{enumerate}
\end{prop}

\proof{$\Lambda^{+} $ and $\Lambda$ clearly are
$K$-algebras, so it  suffices to show that the limit of
any Cauchy sequence of Artinian (resp.\ differential)
operators is again an Artinian (resp.\ differential) operator.

Let $(\theta _{n})_{n\geq 0}$ be such a Cauchy sequence, and
let $\theta $ be its limit. Now,
\begin{eqnarray*}
\theta E^{z}p(x)
&=&\lim_{n\rightarrow +\infty }\theta _{n}E^{z}p(x)\\*
&=&\lim_{n\rightarrow +\infty }E^{z}\theta _{n}p(x)\\*
&=&E^{z}\lim_{n\rightarrow +\infty }\theta _{n}p(x)
\end{eqnarray*}
since $E^{z}$ is a continuous operator. This in turn equals
$E^{z}\theta p(x)$, so $\theta $ is $\D $-invariant.
{\em Mutatis mutandis}, we have $E_{st}\theta =\theta
E_{st}$, so $\theta $ is regular.}

Our objective is to obtain structural characterizations of
Artinian and differential operators. 

\ssect{Taylor's Formula}
 We shall derive analogs of Taylor's
formula in the Logarithmic algebra ${\cal I}$.
We begin by giving the following alternate definition of the
shift operator:

\begin{prop} \label{adthm}
For all complex numbers $z$,
$E^{z}$ is a well defined continuous field isomorphism of
${\cal I}$ which fixes all constants..
\end{prop}

\proo{Need only check that
$E^{z}(p(x)q(x))=(E^{z}p(x))(E^{z}q(x))$ for all $p(x),q(x)
\in{\cal I}$. 
\begin{eqnarray*} 
E^{z}(p(x)q(x)) &=&\sum _{k\geq
0}\frac{z^{k}}{k!}\D ^{k}(p(x)q(x))\\* 
&=&\sum _{k\geq 0}\frac{z^{k}}{k!} \sum _{n+m=k}{k \choose
n} \left(\D^{n}p(x)\right) \left(\D ^{m}q(x)\right)\\ 
&=&\sum_{n,m\geq 0}\frac{z^{n+m}}{n!m!}\left(\D ^{n}p(x)\right)
\left(\D ^{m}q(x)\right)\\*
&=&\left(E^{z}p(x)\right)
\left(E^{z}q(x)\right).\Box 
\end{eqnarray*}}

\digress{Proposition~\ref{adthm} shows that $E^{z}$ satisfies the conditions of
the characterization of Artinian composition in \cite{ch1}.
Any continuous field isomorphism of ${\cal I}$
which fixes all constants is
determined by its values on the iterated logarithms $\ell _{k}$.
Thus, $E_{z}$ is the 0-composition associated with the the substitution of
$x+z$ for $x$, and $\log x+\sum_{j>0}(-1)^{j+1}z^{j}/jx^{j}$, and so on. 
In general,
\begin{eqnarray*}
E^{z}\ell _{k}&=& \sum_{n\geq 0}\frac{z^{n}}{n!}\D ^{n}\ell _{k}\\*
	&=& \sum_{n\geq
0}\frac{z^{n}}{\rn{-n}!n!}\lambda_{-n}^{(0,\ldots
,0,1,0,\ldots )}(x)\\
	&=&  \ell _{k} - \sum_{n\geq 0} 
\sum_{\scriptstyle \ell(\rho )=k \atop \scriptstyle \rho _{k}=1}
\frac{(-z)^{n}s(-n,\rho _{1})}{\rn{n}!n} 
\left[ \prod_{j=1}^{k-2}
e_{\rho _{j}-\rho _{j+1}} (-1,\ldots ,-\rho _{j}+1)
\right] \ell ^{(-n),\rho ^{*}}\\
	&=& \ell _{k} + \sum_{\scriptstyle \ell (\rho )=k\atop
\scriptstyle \rho_{k}=1 } (-1)^{\rho _{1}+|\rho |} 
\left[ \prod_{j=1}^{k-2}
e_{\rho _{j}-\rho _{j+1}}(-1,\ldots ,-\rho _{j}+1) \right] 
\sum_{n\geq 0}\frac{z^{n}(n-1)!s(-n,\rho _{1})}{n} \ell ^{(-n),\rho ^{*}}\\*
	&=& \ell _{k} + \sum_{\scriptstyle \ell (\rho )=k\atop
\scriptstyle \rho_{k}=1 } (-1)^{\rho _{k-1}+|\rho |} 
\left[ \prod_{j=1}^{k-2}
e_{\rho _{j}-\rho _{j+1}}(1,\ldots ,\rho _{j}-1) \right] 
\sum_{n\geq 0}\frac{z^{n}(n-1)!s(-n,\rho _{1})}{n} \ell ^{(-n),\rho ^{*}}
\end{eqnarray*}
where $\rho ^{*}=(\rho _{1},\ldots ,\rho _{k-1})$.}

Observe that $E^{z_{1}}E^{z_{2}}=e^{z_{1}\smallD}
e^{z_{2}\SB{D}} = e^{(z_{1}+z_{2}) \SB{D}}
=E^{z_{1}+z_{2}}.$   

\discrete{By applying the shift operator to the  harmonic logarithms
of order $(1)$ and  degree $n$ a nonnegative integer, $\lambda
_{n}^{(1)}(x)$, we obtain the following 
identity:
\begin{eqnarray*}
\lefteqn{(1+a)^{n}\left(\log (1+a)-1-\frac{1}{2}-\cdots
-\frac{1}{n}\right)}\\*
 &= &\left[(x+a)^{n}\left(\log (x+a)-
1-\frac{1}{2}-\cdots -\frac{1}{n} \right)\right]_{_{x=1}}\\
&=&\left[\sum
_{i=0}^{n}\romcoeff{n}{i}a^{n-i}x^{i}\left(\log x-
1-\frac{1}{2}-\cdots -\frac{1}{i}\right) 
+\sum _{i>n}\romcoeff{n}{i}a^{i}x^{n-i}\right]_{_{x=1}}\\ *
&=&-\sum
_{i=0}^{n}\romcoeff{n}{i}a^{n-i} 
\left(1+\frac{1}{2}+\cdots +\frac{1}{i} \right)
+ \sum _{i>n}\romcoeff{n}{i}a^{i},
\end{eqnarray*}
 and therefore:
 \begin{eqnarray*}
(1+a)^{n}\log (1+a) 
 &= &-\left((1+a)^{n}-1\right)\frac{s(-n,1)}{\rn{-n}!}
+\sum_{i=0}^{n-1}{n \choose i }a^{n-i}  \frac{s(-i,1)}{\rn{-i}!}
+\sum _{i>n}\romcoeff{n}{i}a^{i} \\*
&=&\left((1+a)^{n}-1\right)
\left(1+\frac{1}{2}+\cdots+\frac{1}{n}\right)
-na\left(1+\frac{1}{2}+\cdots +\frac{1}{n-1}\right) \phil 
-{n\choose 2}a^{2}\left(1+\frac{1}{2}+\cdots
+\frac{1}{n-2}\right) -\cdots
-\frac{3}{2}{n\choose 2 }a^{n-2}-na^{n-1} \phil
+a^{n+1}\romcoeff{n}{-1}+a^{n+2}\romcoeff{n}{-2}+\cdots .
\end{eqnarray*}
Contrast with \cite[Corollary 5.6]{GBC}.}

\ssect{Real Analysis}\index{Real Analysis}\label{RA}
Proposition~\ref{adthm} can be summarized as stating that
 $$E^{z}p(x)=e^{z\D }p(x)=p(x+z) $$
 for all logarithmic series
$p(x)\in {\cal I}$. Note that this identity is tautological,
since we cannot ``evaluate'' the variable $x$. However, in
the case of complex numbers, we can, and we have
\begin{prop}\label{Roman1}
All formal power series of logarithmic type $p(x)\in {\cal I}$ represent
asymptotic expansions of a complex valued function $\tilde{p}(x)$ in a
neighborhood of infinity. Moreover, for all real numbers $c$, $E^{c}p(x)$
represents an asymptotic expansion of $\tilde{p}(x+c).$

Conversely, any asymptotic expansion of a function in a neighborhood of
infinity relative to the ladder of comparison given by
$\lambda_{a}^{\alpha}(x)$  is a formal power series of logarithmic type.
\end{prop}

\proo{The asymptotic expansions follow from Proposition~\ref{asyProp}. We
conclude by observing 
\begin{eqnarray*}
E^{c }p(x)|_{x=z} &=& \left. e^{c \smallD}\tilde{p}(x)\right|_{x=z}\\*
&=& \left.\sum_{k\geq 0}\frac{c ^{k}}{k!}\D^{k}\tilde{p}(x)\right|_{x=z}\\*
&=& \left. \sum _{k\geq 0}\frac{x^{k}}{k!}\left[\D ^{k}p(x)
\right]_{x=z}\right|_{x=c}.
\end{eqnarray*}}

\ssect{Characterization of Various Classes of Operators}
\sssect{Artinian Operators}
Every Artinian operator $\theta $ maps ${\cal I}^{\alpha}$
into itself for every nonzero vector $\alpha $. We
denote by $\theta $ its restriction to ${\cal I}^{\alpha}$, by
an abuse of notation, and we say that $\theta $ is an Artinian
operator of ${\cal I}^{\alpha}$ into itself. 

As planned, we can now characterize the algebra of Artinian operators and its
topology. It is isomorphic to an algebra (mentioned in \cite{ch1}) we denote
${\bf C}<x>$; it is in a sense the dual of the Noetherian algebra of 
Definition~\ref{LoAl}. 
\begin{description}
\item[Discrete] In the discrete case, $K<x>$ is the algebra of Laurent series
in the variable $x$ with complex coefficients.
\item[Continuous] In the continuous case, $K<x>$ is the algebra of Artinian
series in the variable $x$ with complex coefficients.
$$ K<x> = \left\{ 
\parbox{4in}{$\sum _{a\in {\bf R}}c_{a}x^{a}$: $c_{a}\in K${ and for all }
$a\in {\bf R}${, there exists finitely many }$b\geq 0${ with }
$c_{a}\neq 0.$}\right\}. $$
Its operations and topology are similar to that of the Noetherian algebra.
\end{description}

\begin{thm} \label{sio1}\index{Artinian Algebra}
 The algebra of Artinian operators $\Lambda ^{+}$ is
naturally isomorphic to the Artinian algebra ${{\bf C}}<x>$ as a
topological algebra
$\Lambda^{+}={\bf C}(\D )^{R_{0}}.$
\end{thm}

\proof{{\bf (${\bf C}(\D)^{R_{0}}\subseteq \Lambda^{+}
$)} Note that $\D ^{a}\D ^{b}=\D ^{b}\D ^{a}$, and $\D
^{a}E_{\alpha \beta}=E_{\alpha \beta}\D ^{a}$.  

{\bf ($\Lambda^{+} \subseteq {\bf C}(\D)^{R_{0}}$)}
Let $\theta\in \Lambda^{+}$ be an Artinian operator. 
By regularity, $\theta $
is determined by its
action on the harmonic logarithms of order $\alpha $ for any
particular $\alpha \neq  (0)$.
$$ \theta \lambda _{a}^{\alpha }(x)= \sum _{b}
e_{ab}\lambda_{b}^{\alpha }(x).$$
Notice, that for any $a$ and $b$, $e_{ab'}\neq 0$ only for
finitely many $b'\geq b$.
Next,
\begin{eqnarray*}
\theta \D^{a}\lambda _{b}^{\alpha}(x)
&=&\theta \frac{\rn{b}!}{\rn{b-a}!}\lambda _{b-a}^{\alpha }(x)\\*
&=&\frac{\rn{b}!}{\rn{b-a}!}\sum_{c}
e_{a-b,c}\lambda _{c}^{\alpha }(x)\\
\D ^{a}\theta \lambda_{b}^{\alpha}(x) &=& \D
^{a}\sum_{c}e_{bc}\lambda_{c}^{\alpha}(x)\\*
&=&\sum_{c}e_{bc}\frac{\rn{c}!}{\rn{c-a}!}\lambda_{c-a}^{\alpha}(x).
\end{eqnarray*}
Thus,
$\romcoeff{b}{a}e_{a-b,c}=
\romcoeff{c+a}{c}e_{b,c+b}$
In particular, setting $a=b$, 
$$e_{0,c}=\romcoeff{c}{b}e_{b,c+b}.$$ 
Hence, $\theta $ is determined by the 
$e_{0,c}$, and therefore equals  the Artinian series
$\sum _{c}e_{c,0}\D ^{c}/\rn{c}!.$

{\bf (Continuity)} It remains only to show that Artinian series in the
derivative are continuous in the operator topology (\S\ref{OpTop}),
and that Artinian operators  are continuous in the Artinian topology. However,
$c_{n}\D ^{a_{n}}$ converges whenever $a_{n}$ and 
$c_{n}$ converge, and $c_{n}\D ^{a_{n}}$ converges to zero whenever $a_{n}$
increases without bound. Moreover, all Cauchy sequences of operators are linear
combinations of the above.}
 
If $f(\D )$ is an Artinian operator which is
represented by a delta series (that is, a series of degree one), then $f(\D )$
is called a {\em delta operator}.\index{Delta Operator}
\begin{description}
\item[Discrete] Recall the
important fact that if $f(\D)$ is a delta operator, then the sequence of
powers $(f(\D )^{n})_{n\in \SB{Z}}$ is a
pseudobasis for the field of Artinian operators
$\Lambda^{+} $. That is, for every Artinian operator $g(\D )\in \Lambda ^{+}$
of degree $d$
there is a sequence $(c_{k})_{k\leq d}$ indexed by integers $k\leq d$
such that $g(\D)= \sum_{k\leq d}c_{k}f(\D )^{k}$.
This sequence can actually be
calculated via Theorem~\ref{GET}. The notation $f(\D )^{a;n}$ comes from the
continuous case (see below); in the context of the discrete theory, it should
be read simply as $f(\D )^{a}$.
\item[Continuous] Recall the important fact from \cite{ch1} that if $f(\D
)=\sum _{a\in {\bf R}}c_{a}x^{a}$ is
an Artinian operator of degree $d\neq 0$ (for example, a delta operator), and
$n$ is an integer, then the set $\{f(\D
)^{a;n}:a\in {\bf R} \}$ is a pseudobasis for the field of Artinian operators
$\Lambda^{+}$ where 
$$ g(x)^{t;n}=m^{t}e^{it(\theta +2n\pi)}c_{d}^{t;n} x^{dt} \sum _{M} {t\choose
M} \left(\prod _{a\in M} \frac{c_{a+d}}{c_{d}}x^{a} \right) $$
and $m=|c_{d}|$ and $\theta = \mbox{arg}(c_{d})$ are the modulus and argument
of the leading coefficient of $f(\D )$.

That is, for every Artinian operator $g(\D )\in \Lambda ^{+}$ there is
a sequence $(c_{a})_{a\in \SB{R}}$ such that $\sum_{a\in
\SB{R}}c_{a}f(\D )^{a;n}$ converges to $g(\D )$. When $f(\D )$ is a delta
operator, this sequence can  be calculated via Theorem~\ref{GET} or
the methods of \cite{ch2}.
\end{description}
Similarly, the series of nonnegative integers powers of any delta operator
$f(\D )$ is a pseudobasis for $\Lambda$.

\sssect{Differential Operators}
In analogy with the preceding result for Artinian operators,
we obtain the following structure theorem for differential
operators: 

\begin{cor}\label{DiffChar}
The topological algebra of differential
operators is naturally isomorphic to the topological algebra of formal power
series in the derivative as a topological algebra:
$\Lambda = {\bf C}[[\D ]]$. 
\end{cor}

\proof{Such series are the only members of
$\Lambda^{+}$ which are well defined on ${\cal I}^{(0)}.$}

As opposed to Artinian operators which are invertible if nonzero, a
differential operator is invertible in $\Lambda$ if and only if it is of degree
$0$.  

\sssapp{$\D $-invariant Operators}

The rings of continuous, linear, $\D $-invariant operators (which are not
necessarily regular) over ${\cal I}$ and ${\cal I}^{+}$ are structured as
follows: 

\begin{prop}\label{sioprop}
Let ${\cal R}$ and ${\cal R}_{0}$ be the rings of $\D $-invariant linear
operators on  ${\cal I}^{+}$ and  ${\cal I}$
respectively which are continuous on each ${\cal
I}^{\alpha }$. Then:
\begin{description}
\item[Discrete]
\begin{enumerate}
\item ${\cal R}$ is the closure ${\cal D}$, in the operator topology, of the 
span of the operators $\D ^{n}E_{\alpha \beta}$ 
where $\alpha $ and $\beta$ are nonzero vectors with finite support of
integers, and $n$ is an integer.
\item  ${\cal R}_{0}$ is the
closure ${\cal D}_{0}$, in the operator topology, of the span of 
the operators $\D ^{n}E_{\alpha \beta}$ where
\begin{enumerate}
\item $\alpha $ and $\beta$ are  vectors with finite support of integers,
\item   $n$ is an integers,
\item  $\alpha \neq (0)$ unless $\beta = (0)$, and
\item  $\beta\neq (0)$ unless $n$ is a nonnegative integer. 
\end{enumerate}
\end{enumerate} 
\item[Continuous]
\begin{enumerate}
\item ${\cal R}$ is the closure ${\cal D}$, in the operator topology, of the 
span of the operators $\D ^{a}E_{\alpha \beta}$ 
where $\alpha $ and $\beta$ are nonzero vectors with finite support of
real numbers, and $a$ is a real number.
\item  ${\cal R}_{0}$ is the
closure ${\cal D}_{0}$, in the operator topology, of the span of 
the operators $\D ^{a}E_{\alpha \beta}$ 
\begin{enumerate}
\item $\alpha $ and $\beta$ are  vectors with finite support of real numbers,
\item   $a$ is a real number,
\item  $\alpha \neq (0)$ unless $\beta = (0)$, and
\item  $\beta\neq (0)$ unless $a$ is a nonnegative integer. 
\end{enumerate}
\end{enumerate}
\end{description}
\end{prop}

\proof{{\bf (${\cal D}\subseteq {\cal R}$ and ${\cal
D}_{0}\subseteq {\cal R}_{0}$)}
Elementary
$\D $-invariant operators commute with the derivative. Thus,
they 
commute with all Artinian and differential operators. Hence,
they are in fact 
$\D $-invariant. 

Observe that $E_{\alpha \beta}\D^{a}$ is continuous,
linear and $\D $-invariant for $\alpha ,\beta\neq  (0)$.
We conclude that every operator in ${\cal
D}$ and ${\cal D}_{0}$
 is continuous, linear, and $\D $-invariant.

{\bf (${\cal R}\subseteq {\cal D}$)} Let
$\theta $ be a continuous, linear, $\D $-invariant operator.
For each pair of vectors $\alpha ,\beta\neq  (0)$, define $\theta
_{\alpha \beta}=\proj_{\alpha }\theta \proj_{\beta}$.
Obviously, $\theta =\sum 
_{\alpha ,\beta\neq  (0)}\theta _{\alpha \beta}$. It 
suffices to show that for all 
nonzero vectors $\alpha $ and $\beta$, there is an Artinian operator
$f_{\alpha \beta}(\D )\in \Lambda^{+} $ such that $\theta
_{\alpha \beta}=f(\D )E_{\beta\alpha }$.

However, $E_{\beta\alpha }\theta _{\alpha \beta}$ is a
continuous, linear, 
$\D $-invariant operator on ${\cal I}^{\alpha}$, so
$E_{\beta\alpha}\theta _{\alpha \beta} =f(\D )\proj
_{\beta}$ for some Artinian operator 
$f(\D )$. Hence, $\theta _{\beta\alpha }=E_{\beta\alpha
}f(\D )E_{\alpha \alpha }=f(\D)E_{\beta\alpha }$ as desired.

{\bf (${\cal R}_{0}\subseteq {\cal D}_{0}$)} Similarly, it
 suffices to show 
that $\theta_{\alpha ,(0)}$ (as defined  above) is equal
to zero for $\alpha \neq  (0)$. Assume not towards contradiction.
By the reasoning above,
$\theta_{\alpha ,(0)}=f(\D )E_{\alpha ,(0)} $ for some nonzero differential
operator $f(\D )$. Let $g(\D )$ be
the inverse of $f(\D )$. Since $g(\D )$ is a $\D $-invariant
operator, the product $g(\D )\theta $ is also
$\D $-invariant. Hence, without loss of generality,
$\theta_{\alpha ,(0)}=E_{\alpha ,(0)}$. We calculate that
$E_{\alpha ,(0)}\D 1=E_{\alpha ,(0)}0$. However, we also know that
$\D E_{\alpha ,(0)}1=E\ell ^{(0),\alpha }\neq 0$.
Contradiction.}

We omit a discussion of non-$\D $-invariant continuous, linear, regular
operators and their expansions in terms of $\D$ and $\sigma $ (to be defined
later), since this would be an unnecessary digression. 

However, there are other interested related questions left unsolved.
\begin{prob}
Is there a simple characterization of non-linear, continuous, shift-invar\-iant
operators on polynomials? Or of the logarithmic algebra? For example,
$p(x)\mapsto p(x)^{2}$.
\end{prob}

\sssapp{Differential Equations}\index{Differential
Equations}\label{diffEqSec} 

We next make some general remarks about solutions of
differential equations of infinite order (and thus, in
particular, of difference equations).
We have seen (Theorem~\ref{sio1}) that the algebra of
Artinian operators is isomorphic to the 
Artinian algebra in the ``variable'' $\D $ as a topological algebra over the
complex numbers. In particular, every Artinian operator is invertible.
Hence, every differential equation of the form
\begin{equation}\label{diff-eq}
 f(\D )p(x)=q(x) 
\end{equation}
where $f(\D )\in \Lambda^{+} $ is an Artinian
operator and $q(x)\in {\cal I}^{+}$ has
a unique solution $p(x)$ in ${\cal I}^{+}$.
For example, 
$q(x)$ may be any rational function of $x$ whose numerator
is of smaller degree than the denominator. 

If $q(x)\in{\cal I}$ (for example, if $q(x)$ is an arbitrary
rational function) the solution may not be unique. We can
nevertheless define the {\em natural solution}\index{Natural Solution of
Differential Equations} of
\eref{diff-eq} as follows. Let $\pi $ be any bijection between the
equipotent sets $A$ and $A-\{(0) \}$ where $A$ is the set of all vectors 
$\alpha $ with finite support of reals (resp. integers).
Define $P_{\pi }$ to be the linear map (continuous on each ${\cal I}^{\alpha
}$) defined by
$$P_{\pi } = \sum _{\alpha } E_{\pi \alpha ,\alpha }$$
or equivalently
$$ P_{\pi }\lambda_{a}^{\alpha}(x) =\lambda_{a}^{\pi \alpha}(x)$$
For example, we might have
$$P_{\pi }=\sum _{n\geq 0} E_{(n)(n+1)} + \sum _{\alpha \neq (n)}
E_{\alpha \alpha }$$
where the second sum in each equation is over vectors which
do not consist of a single nonnegative integer.
Thus, $P_{\pi }:{\cal I}\rightarrow {\cal I}^{+} $ and its inverse $P_{\pi
}^{-1}:{\cal I}^{+}\rightarrow {\cal I}$ is given by $P_{\pi
}^{-1}\lambda_{a}^{\alpha}(x)=\lambda_{a}^{\pi ^{-1}\alpha}(x)$.
Now, let $s(x)=P_{\pi }q(x)$, and
consider the differential equation
$$ f(\D )r(x)=s(x). $$
Since $s(x)\in {\cal I}^{+}$, this differential equation
has a unique solution $r(x)\in {\cal I}^{+}$. Now, set
$p(x)=P_{\pi }^{-1}r(x)$ to obtain a solution of \eref{diff-eq}.
The present definition agrees with (and is simpler than) all
other definitions given of a natural solution over the
complex numbers. Moreover, since $p(x)$ does not depend on the choice of $\pi
$, the solution has been chosen naturally.

A notable example is the difference equation
$$ \Delta p(x)=1/x. $$
By the above remarks, it has a unique solution in ${\cal
I}^{+}$, which turns out to be the $\psi$-function $\psi
(x)$, the logarithmic derivative of the gamma function.
Thus, the theory of the $\psi$-function can be developed
purely formally. (See \cref{psisec}.)

So far, we have only consider linear differential equations. This leads us to
the following open problem.

\begin{prob}
In general, what differential equations have solutions? And when do they have
canonical solutions?
\end{prob}

The logarithmic algebra is not a {\em differentially closed
field}\index{Differentially Closed Fields} since $\D p(x)= p(x)$ only has one
solution---$p(x)=0$. However,  we can redefine the degree of a differential
operator so that we are consistent with \cite{ch1} by insisting
that its degree is the lowest (rather than highest) exponent of $\D $ with a
nonzero coefficient. Under these circumstances, $\D p(x)=p(x)$ is no longer a
counterexample.

\thmplus{C}{\begin{prob}
Under this definition of degree, is the logarithmic algebra differentially
closed? And if so, can we modify the logarithmic algebra so that it is
differentially closed under the usual definition?
\end{prob}}

Recall that differentially closed fields are known to exist in every character;
however, no example of a differentially closed field of character zero has been
found. 

\ssect{Augmentation}

Although it is impossible to evaluate at zero any expression involving
logarithms, the following definition  nevertheless serves as the logarithmic
analog of evaluation at zero.

\begin{defn}\bold{Augmentation}\label{augment} 
For each vector $\alpha $, we define the {\em augmentation}\index{Augmentation}
of  order $\alpha $ to be the linear functional $\action{\alpha } {}$
(continuous on each ${\cal I}^{\alpha }$) from the logarithmic
algebra to the complex numbers such that 
$$\action{\alpha }{\lambda _{a}^{\beta}(x)} = \delta
_{\alpha \beta}\delta _{a,0},$$
or equivalently using the notation introduced after Corollary~\ref{Lsum},
$$\action{\alpha }{p(x)}=[\lambda_{0}^{\alpha}(x)]p(x).$$
\end{defn}

\digress{We {\em digress} to indicate the algebraic significance of
augmentation in the discrete case. First, note that in this case
$\action{\alpha }{}$ is continuous.

When $\action{0}{}$ is restricted to ${\cal I}^{(0)}$, the
augmentation reduces the
evaluation of a polynomial at $x=0$.  That is,
$\action{0}{p(x)}=p(0)$ for $p(x)\in {\cal I}^{(0)}$.
An augmentation of nonzero order can be
viewed as a 
generalization of evaluation at $x=0$; it is closely related
to the residue of complex variable theory. (Recall that the
residue of a Laurent series is its coefficient of $x^{-1}$.)
For instance, for $p(x)\in {\cal I}^{(1)}$,
$$ \action{1} {p(x)}=\Res{\D p(x)}.$$
In fact, for $\alpha$ a vector of $j$ nonnegative integers
and $p(x)\in {\cal I}^{\alpha }$, 
\begin{equation}\label{xd}
 \action{\alpha } {p(x)}= \mbox{Res} \left(\D \left(\prod_{k=1}^{j}
\frac{\ell ^{(0),\delta _{k}}}{\alpha _{k}!} 
\right)p(x)\right).
\end{equation}
Note that \eref{xd} also holds  for all $p(x)\in
\bigoplus _{\beta\leq \alpha } {\cal I}^{\beta}$.}

We derive formulas relating the  augmentation to the derivative. 

\begin{prop}\label{simplest}\label{action}
\begin{enumerate}
\item For all $a$ and $b$, and all $\alpha ,\beta\neq  (0)$
$$\action{\beta} {\D^{b} \lambda_{a}^{\alpha}(x)}=
\rn{a}!\delta_{ab}\delta_{\alpha \beta}.$$
\item For all $a$ and all nonnegative integers $n$,
$$\action{(0)} {\D ^{n}\lambda
_{a}^{(0)}(x)}= n!\delta_{an}.$$
\item If $f(\D )$ is an Artinian
operator with coefficients $[\D ^{a}]f(\D )=c_{a},$
 and $p(x)$ is a formal power series of logarithmic type with coefficients
$$b_{a}^{\alpha }=[\lambda _{a}^{\alpha }(x)]p(x),$$
then the augmentation  $\action{\alpha }
{f(\D )p(x)}$ is given by the finite sum
$$\action{\alpha } {f(\D )p(x)}=\sum _{a}\rn{a}! c_{a}b_{a}^{\alpha}.$$
\item If $f(\D )$ is a differential
operator with coefficients $c_{k}=[\D ^{k}]f(\D ),$
 and $p(x)$ is a formal power series of logarithmic type with coefficients
$$b_{a}^{\alpha }=[\lambda _{a}^{\alpha }(x)]p(x),$$
then the augmentation $\action{\alpha } {f(\D )p(x)}$ is
given by the finite sum
$$\action{\alpha } {f(\D )p(x)}=\sum _{k\geq 0}k!c_{k}b_{k}^{\alpha}.\Box $$
\end{enumerate}
\end{prop}

Some special cases are of interest.
The augmentation of the derivative of a formal power series of
logarithmic type can
be described by
$\action{\alpha }{\D ^{a}p(x) }=\rn{a}!b_{a}^{\alpha}$ where
$p(x)$ is as  in Proposition~\ref{action}.
Similarly, the
augmentation of the action of a differential operator on a harmonic
logarithm is $\action{\alpha }{f(\D )\lambda _{a}^{\alpha}(x)}=\rn{a}!c_{n}$
where $f(\D)$ is as  in Proposition~\ref{action}.

The augmentation leads us to a version of Taylor's formula
for formal power series of logarithmic type:
\begin{thm}\bold{Taylor's Formula}\label{ltf}
Let $p(x)\in {\cal I} ^{+ }$.
Then we have the following convergent expansion in ${\cal
I}^{\alpha}$, 
$$ p(x) = \sum_{\alpha \neq (0)}
\sum_{a}\frac{\action{\alpha } {\D ^{a}p(x)}}{\rn{a}!} \lambda
_{a}^{\alpha}(x).$$ 
\end{thm}

\proof{By linearity and continuity, 
it suffices to consider the case
$p(x)=\lambda _{a}^{\alpha}(x)$. This special case was  treated
above (Proposition~\ref{simplest}).}

Note that the series of the right hand side is convergent. 

\begin{description}
\item[Discrete] Observe that in the discrete case, any formal power series of
logarithmic type is determined by the augmentations of its
derivative. 
\item[Continuous] Nevertheless,  in the continuous case,  members of ${\cal
I}^{(0)}$ can not be recovered from the augmentations of
their derivatives. For example, $\action{(0)} {f(\D )\sqrt{x}}=0$ for all
differential operators $f(\D )$.
\end{description}

On the other hand, in either case, an Artinian operator can be recovered from
the augmentations of its 
action on the harmonic logarithms of any particular order
$ \alpha \neq  (0)$ and a differential operator can be recovered from the
augmentations of the harmonic logarithms of any particular order as is
demonstrated by the following theorem.

\begin{thm} \bold{Expansion Theorem}\label{exp1}
\begin{enumerate}
\item Let $f(\D )$ be an Artinian operator, and let $\alpha
\neq  (0)$, then we have the following convergent expansion
$$f(\D )=\sum _{a}\frac{\action{a}{f(\D)\lambda
_{a}^{\alpha}(x)}}{\rn{a}!}\D ^{a}.$$ 
\item Similarly, if $f(\D )$ is a differential operator,
and $\alpha $ is a vector, then we have the following
convergent series
$$f(\D )=\sum _{n\geq 0} \frac{\action{\alpha }{f(\D)\lambda _{n}^{\alpha}(x)}}
{\rn{n}!} \D ^{n}.\Box $$
\end{enumerate}
\end{thm}

The following argument is used repeatedly in the next
two sections---often implicitly.

\begin{prop}
\bold{Spanning Argument} \label{spanning}
\begin{enumerate}
\item Let $p(x)\in {\cal I}^{+}$. If $\action{\alpha }{f(\D )p(x)} =
0$ for all vectors $\alpha \neq  (0)$
and all Artinian operators $f(\D )$, then $p(x)=0$. 
\item Let $\alpha \neq  (0)$ be a vector, and let $f(\D )$ 
be an Artinian operator. If
$\action{\alpha }{f(\D )p(x)}=0$ for all $p(x)\in {\cal
I}^{\alpha}$, then $f(\D ) =0.$ 
\item \discrete{Similarly, in the discrete case, for ${\cal I}^{(0)}$,
let $p(x)\in {\cal I}^{(0)}$. If $\action{(0)}{f(\D )p(x)} = 0$,
for all  differential operators 
$f(\D )\in \Lambda$, then $p(x)=0$.}
\item Let $\alpha $ be a vector, and let $f(\D )$ be a differential
operator. If  for all $p(x)\in {\cal I}^{\alpha }$,
$\action{\alpha }{f(\D )p(x)}=0$, then $f(\D ) =0.\Box $ 
\end{enumerate}
\end{prop}

\begin{prob}
Is there a simple formula expressing any monomial $\ell^{\alpha }$ or product
of harmonic logarithms $\lambda_{a}^{\alpha}(x)\lambda_{b}^{\alpha}(x)$ in
terms of harmonic logarithms? That is, is there a simple way to calculate 
the coefficients of $\ell ^{\alpha }$ or
$\lambda_{a}^{\alpha}(x)\lambda_{b}^{\alpha}(x)$ given by Theorem~\ref{ltf}?
\end{prob}

\discrete{Note that the expansion of $\lambda_{n}^{(k)}(x)\lambda_{m}^{(k)}(x)$
is known.}

\sect{Roman Graded Sequences} \label{seqsec}
\ssect{Graded Sequences}\label{appellsect}

This chapter is devoted to the study of the logarithmic analog of sequences of
polynomials of binomial type. However, before one can walk; one must crawl. We
must first define the logarithmic analog of a sequence of polynomials.
Classically, a sequence of polynomials $(p_{n}(x))_{n \geq 0}$ is subject to
the 
requirements that $\deg(p_{n}(x))=n$ for all $n$. Here the requirements are
slightly more complicated, yet the idea remains the same.

\begin{defn}[Graded Sequences of Logarithmic Series]
\label{sfpslt}
The sequence $\{ \seq{p}{a}{\alpha }:$ $a$ real
(resp.\ an integer) and $\alpha $ a vector of reals
(resp.\ integers)$\}$
\figx{Part of a Typical Graded Sequence of Formal Power Series of Logarithmic
Type}{\fbox{\begin{tabular}{ccccccccc}\label{sugFig2}
 & $\vdots$&$\vdots$&$\vdots$&$\vdots$&$\vdots$&$\vdots$&$\vdots$ &\\
$\cdots $ & $p_{-3}^{(-2)}(x)$ & $p_{-2}^{(-2)}(x)$ &
$p_{-1}^{(-2)}(x)$ & $p_{0}^{(-2)}(x)$ & $p_{1}^{(-2)}(x)$ & 
$p_{2}^{(-2)}(x)$ & $p_{3}^{(-2)}(x)$ & $\cdots$ \\
$\cdots $ & 0 & 0 & 0 & $\tilde{p}_{0}(x)$ & $\tilde{p}_{1}(x)$ &
$\tilde{p}_{2}(x)$ & $\tilde{p}_{3}(x)$ & $\cdots $ \\
$\cdots $ & $\tilde{p}_{-3}(x)$ & $\tilde{p}_{-2}(x)$ &
\fbox{$\tilde{p}_{-1}(x)$}  & 
$p_{0}^{(1)}(x)$ & $p_{1}^{(1)}(x)$ & 
$p_{2}^{(1)}(x)$ & $p_{3}^{(1)}(x)$ & $\cdots $ \\
$\cdots $ & $p_{-3}^{(2)}(x)$ & $p_{-2}^{(2)}(x)$ &
$p_{-1}^{(2)}(x)$ & 
$p_{0}^{(2)}(x)$ & $p_{1}^{(2)}(x)$ & 
$p_{2}^{(2)}(x)$ & $p_{3}^{(3)}(x)$ & $\cdots$ \\
 & $\vdots$&$\vdots$&$\vdots$&$\vdots$&$\vdots$&$\vdots$&$\vdots$
\end{tabular}}}
 is called 
a continuous (resp. discrete) {\em graded sequence of formal power series of
logarithmic type} 
if the following six (resp.~five) conditions hold. (Part of a typical such
sequence is illustrated by Figure~\ref{sugFig2}
\begin{enumerate}
\item For all $a$ and $\alpha \neq  (0)$,
$p_{a}^{\alpha}(x)$ is a homogeneous formal power series 
of logarithmic type of order $\alpha $. For example, the series in the $k\th $
row of Figure~\ref{sugFig2} are homogeneous of order $(k)$.
\item For $n$ a negative integer, $p_{n}^{(0)}(x)=0$. That is, the left side of
the middle row of Figure~\ref{sugFig2} is filled with zeroes.
\item For $a$ not a negative integer, $p_{a}^{(0)}(x)$
is a logarithmic series of degree $a$. In other words, in the discrete case
for $n$ a nonnegative integer 
$p_{n}^{(0)}(x)$ is a polynomial of degree $n$. That is, the right side of the
middle row of Figure~\ref{sugFig2} is filled with polynomials in order of
degree. 
\item For $\alpha \neq  (0)$, 
$p_{a}^{\alpha}(x)$ is of degree $a$. For example, the $n\th $ column of
Figure~\ref{sugFig2} consists of series of degree $n$.
\item {\bf (Regularity)}\index{Regularity} For all $a$ and $\beta$, and all
$\alpha \neq (0)$,
$E_{\alpha \beta}p_{a}^{\beta}(x)=p_{a}^{\alpha}(x)$. For example, the operator
$E_{(j)(k)}$ sends the $k\th$ row of Figure~\ref{sugFig2} to the $j\th$ row.
\item \cont{{\bf (Continuity)} In the continuous case only, the map from
the reals to ${\cal I}^{\alpha }$ defined by 
\begin{equation}\label{star63}
a\mapsto \left\{ \begin{array}{ll}
p_{a}^{\alpha }(x)&\mbox{$a$ not a negative integer}\\[0.1in]
p_{a}^{\alpha +(1)}&\mbox{$a$ is a negative integer}
\end{array} \right.
\end{equation}
must be continuous for all $\alpha $.}
\end{enumerate}
When a sequence $p_{a}^{\alpha }(x)$ has the property that
its leading coefficients are all positive real numbers, then it is
said to be {\em standard}\index{Standard Sequence}.

The logarithmic series $p_{-1}^{(1)}(x)$ is called the
{\em residual series}\index{Residual Series} of the graded sequence
$\seq{p}{a}{\alpha}$. (Indicated by a box in Figure~\ref{sugFig2}.)

The {\em principle subsequence}\index{Principle Subsequence} of
$\seq{p}{a}{\alpha }$ is the 
subsequence $\left(\tilde{p}_{a}(x)\right) _{n\in \SB{R}}$ defined by 
$$\tilde{p}_{a}(x) =\left\{\begin{array}{ll}
p_{a}^{(0)}(x)&\mbox{for $a$ not a negative integer, and}\\[0.1in]
p_{a}^{(1)}(x)&\mbox{for $a$ a negative integer.}
\end{array}\right. $$
\end{defn}

For example, 
\begin{prop}
\label{harmseq}
$\seq{\lambda}{a}{\alpha }$ is a standard graded sequence of
formal power series of logarithmic type. Its residual series
is $1/x$, and its principal subsequence is the sequence of
powers of $x$, $\left(x^{n}\right) _{n\in \SB{Z}}$.
\end{prop}

\proof{The only thing worth checking is that the
sequence is continuous. However, since the Gamma function is
infinitely differentiable at all points other than the
nonnegative integers, the map $a\mapsto s(a,k)$ is
continuous. Hence, the map defined by \eref{star63} is continuous.}

The nonzero elements of any
graded sequence form a pseudobasis for ${\cal I}$, and  the 
restriction of a graded sequence
to any particular $\alpha $ forms a pseudobasis for
${\cal I}^{\alpha}$.
Thus, for any pair of graded sequences $\seq{p}{a}{\alpha }$ and
$\seq{q}{a}{\alpha }$,
there is a unique continuous linear operator $\theta $ such
that $\theta p_{a}^{\alpha}(x) = q_{a}^{\alpha}(x)$. We study such
operators in detail in \cref{main-end}.

\begin{thm}\label{Roman2}
Let $\seq{p}{a}{\alpha}$ be a graded sequence. Then
every logarithmic series $g(x)\in {\cal I}$ can be uniquely written as an
expression 
$$ g(x)=\sum _{a,\alpha }b_{a}^{\alpha }p_{a}^{\alpha}(x) $$
where $g^{\alpha}(x)\in {\cal I}^{\alpha}$. 
This is a convergent expansion in the topology of Noetherian series, and an
asymptotic expansion in the topology of the complex numbers as $x$ tends
towards infinity.$\Box $
\end{thm}

We are able to actually calculate the constants $b_{a}^{\alpha }$ in many
cases using such results as Theorems~\ref{GLTF} and the corresponding result in
\cite{ch5}.

\ssect{Roman Graded Sequences}
In this section, we introduce the central concept of this work. It
is known that the operational
calculus of formal differential operators is intimately
associated with sequences of polynomials of binomial type,
that is, with sequences of polynomials $p_{n}(x)$ satisfying
the binomial identity (Definition~\ref{CDefBT})
$$p_{n}(x+a)= \sum_{k=0}^{n}{n\choose k }p_{k}(x)p_{n-k}(a)$$
A good many sequences of polynomials occurring in
combinatorics and in the theory of special functions turn
out to be of binomial type. For example, 
 the powers $x^{n}$,
   the lower factorial $(x)_{n}$,
 the upper factorial $(x)^{n}$, 
   the Abel polynomials $A_{n}(x)$,
 the LaGuerre polynomials $L_{n}(x)$, and
 the inverse-Abel polynomials $\mu _{n}(x)$  are all sequences of polynomials
of binomial type. We give here the logarithmic
generalization of this notion;  such graded sequences of
formal power series of logarithmic type are called Roman
graded sequences. We derive five equivalent characterizations
of such graded sequences.
We anticipate the fact (Theorem~\ref{TFAE}) that the five
notions introduced below coincide:
\begin{enumerate}
\item Roman graded sequence (Definition~\ref{defRoman}),
\item associated graded sequence (Definition~\ref{defAss}),
\item basic graded sequence (Definition~\ref{defBas}), and
\item conjugate graded sequence (Definition~\ref{defConj}),
\item graded sequence of logarithmic binomial type.
(Definition~\ref{defBT})
\end{enumerate}

Wee motivate the definition of a Roman graded sequence, by
 deriving a formula for the action of a product of two 
Artinian operators on the harmonic logarithm.

\begin{prop}
\begin{enumerate}
\item Let $f(\D )$ and $g(\D )$ be Artinian operators.
Then for all $a$ and for all $\alpha \neq  (0)$ we have the
following finite sum:
\begin{equation}\label{finitesum}
\action{\alpha }{f(\D )g(\D )\lambda _{a}^{\alpha}(x)}
=\sum_{b}\romcoeff{a}{b} \action{\alpha }{f(\D
)\lambda_{b}^{\alpha}(x)} \action{\alpha }{g(\D )\lambda
_{a-b}^{\alpha}(x)}.
\end{equation}
\item Similarly, let $f(\D )$ and $g(\D )$ be differential
operators. Then \eref{finitesum} holds for all  vectors $\alpha . $
\end{enumerate} 
\end{prop}

\proo{Let  $c_{a}=[\D ^{a}]f(\D )$, and $d_{a}=[\D
^{a}]g(\D )$. The product of the two series is given by the sum
$$f(\D )g(\D )= \sum_{b}\left( \sum _{a} c_{a}d _{a-b} \right) \D^{b}.$$
Hence, by Theorem~\ref{action}, 
\begin{eqnarray*}
\action{\alpha }{f(\D )g(\D )\lambda _{b}^{\alpha}(x)}
&=& \rn{b}! \sum _{a}c_{a}d_{b-a}\\*
&=& \sum_{a} \frac{\rn{b}!}{\rn{a}!\rn{b-a}!} (\rn{a}!c_{a})
(\rn{b-a}!d_{b-a})\\* 
&=&\sum_{a} \romcoeff{b}{a}
\action{\alpha}{f(\D)\lambda_{a}^{\alpha}(x)}
\action{\alpha}{g(\D )\lambda _{b-a}^{\alpha}(x)}.\Box
\end{eqnarray*}}

The extension to  products of more than two operators
  follows easily by induction.

We introduce Roman graded sequences by the following definition. It
will shortly be seen that simpler alternate definitions can
be given.

\begin{defn} \label{defRoman}
\bold{Roman Graded Sequences} Let $\seq{p}{a}{\alpha}$
be a graded sequence of
formal power series of logarithmic type.
The graded sequence is a {\em Roman graded sequence}
 if for all $a$ and $\alpha $, we have the following finite
sum 
\begin{equation}\label{*}
\action{\alpha}{f(\D )g(\D )p_{a}^{\alpha}(x)}=\sum _{b}
\romcoeff{a}{b} 
\action{\alpha}{f(\D)p_{b}^{\alpha}(x)}
\action{\alpha}{g(\D )p_{a-b}^{\alpha}(x)}
\end{equation}
where $f(\D )$ and $g(\D  )$ are Artinian operators if
$\alpha \neq  (0)$, and are differential operators if $\alpha
=(0)$.
\end{defn}

For example, $\seq{\lambda }{a}{\alpha}$ is a Roman graded sequence.

\begin{prop}
\label{easier}\label{coeff-Rom}
A graded sequence is Roman if and only if \eref{*} 
holds when $f(\D )=\D ^{a}$ and $g(\D )=\D ^{b}$.
\end{prop}

\proof{Linearity and continuity.}

\begin{prop} \label{starstar}
A graded sequence $\seq{p}{a}{\alpha}$ 
 with coefficients
$d_{ab}=[\lambda_{b}^{\alpha}(x)]p_{a}^{\alpha }(x)$
is Roman
if and only if for all $a,$ $b,$ and $c$:
\begin{equation}\label{star2}
\romcoeff{a+b}{a}d_{c,a+b} = \sum_{e} \romcoeff{c}{e} d_{eb}d_{c-e,a}.
\end{equation}
\end{prop}

\proof{We demonstrate the ``only if;'' the reasoning for the
other implication is similar.
\begin{eqnarray*}
\action{\alpha}{\D ^{a}\D ^{b}p_{c}^{\alpha}(x)}
&=&\rn{a+b}!d_{c,a+b}\\*
&=&\rn{a}!\rn{b}! \sum _{e}\romcoeff{c}{e}d_{e,a}d_{c-e,b}\\*
&=& \sum_{e}\romcoeff{c}{e}\action{\alpha}{\D ^{a}p_{k}^{\alpha}(x)}
\action{\alpha}{\D ^{b}p_{c-e}^{\alpha}(x)}.
\end{eqnarray*}
Note that in the case $\alpha =(0)$, we need only consider
nonnegative integers $a$ and $b$, so the same
argument applies.}

\ssect{Associated Graded Sequences}\label{sec-ass} 

We proceed to derive an altogether different characteristic
property of Roman graded sequences, which uses delta operators.
Such 
differential operators may be viewed as playing the role of
the derivative---much like the forward difference operator
$\Delta 
$ (Definition~\ref{FDO}) acts on the polynomial sequence
of lower factorials $(x)_{n}=x(x-1)\cdots (x-n+1)$.

\begin{prop}
\label{assthm}
Let $f(\D )$ (and $n$ an integer) be a delta operator in $\Lambda^{+}$. 
Then there is a unique graded sequence of formal power series of
logarithmic type, $\seq{p}{a}{\alpha}$, such that for all $a,b$ and $\alpha $
\dclab{ass}{ \action{\alpha}{f(\D)^{a}p_{b}^{\alpha}(x)} }{ \rn{a}!\delta _{ab}
}{\action{\alpha}{f(\D)^{a;n}p_{b}^{\alpha}(x)}}{\rn{a}!\delta _{ab} .} 
\end{prop} 

\proof{ {\bf (Uniqueness)} Spanning argument
(Proposition~\ref{spanning}.)  

{\bf (Existence)} We need only consider the continuous case. Let $d_{b}^{(a)}=
[\D ^{b}]f(\D )^{a;n}$ be the coefficients of $f(\D )^{a;n}$, and let
$c_{b}=c_{b}^{(1)}=[\D^{b}] f(\D )$. Thus, $c_{a}^{(a)}\neq 0$ for all real
numbers $a$. Similarly, let $d_{ab}^{\alpha
}=[\lambda_{b}^{\alpha}(x)]p_{a}^{\alpha }(x)$ be the coefficients of
$p_{a}^{\alpha }(x)$.
In this notation, \eref{ass} is equivalent to the equation
$$ \sum _{b}\rn{b}!c_{b}^{(a)}d_{eb}^{\alpha} = \rn{e}!\delta _{eb}. $$ 
Hence, the  coefficients
$d_{ab}^{\alpha}$ can be computed recursively as:
 $$d_{ab}^{\alpha}=
\frac{1}{\rn{b}!a_{b}^{(b)}}
\left(\rn{a}!\delta _{ab}
-\sum_{e>b}\rn{e}!c_{e}^{(b)}d_{ab}^{\alpha}\right).$$
 The recursion is well
defined since $d_{ab}^{\alpha}$ is a Noetherian sequence in $b$.
Finally,  notice that $\seq{p}{a}{\alpha}$ is regular by
symmetry, and that $\deg(p_{n}^{\alpha}(x))=n$
since $b_{nn}^{\alpha}= {\rn{n}!}/{\rn{n}!a_{n}^{(n)}}\neq
0$. 
Finally, note that the map defined by \eref{star63} 
is continuous since $a\mapsto f(\D )^{a;n}$ is continuous.}
 
\begin{defn}\bold{Associated Graded Sequence}\label{defAss}
Let $f(\D )$ be a delta operator (and $n$ an integer). The unique
graded sequence mentioned in Proposition~\ref{assthm} is
called the {\em ($n\th$) associated 
graded sequence} of the delta operator $f(\D)$.  We  also say that
the sequence is $n-$associated with $f(\D )$.
(The term $n\th $ applies only in the continuous case.)
\end{defn}

For example, by Proposition~\ref{action}, $\seq{\lambda }{a}{\alpha}$ is the
standard graded sequence associated with 
the delta operator $\D $. 

We can now generalize Theorem~\ref{exp1} to explicitly
determine the coefficients in the expansion of an arbitrary
Artinian operator in terms of the powers of a delta operator:

\begin{thm}\bold{Expansion Theorem}
\label{GET}
Let the graded sequence $\seq{p}{a}{\alpha}$ be
$n$-associated with the delta operator
$f(\D)$. Then for all Artinian operators
$g(\D )$  and  vectors $\alpha \neq  (0)$
we have the following convergent sum.
\dclab{GETeq}{ g(\D ) }{\sum_{k} \frac{\action{\alpha}
{g(\D )p_{k}^{\alpha}(x)}}{\rn{k}!} f(\D )^{k}}{g(\D )}{\sum_{a}
\frac{\action{\alpha} 
{g(\D )p_{a}^{\alpha}(x)}}{\rn{a}!} f(\D )^{a;n}. }

When $\alpha =(0)$, \eref{GETeq} holds for all differential operators $g(\D )$.
\end{thm} 

\proof{By Definition~\ref{augment} we have,
 \begin{eqnarray*}
\action{\alpha}{\sum _{a}
 \frac{\action{\alpha}{g(\D )p_{a}^{\alpha}(x)}}{\rn{a}!}f(\D )^{a;n}
p_{b}^{\alpha}(x)}
 &=&\sum _{a}
\frac{\action{\alpha}{g(\D
)p_{a}^{\alpha}(x)}}{\rn{a}!}\action{\alpha}{f(\D
)^{a;n}p_{b}^{\alpha}(x)}\\*
 &=& \action{\alpha}{g(\D)p_{n}^{\alpha}(x)}.
\end{eqnarray*} 
The conclusion  follows by the spanning argument.}

Dually, we obtain the explicit form of the expansion of an
arbitrary formal series of logarithmic type as a linear
combination of elements of a Roman graded sequence. This gives
a useful generalization of Theorem~\ref{ltf}.

\begin{thm}\bold{Logarithmic Taylor's Theorem}
\label{GLTF}
Let $\seq{p}{a}{\alpha}$ be the ($n\th $) graded sequence
associated with the delta 
operator $f(\D )$. Then for every formal power series of
logarithmic type $p(x)\in {\cal I}^{+}$ we have the
following convergent sum
\dc{ p(x)}{\sum_{\alpha\neq (0) }\sum_{k}
{\frac{\action{\alpha}{f(\D)^{k} p(x)}}
{\rn{k}!}}  p_{k}^{\alpha}(x)} {p(x)}{\sum_{\alpha\neq (0) }\sum_{a}
{\frac{\action{\alpha}{f(\D)^{a;n} p(x)}}
{\rn{a}!}}  p_{a}^{\alpha}(x).}
\end{thm}

\proo{We need only consider the continuous case, but first a simple lemma. For
all $p(x),q(x)\in {\cal I}^{+}$, by Part {\bf (1)} of the spanning argument
(Proposition~\ref{spanning}), it is clear that  $ p(x)=q(x) $
if and only if for all $a$ and $\alpha \neq (0)$,
$\action{\alpha}{\D ^{a}p(x)}=\action{\alpha}{\D ^{a}q(x)}.$

When we apply the Expansion Theorem where $g(\D )$ is
set equal to $\D ^{a}$, we find that
$$ \D ^{a}=\sum _{b}
\frac{\action{\alpha}{\D^{a}p_{b}^{\alpha}(x)}}{\rn{b}!}f(\D )^{b;n}.$$
Thus,
\begin{eqnarray*}
\action{\alpha}{\D^{a}p(x)}
&=&\sum_{b}\frac{\action{\alpha}{\D ^{a}p_{b}^{\alpha}(x)}}{\rn{b}!}
\action{\alpha}{f(\D)^{b;n}p(x)}\\*
&=&\action{\alpha}{\D ^{a}\sum _{b}
p_{b}^{\alpha}(x)\frac{\action{\alpha}{f(\D )^{b;n}p(x)}}{\rn{b}!}}.
\end{eqnarray*}
Therefore, the above remarks, we conclude that
$$ p(x)=\sum_{\alpha\neq (0) }\sum_{a}
{\frac{\action{\alpha}{f(\D)^{a}p(x)}}{\rn{a}!}}
p_{a}^{\alpha}(x).\Box $$ }

In the continuous case, where this is relevant, we may now classify standard
associated sequences (Definition~\ref{sfpslt}).

\thmplus{C}{\begin{prop}\label{Stand}
The $n\th $ associated sequence of a delta operator $f(\D )$
is standard if and only if the leading coefficient of $f(\D
)$ is a positive real number and $n=0.\Box $
\end{prop}}

Such delta operators are be called {\em standard} operators.\index{Standard
Operators}

\ssect{Basic Graded Sequences}

\begin{defn} \label{defBas}
\bold{Basic Graded Sequence} Let $f(\D )$ be a delta operator. A
 graded sequence of formal power series of logarithmic type
$\seq{p}{a}{\alpha}$  
is called the {\em ($n\th$) basic graded sequence} for $f(\D) $ if
for all $\alpha $,
\begin{enumerate}
\item \label{prop1} 
$\action{\alpha}{p_{0}^{\alpha}(x)}=1$;
\item \label{prop2}
For all $a>0$ (and thus all $a\neq 0$),
$\action{\alpha}{p_{a}^{\alpha}(x)}=0$; and
\item\label{prop3}
\begin{description}
\item[Discrete]
For all integers $n$, and all vectors $\alpha $ of integers with finite
support, 
$f(\D)p_{n}^{\alpha}(x)=\rn{n}p_{n-1}^{\alpha}(x)$. 
\item[Continuous] For all real
numbers $a$ and $b$, and all nonzero vectors $\alpha $ of real numbers with
finite support,
\begin{equation}\label{star69}
f(\D )^{b;n}p_{a}^{\alpha }(x)=\frac{\rn{a}!}{\rn{a-b}!}
p_{a-b}^{\alpha }(x).
\end{equation}
\end{description}
\end{enumerate}
\end{defn}

\begin{thm}\label{bas-ass}
\begin{enumerate}
\item Let $\seq{p}{a}{\alpha}$ be a logarithmic graded sequence. Such a
graded sequence is the  ($n\th $) basic graded sequence
for the delta operator $f(\D )$ if and only if it is
the ($n\th $) associated graded sequence for the delta operator $f(\D )$.
\item Every delta operator has a unique ($n\th$) basic sequence.
\item Every basic sequence is basic for a unique delta
operator (and integer $n$).
\end{enumerate}
\end{thm}

\proof{{\bf  (1: If)} Properties~\ref{prop1} and \ref{prop2}
of Definition~\ref{defBas} follow from Definition~\ref{defAss}.
Property~\ref{prop3}
follows from the following series of
equalities, and the spanning argument as shown here in the continuous case:
\begin{eqnarray*} 
\action{\alpha} {f(\D )^{a;n}\left(f(\D
)^{b;n}p_{c}^{\alpha}(x)\right)} 
&=&\action{\alpha} {f(\D )^{a+b;n}p_{c}^{\alpha}(x)}\\*
&=&\rn{c}!\delta _{c,a+b}\\
&=&\frac{\rn{c}!}{\rn{c-b}!}\left(\rn{c-b}!\delta _{c-b,a}\right)\\*
&=&\frac{\rn{c}!}{\rn{c-b}!} \action{\alpha} {f(\D
)^{a;n}p_{c-b}^{\alpha}(x)} 
\end{eqnarray*}

{\bf (1: Only if)} We proceed differently in the continuous and discrete cases:
\begin{description}
\item[Discrete] By induction we have, $f(\D
)^{k}p_{n}^{\alpha}(x)={\rn{n}!}p_{n-k}^{\alpha}(x)/{\rn{n-k}!}$. Hence,
\begin{eqnarray*} 
\action{\alpha} {f(\D )^{k}p_{n}^{\alpha}(x)}
&=&\action{\alpha} {\frac{\rn{n}!}{\rn{n-k}!}p_{n-k}^{\alpha}(x)}\\*
&=&\rn{n}!\delta _{n-k,0}\\*
&=&\rn{n}!\delta _{n,k}.
\end{eqnarray*}
Hence, $\seq{p}{n}{\alpha}$ is the associated graded sequence of $f(\D
)$ by definition. 
\item[Continuous] By \eref{star69}, $\action{\alpha
}{f(\D )^{b;n}p_{a}^{\alpha }(x)}=\frac{\rn{a}!}{\rn{a-b}!}
\delta _{a-b,0}=\rn{a}!\delta _{ab}.$
\end{description} 

{\bf (2 and 3)} Immediate from 1.}

The simplest example of a basic  graded sequence is the
graded sequence of 
harmonic logarithms $\seq{\lambda}{a}{\alpha}$; 
by Theorem~\ref{deriv}, it is the standard basic graded sequence for the delta
operator $\D $. It can be viewed as the natural
 logarithmic extension of the
sequence of powers of $x$. More generally, every sequence of
binomial type (and every factor sequence) has a natural
logarithmic extension into a basic
graded sequence related to  the same delta
operator.

$\D $ is invertible on ${\cal I}^{+}$. Hence,
just as in the continuous case (\eref{star69}), we have:

\thmplus{D}{\begin{por}
Let $\seq{p}{n}{\alpha}$ be the associated graded sequence
of formal power series of logarithmic type
of the delta operator $f(\D )$, then for all integers $n$ and $k$, and all
nonzero vectors $\alpha$ with finite support of integers.
\begin{equation}\label{justas}
f(\D)^{k}p_{n}^{\alpha}(x) = \frac{\rn{n}!}{\rn{n-k}!}p_{n-k}^{\alpha}(x)
\end{equation}
Moreover, when $k$ is a nonnegative integer,
\eref{justas} holds for all $\alpha .\Box $
\end{por}}

We may now connect the notion of an associated graded sequence with
that of a Roman graded sequence.

\begin{prop}\label{consqnt}
Let $\seq{p}{a}{\alpha}$ be an associated graded sequence, 
and let $f(\D )\in \Lambda^{+} $ be an Artinian operator. Then we
have the following convergent sum 
$$ f(\D )p_{a}^{\alpha}(x)= 
\sum _{b}\romcoeff{a}{b}
\action{\alpha} {f(\D)p_{b}^{\alpha}(x)}p_{a-b}^{\alpha}(x)$$
for all $\alpha \neq (0)$, and also for $\alpha =(0)$ if
$f(\D )$ is in fact a differential operator. 
\end{prop}

\proof{We need only consider the continuous case, suppose $\seq{p}{a}{\alpha}$
is $n$-associated with the delta operator $g(\D )$. Then for all $b$
\begin{eqnarray*}
g(\D )^{b;n}p_{a}^{\alpha}(x) &=&
\frac{\rn{a}!}{\rn{a-b}!}p_{a-b}^{\alpha}(x)\\* 
&=& \sum _{c}\romcoeff{a}{c}\action{\alpha}
{g(\D)^{b;n}p_{c}^{\alpha}(x)} p_{a-c}^{\alpha}(x).
\end{eqnarray*}
Since $\{g(\D )^{b;n}: b\in {\bf R}\}$ is a pseudobasis, we
can use 
continuity and linearity to 
replace $g(\D)^{b;n}$ by $f(\D ).$}

The following identity is classically proven \cite{RR} by introducing a new
variable. However, we have a much simpler proof.

\begin{cor}\label{questionable}
Let $\seq{p}{a}{\alpha}$ be  an associated graded sequence with coefficients 
$$c_{ab}=[\lambda _{b}^{\alpha}(x)]p_{a}^{\alpha}(x).$$
Then for all $a,$ $b,$ and $\alpha $
we have
$$ \D ^{b}p_{a}^{\alpha}(x)= \sum _{d\geq b}
\frac{\rn{a}!c_{db}}{\rn{a-d}!} p_{a-d}^{\alpha}(x) .$$
\end{cor}

\proof{Proposition~\ref{consqnt}.}

As another consequence of Proposition~\ref{consqnt} we obtain:

\begin{prop}\label{spl-ass}
A  graded sequence of
formal power series of logarithmic type $\seq{p}{a}{\alpha}$
is an associated graded sequence if and only if it
is a Roman graded sequence. 
\end{prop}

\proo{{\bf  (Only if)} Let $f(\D )$ and $g(\D )$ be Artinian operators. By
Proposition~\ref{consqnt},
$$ f(\D )p_{a}^{\alpha}(x)= 
\sum _{b}\romcoeff{a}{b}
\action{\alpha} {f(\D)p_{b}^{\alpha}(x)}p_{a-b}^{\alpha}(x).$$
Thus,
\begin{eqnarray*}
\action{\alpha} {f(\D )g(\D )p_{a}^{\alpha }(x)} &=& 
\action{\alpha} {g(\D )\left(\sum _{b}\romcoeff{a}{b}
\action{\alpha} {f(\D)p_{b}^{\alpha}(x)}p_{a-b}^{\alpha}(x)
 \right)}  \\*
&=& \sum _{b}\romcoeff{a}{b}\action{\alpha} {f(\D )p_{b}^{\alpha }(x)}
\action{\alpha} {g(\D )p_{a-b}^{\alpha }(x)}.
\end{eqnarray*}

{\bf (If)} Conversely, let $\seq{p}{a}{\alpha}$ be a Roman graded
sequence. We define a
sequence  of Artinian operators $f_{b}(\D )$ by the relation 
$$\action{(1)}{f_{b}(\D )p_{a}^{(1)}(x)}=\rn{a}!\delta _{ab}.$$ 
It suffices to show that 
\dc{f_{b}(\D )}{f(\D )^{b;n}}{f_{b}(\D )}{f(\D )^{b}}
 for some delta operator $f(\D )$ (and some integer $n$).

 By the spanning argument, $f_{b}(\D )$ is well defined. Now, 
$$\action{(1)}{f_{b}(\D )\lambda _{a}^{(1)}(x)}=0$$ 
for $a<b$, and 
$$\action{(1)}{f_{b}(\D )\lambda _{b}^{(1)}(x)}\neq 0,$$
hence
$\deg(f_{b}(\D ))=b$. In particular, $f_{1}(\D )$ is a delta
operator. Since 
$\seq{p}{a}{\alpha}$ is a Roman graded sequence, we infer that 
\begin{eqnarray*}
\action{(1)} {f_{b}(\D )f_{c}(\D )p_{a}^{(1)}(x)}
&=&\sum _{d}\romcoeff{a}{d} 
\action{(1)} {f_{b}(\D)p_{d}^{(1)}(x)} 
\action{(1)} {f_{c}(\D )p_{a-d}^{(1)}(x)}\\*
&=&\rn{a}!\delta _{a,b+c}\\*
&=&\action{(1)} {f_{b+c}(\D )p_{a}^{(1)}(x)}. 
\end{eqnarray*}
By the spanning argument,
$f_{i}(\D )f_{j}(\D )=f_{i+j}(\D )$. Hence, 
\begin{description}
\item[Discrete] By induction, $f_{n}(\D )=f_{1}(\D )^{n}$.
\item[Continuous] By the characterization of exponentiation in \cite{ch1},
there is an integer $n$ such that $f_{b}(\D )=f_{1}(\D )^{b;n}$ for all $n.\Box
$ 
\end{description}}

Thus, every Roman graded sequence is
associated  with a unique delta operator (and integer) and {\em visa versa}.

\ssect{Conjugate Graded Sequences}
Each delta operator (and integer) has another graded sequence associated with
it: its conjugate graded sequence. We see that these sequences are Roman;
however, they are not associated with the same delta operator for which they
are conjugate.

\begin{defn} \label{defConj}
\bold{Conjugate Graded Sequence} 
Let $f(\D )\in \Lambda^{+} $ be a delta operator. 
Its ($n\th$) {\em conjugate graded sequence} $\seq{q}{n}{\alpha}$ is
defined as 
\dclab{star71}{q_{n}^{\alpha}(x)}{\sum _{k\leq n}\frac{\action{\alpha}{f(\D
)^{k}\lambda _{n}^{\alpha}(x)}}{\rn{k}!}\lambda
_{k}^{\alpha}(x)}{q_{a}^{\alpha}(x)}{\sum _{b\leq a}
\frac{\action{\alpha}{f(\D)^{b;n} \lambda_{b}^{\alpha}(x)}}
{\rn{b}!} \lambda _{b}^{\alpha}(x)}
for all $a$ and $\alpha $.
\end{defn}

Indeed, for all delta operators $f(\D )$ (and integers $n$),
the graded sequence $\seq{q}{a}{\alpha}$ as defined by \eref{star71}
automatically meets the conditions of Definition~\ref{sfpslt}.

The canonical example of a conjugate  graded sequence is the
graded sequence of
harmonic logarithms; it is the standard conjugate
graded sequence of the delta operator $\D $, since by Theorem~\ref{ltf}, 
$$\lambda _{a}^{\alpha}(x) = 
\sum _{b}\action{\alpha}{\D ^{b}\lambda _{a}^{\alpha}(x)}
\lambda_{b}^{\alpha}(x)/ \rn{b}!.$$

In the above example, the $n\th $ conjugate graded sequence
of a delta 
operator is a Roman graded sequence. This fact is true in
general: 

\begin{prop}\label{spl-conj}
A graded sequence
of formal power series of logarithmic type
$\seq{q}{a}{\alpha}$ is Roman 
if and only if it is the ($n\th$) conjugate graded sequence of a delta
operator. Moreover, the delta operator (and integer) are unique.
\end{prop}

\proof{{\bf (If)} We need only consider the continuous case.
Let $\seq{q}{a}{\alpha}$
be the $n\th$ conjugate graded sequence of the delta
operator $f(\D )$. Now, by \eref{star71} the coefficients $c_{ab}$ of
$p_{a}^{\alpha }(x)$ are given by
$$c_{ab}=[\lambda_{b}^{\alpha}(x)]p_{a}^{\alpha }(x) =
\frac{\action{\alpha}{f(\D )^{b;n}\lambda _{a}^{\alpha}(x)}}{\rn{b}!}.$$ 
It suffices to show that the
$c_{ab}$ satisfy \eref{star2}. For any $b_{1},b_{2}$,
\begin{eqnarray*} 
\sum _{d}\romcoeff{a}{d}c_{db_{2}}c_{a-d,b_{1}}
&=&\sum _{d}\romcoeff{a}{d}
\frac{\action{\alpha}{f(\D)^{b_{2};n} \lambda
_{d}^{\alpha}(x)}} {\rn{b_{2}}!} 
\frac{\action{\alpha}{f(\D)^{b_{1};n}\lambda
_{a-d}^{\alpha}(x)}}{\rn{b_{1}}!}\\* 
&=&\frac{\action{\alpha}{f(\D )^{b_{1}+b_{2};n}
\lambda_{a}^{\alpha}(x)}}{\rn{b_{1}}!\rn{b_{2}}!} 
\end{eqnarray*}
since $\seq{\lambda}{a}{\alpha}$ is a Roman
graded sequence. The last expression equals
$\romcoeff{b_{1}+b_{2}}{b_{1}}c_{a,b_{1}+b_{2}}.$
Hence, \eref{star2} holds.

{\bf (Only if) } Conversely, suppose that $\seq{q}{a}{\alpha}$ is
a Roman graded sequence. 
Let $c_{ab}$ denote the coefficients of $q_{a}^{\alpha }(x)$
$$c_{ab}= [\lambda_{b}^{\alpha}(x)] q_{a}^{\alpha}(x).$$ 
Define the Artinian operator $f_{d}(\D )$ by 
$$\action{\alpha}{f_{d}(\D )\lambda
_{a}^{\alpha}(x)}=\rn{d}!c_{ad}.$$
 By the spanning argument,
this condition defines $f_{d}(\D )$. As in the proof of
Proposition~\ref{spl-ass} which employs a similar technique, it suffices to
observe that $f_{d}(\D )=f(\D )^{d}$ (resp. $f(\D )^{d;n}$) for some delta
operator $f(\D )$ (and integer $n$).

In fact, we have:
$$ f_{1}(\D )=\sum _{a}\frac{c_{a,1}}{\rn{a}!}\D ^{a}.$$ 
Thus, $f_{1}(\D )$ is a
delta operator. Now,
$$\action{\alpha}{f_{b_{1}+b_{2}}(\D)\lambda _{a}^{\alpha}(x)} 
=\rn{b_{1}+b_{2}}!b_{a,b_{1}+b_{2}}$$ 
by definition, and that expression equals 
$$\sum _{d}\romcoeff{a}{k}
\action{\alpha}{f_{b_{1}}(\D )\lambda_{d}^{\alpha}(x)}
\action{\alpha}{f_{b_{2}}(\D )\lambda _{a-d}^{\alpha}(x)}$$ by
\eref{star2}. This in turn equals
$$\action{\alpha}{f_{b_{1}}(\D )f_{b_{2}}(\D )\lambda
_{a}^{\alpha}(x)}$$ since $\seq{\lambda }{a}{\alpha}$ is a
Roman graded sequence. In other words,
$$ \action{\alpha } {f_{b_{1}+b_{2}}(\D )\lambda _{a}^{\alpha}(x)}
=\action{\alpha } {f_{b_{1}}(\D )f_{b_{2}} (\D  )\lambda
_{a}^{\alpha}(x)} ,$$ 
so by the spanning argument
 $f_{i+j}(\D )=f_{i}(\D )f_{j}(\D )$.

\begin{description}
\item[Discrete] Hence,
by induction, $f_{n}(\D )=f_{1}(\D )^{n}$.
\item[Continuous] Notice first that $d\mapsto
f_{d}(\D )$ is a continuous map.
Hence, by the characterization of exponentiation in \cite{ch1}, there is an integer $n$ such that
$f_{b}(\D )=f_{1}(\D )^{b;n}$ for all $n$. 
\end{description}

In other words, $\seq{q}{a}{\alpha}$ is
the ($n\th $) conjugate graded sequence of the delta operator $f_{1}(\D )$ (and
no other).}

In the continuous case where this is relevant, we may now classify standard
conjugate sequences.

\thmplus{C}{\begin{prop}\label{Stand-c}
The $n\th $ conjugate sequence of a delta operator $f(\D )$
is standard if and only if  $f(\D)$ is a standard operator and $n=0$.$\Box $
\end{prop}}

\ssect{Graded Sequences of Logarithmic Binomial Type}

We motivate our final reformulation of the Roman condition by  briefly
recalling some facts about sequences of polynomials of binomial type.

\begin{defn}
\bld{Polynomial Sequence of Binomial Type} 
 A sequence $\left(p_{n}(x)\right)_{n\geq 0}$ of polynomials\label{CDefBT} 
 is of {\em binomial type}\index{Binomial Type,Polynomial Sequence of}
if each polynomial $p_{n}(x)$ is of
degree $n$, and
for all field elements $a$, and all nonnegative integers $n$, 
$$ p_{n}(x+a)= \sum _{k=0}^{n}{n\choose k}p_{k}(a)p_{n-k}(x).$$
\end{defn}

It is a  basic result of Umbral calculus (see for example \cite{RR})
that every sequence
of binomial type is associated with a unique delta operator
and {\it vice versa}. $(p_{n}(x))_{n\geq 0}$ is  the
associated
sequence of the delta operator $f(\D )$ if and only if
$f(\D)p_{n}(x)=np_{n-1}(x) $ for positive $n$,
and $p_{n}(0)=\delta _{n,0}$ for nonnegative $n$.

In the present theory, the
logarithmic analogs of sequences of binomial type
are the  graded sequences of binomial type defined below. As before, this 
 definition turns out to be equivalent to that of a  Roman graded
sequence. 

\begin{defn} \label{DefBT}
\label{defBT}
\bld{Graded Sequence of Logarithmic Binomial Type}
A graded sequence of formal
power series of logarithmic type $\seq{p}{a}{\alpha}$ is of
{\em logarithmic binomial type}\index{Binomial Type} if for
all  $a$ and $\alpha $, and all complex numbers $z$, $\sum
_{b}\romcoeff{a}{b} \action{(0)}{E^{z} p_{b}^{(0)}(x)}
p_{a-b}^{\alpha }(x)$
is a convergent sum which equals $E^{z}p_{a}^{\alpha }(x)$.
\end{defn}

\discrete{Note that in the discrete case $\action{(0)}{E^{z} p_{b}^{(0)}(x)}$
is merely $p_{b}^{(0)}(z)$. Thus,
for $\alpha =(0)$, we obtain simply  the definition of a sequence of
polynomials of
binomial type. For $\alpha =(1)$ and $a$ a negative integer, we obtain 
factor sequences. \cite{UC}

 Thus, the notion of a 
 graded sequence of logarithmic binomial type  subsumes both the notion of a
sequence of polynomials of binomial type, and the notion of a factor
sequence. In view of the present theory, these older notions can be seen as
obsolete.}

We next prove that, as promised, every 
graded sequence of logarithmic binomial type
is a Roman graded sequence, and conversely.

\begin{thm}\label{bas-bin}
A logarithmic graded sequence
$\seq{p}{a}{\alpha}$ is the ($n\th $) basic graded sequence
for some delta 
operator $f(\D )$ (and integer $n$) if and only if
 it is a graded sequence of logarithmic binomial type.
\end{thm}

\proo{{\bf (Only if)} Proposition~\ref{consqnt}.

{\bf (If)} To prove  that a graded sequence is basic, we must demonstrate each
of the three properties enumerated in Definition~\ref{defBas}. 

{\bf (1)}  By Definition~\ref{DefBT},
$$ p_{0}^{(0)}(x)=p_{0}^{(0)}(0)p_{0}^{(0)}(x), $$
and ${\cal I}$ is a field, so we have
$$ p_{0}^{(0)}(x)=1. $$

{\bf (2)} Whereas for $a$ positive, we have
$$ 0=p_{a}^{(0)}(x)-p_{a}^{(0)}(x)=
\sum_{b>0}\romcoeff{a}{b} \action{(0)}{p_{b}^{(0)}(x)}
p_{a-b}^{(0)}(x),$$
and the $p_{a}^{(0)}(x)$ form a pseudobasis for ${\cal I}^{(0)}$, so 
$\action{(0)}{p_{b}^{(0)}(x)}=0$ for $b>0$. Thus, by regularity, 
$$ \action{\alpha} {p_{a}^{\alpha}(x)}=\delta _{n,0} $$ 
for all $a$ and $\alpha $.

{\bf (3)} Define a sequence of continuous, linear operators $Q^{b}$ by
the identity
$$ Q^{b}p_{a}^{\alpha}(x)= \frac{\rn{a}}{\rn{a-b}}
p_{a-b}^{\alpha}(x). $$ 
Clearly, $Q^{c}Q^{c}=Q^{b+c}$, and the map $b\mapsto Q^{b}$
is continuous in the operator topology, so by
the characterization of exponentiation in \cite{ch1}, it remains only
to show that $Q=Q^{1}$ is a delta operator.

By inspection $Q$ is regular, and lowers the degree of any
logarithmic series by one. We now demonstrate via the
following string of identities that
$Q$ is $\D $-invariant:
\begin{eqnarray*}
QE^{z}p_{a}^{\alpha}(x) 
&=& Q\sum _{b} \romcoeff{a}{b}
\action{(0)}{E^{z}p_{b}^{(0)}(x)} p_{a-b}^{\alpha}(x)\\*
&=& \sum _{b} \romcoeff{a}{b}\rn{a-b}
\action{(0)}{E^{z}p_{b}^{(0)}(x)} p_{a-b-1}^{\alpha}(x)\\
&=& \sum _{b} \romcoeff{a-1}{b}\rn{a}
\action{(0)}{E^{z}p_{b}^{(0)}(x)} p_{a-b-1}^{\alpha}(x)\\*
&=&E^{z}Qp_{a}^{\alpha}(x).\Box 
\end{eqnarray*}}

We summarize the results of the preceding sections in the
following theorem:

\begin{thm}\label{TFAE}
For any graded sequence of formal power series of logarithmic type
$\seq{p}{a}{\alpha}$, the following statements are equivalent:
\begin{enumerate}
\item $\seq{p}{a}{\alpha}$ is a Roman graded sequence.
(Definition~\ref{defRoman})
\item $\seq{p}{a}{\alpha}$ is the ($n\th$) associated graded sequence for some
unique delta operator $f(\D )$. (Definition~\ref{defAss}) 
\item $\seq{p}{a}{\alpha}$ is the ($m\th$) conjugate graded
sequence for some unique delta operator $g(\D )$. (Definition~\ref{defConj})
\item $\seq{p}{a}{\alpha}$ is the ($n\th$) basic graded sequence for some
unique delta operator $f(\D )$. (Definition~\ref{defBas})
\item $\seq{p}{a}{\alpha}$ is a graded sequence of logarithmic binomial type.
(Definition~\ref{defBT})
\item $\seq{p}{a}{\alpha}$ is a  graded sequence of formal power series
of logarithmic type whose coefficients
$$p_{n}^{\alpha}(x)=\sum _{k\leq n}b_{nk}\lambda _{k}^{\alpha}(x)$$
satisfy \eref{star2}.
\end{enumerate}

Furthermore, in this case, $\seq{p}{a}{\alpha}$ is associated with
and basic for the same delta operator $f(\D )$ (and integer $n$).
\end{thm}

\proof{Condition~1 is equivalent to Condition~6 by
Porism~\ref{coeff-Rom}, it is equivalent to Condition~3
by Proposition~\ref{spl-conj}, and it is equivalent to
Condition~2 by Proposition~\ref{spl-ass}. Next,
Condition~2 is equivalent to Condition~4 by
Theorem~\ref{bas-ass}. Finally, Condition~4 is equivalent to
 Condition~5 by Theorem~\ref{bas-bin}.}

\sect{Relations Among Roman Graded Sequences}\index{Roman
Graded Sequences} 
\label{main-end}
\begin{description}
\item[Continuous] In this chapter, we focus our attention
primarily upon standard operators, and their corresponding  $0\th $
associated and conjugate sequences which we have seen (Propositions~\ref{Stand}
and \ref{Stand-c}) are standard sequences.
\item[Discrete] In this chapter, we maintain the same level
of generality as before. However, our results are still less general than
the corresponding results in the continuous case.
\end{description}
\ssect{Transfer Operators}

Let $\seq{p}{a}{\alpha}$ be
the standard associated graded sequence for the delta
operator $f(\D )$. 
In view of Propositions~\ref{spl-conj}, and \ref{spl-ass}, the standard
conjugate graded sequence for the operator $f(\D )$ 
is in general some other Roman graded sequence, say
$\seq{q}{a}{\alpha}$. By Proposition~\ref{assthm}, this 
graded sequence is in turn the standard associated graded
sequence for another delta operator $g(\D )$. What is the relationship
between $f(\D )$ and $g(\D )$? And between $\seq{p}{a}{\alpha}$
and $g(\D )$?

We shall prove (Corollary~\ref{conj-ass}) the remarkable fact that the formal
power series $f(\D )$ and $g(\D )$ are inverses to each other in
the sense of LaGrange Inversion \cite{ch1}, and $\seq{p}{a}{\alpha}$
is the conjugate graded sequence for $g(\D )$. Actually, the results
we shall obtain are  more sweeping, and lead to a
powerful technique for establishing identities among formal
power series of logarithmic type.

We shall occasionally use the boldface {\bf p} to denote a
logarithmic graded sequence $\seq{p}{a}{\alpha}$. Thus,

\begin{defn}[Umbral Composition] \label{compos}
Let  ${\bf p}$ and ${\bf q}$ be two graded sequences, and
let $c_{ab}$ denote the coefficients of the ${\bf p}$ sequence
$c_{ab}=[\lambda_{b}^{\alpha}(x)]p_{a}^{\alpha }(x)$. 
We define the {\em umbral composition} of  the ${\bf q}$
graded sequence with the  
${\bf p}$ graded sequence
to be the graded sequence defined by the following
convergent summation
$$ \comp{q_{a}^{\alpha}}{p}=
\sum_{b}c_{ab}p _{b}^{\alpha}(x).$$
\end{defn}

\begin{prop}
\label{wduc}
Given two graded sequences $\seq{p}{a}{\alpha }$ and
$\seq{q}{a}{\alpha }$. Their composition $p_{a}^{\alpha
}({\bf q})$ is a well defined graded sequence.
\end{prop}

\proof{We need only consider the continuous case. 
By the Noetherian condition, their composition is well defined. The conditions
on its degree and order follow immediately. Continuity follows
since the composition of two continuous functions is itself continuous.}

By \cite{ch1},  $\D $ is the two-sided identity for the group of
delta operators under (0-)comp\-osi\-tion, by definition the
graded sequence of 
harmonic logarithms is the two-sided identity for the
semigroup formed by the operation of composition of
standard Roman graded sequences. We shall show that this semigroup is
actually a group, and that the two groups are naturally
isomorphic. The isomorphism is given by
the function which associates each delta operator with its
associated graded sequence. The crucial role in 
obtaining these results is played by the notion of a
transfer operator which we proceed to define:

\begin{defn}\bold{Transfer Operator}
Let $\seq{p}{a}{\alpha}$ be a  Roman graded sequence.  The
{\em transfer operator} associated with the graded sequence
is the continuous linear operator 
$\tau_{\SB{p}}:{\cal I}\rightarrow {\cal I}$ defined as
$$ \tau_{\SB{p}} \lambda _{a}^{\alpha}(x) = p_{a}^{\alpha}(x)$$
for all $a$ and $\alpha $. 
\end{defn}

Thus, $\tau _{\SB{p}}q_{a}^{\alpha }(x)=q_{a}^{\alpha }({\bf p}).$
\begin{description}
\item[Discrete] If $\seq{p}{a}{\alpha}$
is the associated graded sequence for the delta operator $f(\D )$, we 
frequently write $\tau _{f}$ for the transfer operator associated with
$\seq{p}{a}{\alpha}$.
\item[Continuous] If $\seq{p}{a}{\alpha}$ is the $n\th $ associated graded
sequence for the delta operator $f(\D )$, we  
frequently write $\tau _{f;n}$ for the transfer operator associated with
$\seq{p}{a}{\alpha}$. We also write $\tau_{f}$ for $\tau _{f;0}.$ 
\end{description}

\begin{prop}
\label{wdto} The transfer operator is a well defined regular
operator. 
\end{prop}

\proof{Proposition~\ref{wduc}.}

\ssect{Adjoints}

We must now define the adjoint of a regular operator. This definition is
equivalent to the usual adjoint with respect to the inner product\index{Inner
Product} defined in \S\ref{orth}. 
\begin{defn}
\bold{Adjoint}\label{adjoint}
If $\theta $ is a continuous linear operator on 
${\cal I}^{\alpha}$ for some $\alpha \neq  (0)$; 
then  the {\em adjoint}
of $\theta $ is defined to  be the linear operator $\adj{\theta}$ on Artinian
operators defined by
\begin{equation}\label{adjeq}
\action{\alpha} { {\left[\adj{\theta}f(\D )\right]p(x)} }=  \action{\alpha}
{{f(\D )\theta p(x)}} 
\end{equation}
for all formal power series of logarithmic type
$p(x)\in {\cal I}^{\alpha}$, 
and for all Artinian  operators $f(\D )\in \Lambda^{+}$.

\discrete{In the discrete case, if $\theta $ is a continuous, linear operator
on ${\cal I}^{(0)}$,  then the {\em adjoint}
of $\theta $ is defined to be the linear operator
$\adj{\theta}:\Lambda \rightarrow \Lambda$ such that 
\eref{adjeq} holds for all formal power series of logarithmic type
$p(x)\in {\cal I}^{(0)}$, and for all  differential operators $f(\D )\in
\Lambda$.}
\end{defn}

\begin{prop}
\label{wdadj}
Let $\theta $ be as in Definition~\ref{adjoint}. Then
$\adj{\theta}$ is a well defined continuous linear operator.
\end{prop}

\proof{{\bf (Well Defined)} Spanning argument (Proposition~\ref{spanning}).

{\bf (Linear)} Consider the following string of equalities
\begin{eqnarray*}
\action{\alpha} { {\adj{a_{1}\theta _{1} +
a_{2}\theta _{2}} f(\D )p(x)}}&=& \action{\alpha} {f(\D )(a_{1}\theta
_{1}+a_{2}\theta _{2})p(x)}\\* 
&=& a_{1}\action{\alpha} {f(\D )\theta _{1}p(x)}+ 
a_{2}\action{\alpha} {f(\D )\theta _{2}p(x)}\\
&=& a_{1}\action{\alpha} {\adj{\theta_{1}} f(\D )p(x)}+ 
a_{2}\action{\alpha} {\adj{\theta_{2}} f(\D ) p(x)}\\*
&=& \action{\alpha} {(a_{1}\adj{\theta_{1}}+
a_{2}\adj{\theta_{2}}) f(\D ) p(x)}  .
\end{eqnarray*}

{\bf (Continuity)} Note that the expression $\action{\alpha} {f(\D )\theta
p(x)}=\action{\alpha} {\adj{\theta} f(\D )p(x)}$ is
continuous in $f(\D )$ and $p(x)$.}

Let $\theta $ be a continuous linear operator on ${\cal I}$
or ${\cal I}^{+}$ which maps each ${\cal I}^{\alpha}$ to
itself. 
The adjoint of the restriction of  $\theta $ to ${\cal
I}^{\alpha}$ is denoted by $\adj{\theta }^{\alpha}$.
If $\theta $ is a regular, continuous, linear operator on
${\cal I}^{+}$ or
all of ${\cal I}$, the adjoint of $\theta $ is  the
unique operator which coincides with $\adj{\theta}^{\alpha}$
 for $\alpha \neq (0)$.

The adjoint of nonregular operators on ${\cal I}$ or ${\cal
I}^{+}$ can be
similarly defined; however, we omit the discussion, since
it would be an unnecessary digression.

Let us consider several important operators, and as an exercise compute their
adjoints. 

\begin{ex}
Let $f(\D )$ be an Artinian operator. Then 
$$\left[\adj{f(\D )}\right](g(\D ))=f(\D )g(\D ).$$ 
In other words, the adjoint of an Artinian operator $f(\D )$ is the operator of
 multiplication by $f(\D)$. 
\end{ex}

The Roman shift is an operator which is crucial to our work in the remainder of
this chapter; moreover, it is not a $\D $-invariant operator, so its adjoint is
of particular interest.

\begin{ex}
\label{ex1} We define a regular,
linear operator $\sigma$ on ${\cal I}$ continuous on each ${\cal I}^{\alpha }$
by requiring that for all $a$, 
$$ \sigma \lambda _{a}^{\alpha}(x) = \left\{\begin{array}{ll}
\lambda _{a+1}^{\alpha}(x)&\mbox{if }a\neq -1\mbox{, and}\\[0.1in]
0&\mbox{if }a=-1
\end{array}\right.$$
We call the operator $\sigma_{}$ the {\em standard Roman
shift}.\index{Roman Shift,Standard}\index{Standard Roman Shift}
It is not a
$\D$-invariant operator. For example, $\D \sigma -\sigma \D =\I \neq 0$.
The standard  Roman shift has as its adjoint the Pincherle derivative
operator $ \adj{\sigma}(\D^{a})=a\D ^{a-1}$ for all integers $a$, since 
$$ \action{\alpha} {b\D ^{b-1}\lambda _{a}^{\alpha}(x)}=
\rn{b}!(1-\delta_{0,b}) \delta _{a+1,b}=
 \action{\alpha} {\D ^{b}\sigma_{} \lambda_{a}^{\alpha}(x)}.$$

\cont{In the continuous case, $\sigma $ is not continuous over all of
${\cal I}^{\alpha }$ since  the limit of $\sigma p(x)$ as $p(x)$ approaches
$1/x$ is not well defined.}
\end{ex}

The standard Roman shift is a logarithmic
generalization of the operator {\bf x} of multiplication by $x$, since
$\sigma_{}x^{a}=x^{a+1}$ for $a\neq -1$. 
See Theorem~\ref{esig} for an amazing property obeyed by the Roman shift.

The Pincherle derivative has several equivalent definitions.
\begin{defn}\bold{Pincherle Derivative}\label{pinch}
Let $f(\D )$ be an Artinian operator. We can define its {\em  Pincherle
derivative} $f'(\D )$ by any of the following equivalent formulations:

\begin{enumerate}
\item The Pincherle derivative is the continuous, linear map
$$ \begin{array}{rrcl}
': & \Lambda^{+} &\rightarrow & \Lambda^{+} \\*
&\D ^{a}&\mapsto &a\D ^{a-1}.
\end{array} $$
\item  $f'(\D )=f(\D )\sigma -\sigma f(\D )$.
\item The Pincherle derivative is the adjoint of the
standard Roman shift $\sigma_{}$. In other words,
 $f'(\D )=\adj{\sigma}f(\D)$.
\item The Pincherle derivative of an Artinian operator is its derivative as an
Artinian series. (See the section on derivatives in \cite{ch1}, and
Theorem~\ref{sio1} here.)
\end{enumerate}
\end{defn}

Before considering the adjoint of an arbitrary transfer operator. Let us
consider the following transfer operator which illustrates an interesting quirk
of the continuous theory; it clarifies the meaning of the ``$;n$'' in the
definition of the exponentiation of Artinian series, and their composition.
\explus{C}{\begin{ex}\label{psi}
Let $\psi_{n}$ be the transfer operator
$\tau _{\smallD ;n}$. Thus, $\psi_{n}$ is the continuous
linear operator 
$$ \psi_{n}\lambda_{a}^{\alpha}(x) = e^{2\pi
ian}\lambda_{a}^{\alpha}(x) .$$  
$\psi_{n}$ is $\D $-invariant; however, 
$\psi_{n}$ is not $\D $-invariant for $n\neq 0$. For example,
$\psi_{1}\D^{1/2}=-\D ^{1/2}\psi_{1}$. On the other hand, 
$\psi_{n}\psi_{m}=\psi_{n+m}=\psi_{m}\psi_{n}$.

Now, $\adj{\psi_{n}}f(\D )=f(\D ;n)$ where for $f(x )=\sum _{a}c_{a}x^{a}$,
$f(g;n)=\sum _{a}c_{a}g(x)^{a;n}$. 
\end{ex}}

Finally, we consider the adjoint of an arbitrary transfer operator.
\begin{prop}\label{trans-bij}
If $\tau $ is a  transfer operator, then its adjoint $\adj{\tau}$ is an
automorphism of $\Lambda^{+} $.

\discrete{Moreover, $\adj{\tau }^{(0)}$ is an
automorphism of $\Lambda$.}
\end{prop}

\proof{We show that $\adj{\tau}$ is an monomorphism
and that it acts on delta operators. This  implies that
$\adj{\tau }$ preserves degree, and thus, that
$\adj{\tau} $ is an automorphism.

Let $\tau $ be associated with the  Roman graded sequence
$\seq{p}{a}{\alpha}$.

{\bf (Injectivity)} Assume that $\adj{\tau}f(\D )=\adj{\tau}g(\D )$
for some pair of Laurent operators $f(\D )$ and $g(\D ).$
Thus, for all $p(x)\in {{\cal I}}^{(1)}$,
$$\action{(1)} {f(\D )\tau p(x)}=\action{(1)} {g(\D )\tau
p(x)},$$
so by the spanning argument (Proposition~\ref{spanning}), we infer $f(\D)=g(\D
)$. 

{\bf (Morphism)} 
By Proposition~\ref{wdadj}, $\adj{\tau}$ is continuous and
linear, so we 
need only confirm that $\adj{\tau}$ preserves
multiplication. 
Let  $f(\D )$ and $g(\D )$ be
Artinian operators. Then we have for all a,
\begin{eqnarray*} 
\action{(1)} {\adj{\tau}(f(\D )g(\D ))\lambda_{a}^{(1)}(x)}
&=&\action{(1)} {f(\D )g(\D )\tau \lambda_{a}^{(1)}(x)}\\*
&=&\action{(1)} {f(\D )g(\D )p _{a}^{(1)}(x)}\\ 
&=&\sum_{b}\romcoeff{a}{b}\action{(1)} {f(\D )p_{b}^{(1)}(x)}
\action{(1)} {g(\D )p_{a-b}^{(1)}(x)}\\ 
&=&\sum_{b}\romcoeff{a}{b} 
\action{(1)} {f(\D )\tau \lambda_{b}^{(1)}(x)} 
\action{(1)} {g(\D )\tau \lambda_{a-b}^{(1)}(x)}\\   
&=&\sum _{b}\romcoeff{a}{b}\action{(1)}
{\left[\adj{\tau}f(\D ) \right] \lambda_{b}^{(1)}(x)}
\action{(1)}{\left[\adj{\tau}g(\D )\right] 
\lambda _{a-b}^{(1)}(x)}\\*
&=&\action{(1)} {\left[\adj{\tau}f(\D )\right]
\left[\adj{\tau}g(\D )\right] \lambda_{a}^{(1)}(x)} 
\end{eqnarray*}
Thus, $\adj{\tau}\left[f(\D )g(\D )\right]=
(\adj{\tau}f(\D ))(\adj{\tau}g(\D ))$ by the spanning argument.

{\bf (Degree)} We need only consider the continuous case. Suppose
$\seq{p}{a}{\alpha}$ is the 
($n\th $) associated graded sequence for the delta operator
$f(\D )$, then for all  $a$,
\begin{eqnarray*}
\action{(1)} {\adj{\tau}f(\D )^{b;n}\lambda _{a}^{(1)}(x)}
&=&\action{(1)} {f(\D )^{b;n}p_{a}^{(1)}(x)}\\*
&=&\rn{a}!\delta _{ab}\\*
&=&\action{(1)} {\D^{b} \lambda _{a}^{(1)}(x)} 
\end{eqnarray*}
and by the spanning argument (Proposition~\ref{spanning}),
$\adj{\tau}f(\D)=\D$. The conclusion now follows by the Expansion
Theorem~\ref{exp1}.}   

The most important properties of transfer operators are
stated in the following proposition.

\begin{prop}\label{transfer}
\begin{enumerate}
\item A transfer operator maps Roman
graded sequences to Roman graded sequences.
\item If $\tau : p_{a}^{\alpha}(x) \mapsto
q_{a}^{\alpha}(x)$ is a continuous 
linear operator, where the $\seq{p}{a}{\alpha}$ and
$\seq{q}{a}{\alpha}$ are the ($n\th $ and $m\th $)
associated graded 
sequences for the delta operators $f(\D )$ and $g(\D )$
respectively, then we have  $\adj{\tau}g(\D )=f(\D)$.
\item \begin{description}
\item[Continuous] Moreover, $\adj{\tau }g(\D )^{a;n}=f(\D )^{b;n}$ for all
$a$, and if $\seq{p}{a}{\alpha}$ and $\seq{q}{a}{\alpha}$ are
standard then the operator $\tau $ above is a transfer
operator.
\item[Discrete] The operator $\tau $ above is a transfer
operator. 
\end{description}
\end{enumerate}
\end{prop} 

\proof{We need only consider the continuous case.

{\bf (1)} Let $\tau :\lambda
_{a}^{\alpha}(x) \mapsto p_{a}^{\alpha}(x)$ be a transfer
operator. By 
Theorem~\ref{trans-bij}, $\adj{\tau}$ is an isomorphism of
$\Lambda^{+} $. 
Let $\seq{q}{a}{\alpha}$ be the $m\th $ associated graded
sequence 
for the delta operator 
$g(\D )$. Then,
$$\action{\alpha} {\adj{\tau}^{-1}g(\D )^{b;m}\tau
q_{n}^{\alpha}(x)} =
 \action{\alpha} {g(\D)^{b;m}q_{a}^{\alpha}(x)} =
 \rn{a}!\delta_{ab}.$$ 
By  \cite{ch1} and 
Proposition~\ref{trans-bij}, 
$\adj{\tau}^{-1}g(\D)^{b;m}=(\adj{\tau }^{-1}g(\D ))^{b;n}$  
for some $n$.

Hence, $\tau q_{n}^{\alpha}(x)$ is {\em an} associated
graded sequence for the delta operator $\adj{\tau}^{-1}g(\D).$ 

{\bf (2)} We have the following sequence of
equalities:
\begin{eqnarray*} 
\action{(1)} {\adj{\tau}g(\D )p_{a}^{(1)}(x)}
&=&\action{(1)} {g(\D )\tau p_{a}^{(1)}(x)}\\* 
&=&\action{(1)} {g(\D)q_{a}^{(1)}(x)}\\ 
&=&\delta _{a,1}\\*
&=&\action{(1)} {f(\D)p_{a}^{(1)}(x)}.
\end{eqnarray*} 
Hence, by the spanning argument,
$\adj{\tau}g(\D )=f(\D )$.

{\bf (3)} Similar to {\bf (2)}.

{\bf (4)} More generally, for  an arbitrary Artinian operator
$\sum _{b}c_{b}g(\D )^{b;n}$, we have by
Proposition~\ref{trans-bij} and \cite{ch1},
$\adj{\tau}\sum_{b} c_{b}g(\D)^{b;n}=\sum c_{b}g(\D )^{b;m}$
for some $m$.
Specialize to the case of $\sum_{b}
c_{b}\D ^{b}=f^{(-1)}(\D )^{b;n}$ to get
$$\adj{\tau}\left(f^{(-1)}(g)\right)^{b;n} = e^{2\pi ibm}\D ^{b},$$ 
and hence  
$$\action{\alpha} {\left(f^{(-1)}(g)\right)^{b;n} \tau 
\lambda _{a}^{\alpha}(x)} = \action{\alpha}
{e^{2\pi ian}\D ^{b}\lambda _{a}^{\alpha}(x)}.$$ 
In other words, $\tau e^{2\pi iam}\lambda_{n}^{\alpha}(x)$ 
is associated with the delta operator
$f^{(-1)}(g).$ Finally, note that $m=0$ when all of the
sequences in question are standard.}

The following results illustrate the applications of the
preceding proposition.

\begin{prop}\label{Umbral}
If $\seq{p}{a}{\alpha}$ and 
$\seq{q}{a}{\alpha}$ are the standard associated graded
sequences of the delta 
operators  $g(\D )$ and $f(\D )$ respectively, then the
composition $g(f)$ (resp. $g(f;0)$) is the delta operator with standard
associated graded sequence $q_{a}^{\alpha}({\bf p}).$ 
\end{prop}

\proof{If $\tau : \lambda _{a}^{\alpha}(x)\mapsto
p_{a}^{\alpha}(x)$ 
is a transfer operator, then as in the proof of Part {\bf (1)} of
Proposition~\ref{transfer}, $\tau q_{n}^{\alpha}(x)$ is the
standard associated graded sequence for
$(\adj{\tau})^{-1}g(\D )$. However, $\tau 
q_{a}^{\alpha}(x)=\comp{q_{a}^{\alpha}}{p}$ by
Definition~\ref{compos}. 
Moreover, Part {\bf  (2)} of Proposition~\ref{transfer} asserts that
$\adj{\tau}f(\D )=\D $, so 
$\adj{\tau}^{-1}\D =f(\D )$, and, thus, $\adj{\tau}^{-1}g(\D)
=g(f)$ (resp. $g(f;0)$). The conclusion follows.}

\begin{cor}
The set of Roman graded sequences is closed under
umbral composition. Similarly, the set of standard
Roman graded sequences is closed under umbral
composition.$\Box $  
\end{cor}

\begin{cor} \label{UmbCor}\bld{Inverses}
Let $\seq{p}{a}{\alpha}$, and $\seq{q}{a}{\alpha}$ be
(standard) Roman graded sequences associated with the (standard) delta
operators $f(\D )$, and $g(\D )$ 
respectively. Then the following statements are equivalent:
\begin{enumerate}
\item $f(g)=\D$ (resp.~$f(g;0)=\D $),
\item $g(f)= \D $ (resp. $f(g;0)=\D $),
\item $q_{n}^{\alpha}({\bf p})=\lambda _{n}^{\alpha}(x),$ and
\item $p_{n}^{\alpha}({\bf q})=\lambda _{n}^{\alpha}(x).\Box $
\end{enumerate}
\end{cor}

\begin{cor}\label{conj-ass}
The standard associated Roman graded sequence for the standard
operator $f(\D)$ is the standard conjugate   
Roman graded sequence for the delta operator $f^{(-1;0)}(\D )$, and visa versa.
\end{cor}

\proof{Theorem~\ref{GLTF} and Corollary~\ref{UmbCor}.}

\ssect{The Roman Shift}

The following definition generalizes to all
Roman graded sequences the notion of the standard Roman
shift (Example~\ref{ex1}).

\begin{defn}\bold{Roman Shift}\label{shiftdef}
If $\seq{p}{a}{\alpha}$ is a Roman graded sequence, the {\em roman shift}
relative to $\seq{p}{a}{\alpha}$ is the linear operator 
$\sigma_{\SB{p}} : {\cal I}\rightarrow {\cal I}$ defined by 
\begin{equation}\label{shifteq}
\sigma_{\SB{p}} p_{a}^{\alpha}(x) = 
\left\{\begin{array}{ll}
p_{a+1}^{\alpha}(x)&\mbox{if }a\neq -1\mbox{, and}\\[0.1in]
0&\mbox{if }a=-1.
\end{array}\right.
\end{equation}
for all $a$ and $\alpha $ where 
$\sigma _{\SB{p}}$ is continuous over each ${\cal I}^{\alpha }$ but not
continuous simultaneously over all of ${\cal I}$.
\end{defn}

As with Transfer operators,
if $\seq{p}{a}{\alpha}$ is ($n$-)associated with the delta operator
$f(\D )$, we also  write $\sigma_{f}$ (or $\sigma _{f;n}$) instead of
$\sigma_{\SB{p}} $. 
When no graded sequence of delta operator has been specified, we
assume $\sigma=\sigma_{\smallD ;0 }= \sigma_{\mbox{\boldmath
\scriptsize$\lambda$}} $, that  
is, that $\sigma_{}$ is the standard Roman shift as
previously defined.

\begin{defn}\label{ArtDer}
\bold{Artinian Derivation}
Let $Q$ be an operator on the Artinian (resp. Noetherian) algebra such that 
$$ Qf(x)^{a;n}g(x)^{b;m} = a(Qf(x))(f(x)^{a-1;n}g(x)^{b;m}
+b(Qg(x))f(x)^{a;n}g(x)^{b-1;m}.$$ 
\end{defn}

Note that an Artinian (resp. Noetherian) derivation is automatically a
derivation, and thus is linear.

\begin{lem}
The derivative is an Artinian (resp. Noetherian) derivation.$\Box$
\end{lem}

\label{probsec}
\begin{prop}\label{problem}
\begin{enumerate}
\item A regular, continuous, linear operator $\theta $ defined on
${\cal I}^{+}$ is
a Roman shift if and
only if $\adj{\theta}$ is a continuous, everywhere defined (Artinian)
derivation of the algebra
of Artinian  operators $\Lambda^{+} $ for which $\adj{\theta f(\D )} =1$ for
some delta operator $f(\D)$. 
\item A regular, continuous, linear operator
$\theta $ defined on the logarithmic algebra ${\cal I}$ is
a Roman shift if and
only if $\adj{\theta}$ is a continuous, everywhere defined, (Artinian)
derivation of the algebra of differential operators $\Lambda$
for which $\adj{\theta f(\D )} =1$ for some delta operator
$f(\D)$. 
\end{enumerate}
\end{prop}

\proof{We need only consider the continuous case.

{\bf (Only if)} Let $\seq{p}{a}{\alpha}$ be a Roman
graded sequence of formal power series of logarithmic type. 
By Proposition~\ref{spl-conj}, 
$\seq{p}{a}{\alpha}$ is ($n-$)associated with a delta operator
$f(\D)$. Now, we have  
\begin{eqnarray*}
\action{(1)} {\adj{\sigma _{f}}f(\D )^{b;n}p_{a}^{(1)}(x)}
&=&\action{(1)} {f(\D )^{b;n}\sigma _{f}p_{a}^{(1)}(x)}\\*
&=&\delta _{a,-1}\action{(1)} {f(\D )^{b;n}p_{a+1}^{(1)}(x)}\\
&=&(a+1)\rn{a}!\delta _{a+1,b}\\
&=&b\rn{a}!\delta _{a,b-1}\\*
&=&\action{(1)} {bf(\D )^{b-1;n}p_{a}^{(1)}(x)}
\end{eqnarray*}
Hence, by the spanning argument,  $\adj{\sigma_{f}}
f(\D)^{b;n} = bf(\D )^{b-1;n}.$  Also,   
$\adj{\sigma _{f}}f(\D )=1$. By the continuity of
$\adj{\sigma_{f}}$, the result follows.

{\bf (If)} Conversely, suppose $\adj{\theta}$ is a continuous,
everywhere defined derivation of $\Lambda^{+} $ with
$\adj{\theta}f(\D )=1.$
 Let $\sigma_{\SB{p}} $ be the Roman shift associated with
$\seq{p}{a}{\alpha}$, the $n\th$
associated graded sequence for $f(\D )$. We shall show that
$\theta =\sigma _{\SB{p}}$. 
\begin{eqnarray*}
\action{(1)} {f(\D )^{a;n}\theta p_{b}^{(1)}(x)}
&=&\action{(1)} {\adj{\theta}f(\D )^{a;n} p_{b}^{(1)}(x)}\\*
&=&\action{(1)} {af(\D )^{a-1;n} p_{b}^{(1)}(x)}\\
&=&a\rn{b}!\delta _{a,b-1}\\* 
&=&\action{(1)} {f(\D )^{a;n}\sigma_{\SB{p}} p_{b}^{(1)}(x)}.
\end{eqnarray*}
Thus, by the spanning argument, $\theta =\sigma_{\SB{p}}$.}

Next, we derive the chain rule for Roman shifts.

\begin{prop}\bold{Chain Rule}\label{chain}
Suppose $\sigma _{f}$, and $\sigma
_{g}$ are Roman shift operators. Then
$$ \adj{\sigma _{f}}=\left(\adj{\sigma_{f}}g(\D
)\right)\adj{\sigma _{g}}. $$ 
\end{prop}

\proo{We need only consider the continuous case. 

For any Artinian operator
$h(\D )=\sum _{b}c_{b}g(\D) ^{b;n}$,
\begin{eqnarray*}
\adj{\sigma _{f}}h(\D ) &=&
\sum _{b} bc_{b} g(\D )^{b-1;n} \adj{\sigma_{f}}g(\D ) \\*
&=& \left[\adj{\sigma _{f}}g(\D )\right] \left[\adj{\sigma
_{g}}h(\D )\right],
\end{eqnarray*}
so 
$$ \adj{\sigma_{f}}= \left(\adj{\sigma _{f}}g(\D)\right)
\adj{\sigma_{g}}.\Box$$}  

The following proposition allows us to relate two Roman
shift  operators.

\begin{prop}\label{chaincor}
 If $\sigma _{f}$, and $\sigma
_{g}$ are regular shift operators, then 
$$ \sigma _{f}=\sigma _{g}\adj{\sigma_{f}}g(\D ). $$ 
\end{prop}

\proof{For any Artinian series $h(\D )$ and
logarithmic series $p(x)$, we have
\begin{eqnarray*}
h(\D )\sigma _{f}p(x)&=&\adj{\sigma_{f}}h(\D )p(x)\\*
&=&(\adj{\sigma _{g}}h(\D ))(\adj{\sigma _{f}}g(\D))p(x)
\end{eqnarray*}
by the chain rule. This in turn equals $h(\D )\sigma
_{g}(\adj{\sigma _{f}}g(\D))p(x).$}

\sssapp{Orthogonal Functions}\label{orth}
We  can now define a
symmetric non-degenerate bilinear form associated with any
Roman graded sequence $\seq{p}{a}{\alpha}$. Consider the
principal subsequence $\tilde{p}_{a}(x)$
(Definition~\ref{sfpslt}).
Define an inner product\index{Inner Product}
$$ \action{} {\tilde{p}_{a}(x)|\tilde{p}_{b}(x)}=\delta
_{ab} \rn{a}!.$$  
One verifies that relative to this
bilinear form, the operator $f(\D )$ is adjoint to the Roman
shift $\sigma_{f}$. For example, for the harmonic logarithm
one has $\tilde{\lambda}_{a}(x)=x^{a}$ for all $a$, and  
$$ \action{} {\D x^{a}|x^{b}} = \action{} {x^{a}|\sigma_{}x^{b}} $$
where $\sigma_{}$ is the standard Roman shift, which
restricts to the operator {\bf x} of  multiplication by 
$x$. The principal sequence $\tilde{p}_{a}(x)$ is the set of
eigenfunctions of 
the operator $\sigma_{f}f(\D )$ with eigenvalue $n$. For
example, the  $x^{a}$ form the set
of eigenvalues of the operator ${\bf x}\D $.

The bilinear form
is not  definite, since
$ \action{} {1+x^{-2}|1+x^{-2}}=0 .$

\discrete{However, in the discrete case, we can define a Hermitian form
$$ \action{c} {p(x)|q(x)}=-\action{} {\overline{p(x)}|q(x)},$$ 
so we have
\begin{eqnarray*}
\action{c} {(ix)^{n}|(ix)^{n}} &=& \action{} {\overline{(ix)^{n}}|(ix)^{n}}\\*
&=& \left\{ \begin{array}{ll}
	-\action{} {x^{n}|x^{n}} &\mbox{for $n$ even, and}\\[0.1in]
	\action{} {x^{n}|x^{n}} &\mbox{for $n$ odd,}
\end{array} \right.
\end{eqnarray*}
and thus for $n$ negative
$$ \action{c} {(ix)^{n}|(ix)^{n}}={1}/{(-n-1)!}>0 .$$
Extending by linearity, we obtain a Hermitian inner product
which is positive definite, where the sequence $h_{n}(x)$
defined as 
$$ h_{n}(x)= \left\{ \begin{array}{ll} 
1x^{n} & \mbox{for $n\geq 0$, and}\\[0.1in]
(ix)^{n}&\mbox{for $n<0$}
\end{array} \right.$$
is a complete orthogonal sequence in the Hilbert space
obtained by completion relative to this Hermitian inner
product. We can obtain this sort of sequence using the Knuth coefficients
\cite{ch3} with $\epsilon =-i$ in place of the Roman coefficients. In this
manner,  
$$ \rn{n}!=\left\{ \begin{array}{ll}
n!&\mbox{if $n\geq 0$, and}\\*
(-1)^{n}i/(-n-1)!&\mbox{if $n<0$.}
\end{array} \right. $$
One can therefore develop a spectral theory of the
operator ${\bf x}\D $ in this Hilbert space (rather than
with the smaller one including only positive powers of $x$,
as is done classically).}

\sssapp{Rodrigues' Formula}

Note the following amazing fact.

\begin{thm}\label{esig}
Let $p(x)$ be a finite linear combination of harmonic logarithms
$\sum _{i=1}^{k}c_{i}\lambda_{a_{i}}^{\alpha^{(i)}}(x)$. Then for any complex
number  $z$, we have the following formal identity
$$ (e^{-z\sigma_{}}\D e^{z\sigma_{}})p(x)=(\D - z\I  )p(x). $$
\end{thm}

By $ e^{-a\sigma_{}}\D
e^{a\sigma_{}}p(x)$ we merely mean to indicate the
expression
$$ \sum_{j\geq 0}\sum_{k\geq 0} \frac{(-1)^{j}a^{j+k}}{j!k!}
\sigma_{}^{j}\D \sigma_{}^{k} p(x)$$
which we assert does in fact converge to the indicated value.
The ``operator'' $e^{z\sigma_{}}$ can not be extended to all of
${\cal I}$, since the series $\sum _{k\geq 0}{z^{k}\sigma ^{k}}/{k!}$.

The corresponding classical identity
$$e^{-{\bf x}}\D e^{{\bf x}}=\D +\I  $$ 
is associated with the classical identity
$$ \D {\bf x}-{\bf x}\D =\I  $$
which corresponds to the logarithmic identity
$$ \D \sigma_{}-\sigma_{}\D  = \D '=\I . $$

\proo{By linearity, it suffices to consider the case
$p(x)=\lambda _{a}^{\alpha}(x)$. There are three cases to consider.

{\bf (When $a$ is not a nonpositive integer)} Classical proof can be
applied {\em mutatis mutandis}.

{\bf (When $a=0$)} We calculate:
\begin{eqnarray*}
e^{-z\sigma_{}}\D e^{z\sigma_{}}\lambda _{0}^{\alpha}(x)
&=&  \sum_{j\geq 0}\sum_{k\geq 0} \frac{(-1)^{j}
z^{j+k}}{j!k!} \sigma_{}^{j} \D \sigma_{}^{k}
\lambda_{0}^{\alpha}(x)\\* 
&=& \lambda _{-1}^{\alpha}(x)- 
z\sum_{j\geq 0}\sum_{k\geq 0} \frac{(-1)^{j}
z^{j+k}}{j!(k+1)!} \sigma_{}^{j} \D \sigma_{}^{k}\lambda
_{1}^{\alpha}(x)\\ 
&=& \lambda_{-1}^{\alpha}(x)- 
z\sum_{j\geq 0}\sum_{k\geq 0} \frac{(-1)^{j}z^{j+k}}{j!k!} \sigma_{}^{j}
\lambda _{k}^{\alpha}(x)\\
&=& \lambda _{-1}^{\alpha}(x)
-z\sum_{k\geq 0}\frac{z^{j}}{j!} 
\left(\sum_{k=0}^{j} (-1)^{j}{k\choose j}\right)
\lambda_{k}^{\alpha}(x)\\ 
&=&\lambda _{-1}^{\alpha}(x)- z\lambda _{0}^{\alpha}(x)\\*
&=&(\D -z\I )\lambda _{0}^{\alpha}(x)
\end{eqnarray*}

{\bf (When $a$ is a negative integer)} We similarly compute
\begin{eqnarray*}
e^{-z\sigma_{}} \D e^{z\sigma_{}} \lambda _{a}^{\alpha}(x)
&=&  \sum_{j\geq 0}\sum_{k\geq 0}
\frac{(-1)^{j}a^{j+k}}{j!k!} \sigma_{}^{j} \D \sigma_{}^{k}
\lambda_{a}^{\alpha}(x)\\* 
&=&\sum_{k\geq 0}\sum_{j=0}^{k}
\frac{(-1)^{j}z^{k}} {j!(k-j)!} \sigma_{}^{j} \D
\sigma_{}^{k-j} \lambda_{a}^{\alpha}(x)\\
&=& \sum_{k=0}^{-a-1} \sum_{j=0}^{k}
\frac{(-1)^{j} z^{k}(a+k-j)}{j!(k-j)!}
\lambda_{a+k-1}^{\alpha}(x)\\ 
&=&a\lambda_{a-1}^{\alpha}(x) +z\lambda_{a}^{\alpha}(x)\\*
&=&(\D -z\I )\lambda _{a}^{\alpha}(x).\Box 
\end{eqnarray*}}

\ssect{Explicit Formulas for Roman Graded Sequences}
We are now ready to derive the following recurrence
formula for Roman graded sequences:

\begin{thm}\bold{Recurrence Formula}\label{rodrigues}
If $\seq{p}{a}{\alpha}$ is an associated graded sequence of
delta operator $f(\D)$, then for all $a\neq -1$ and
for all $\alpha $
$$ p_{a+1}^{\alpha}(x)= \sigma (f'(\D ))^{-1}p_{a}^{\alpha}(x)$$
where $\sigma_{}$ is the standard Roman shift.
\end{thm}

\proo{Let $a\neq -1$, then by
Proposition~\ref{chaincor}.
\begin{eqnarray*}
p_{a+1}^{\alpha}(x)&=&\sigma _{f}p_{a}^{\alpha}(x)\\*
&=&\sigma \left( \adj{\sigma} f(\D )\right)^{-1}
p_{a}^{\alpha}(x)\\*  
&=&\sigma f'(\D )^{-1}p_{a}^{\alpha}(x) .\Box 
\end{eqnarray*}}

Next, we give an explicit formula for the associated graded
sequence of a delta operator in terms of the residual series
of the graded sequence of harmonic logarithms.

\begin{prop}\label{TranForm}
If $\seq{p}{a}{\alpha}$ is the ($n\th $) associated graded sequence of
formal power series of logarithmic type
for the delta operator $f(\D )$, then for all $a$ and all
$\alpha \neq  (0)$:
\dc{p_{a}^{\alpha}(x) }{\rn{a}! f'(\D ) f(\D)^{-1-a}
\lambda_{-1}^{\alpha}(x)}{p_{a}^{\alpha}(x) }{\rn{a}! f'(\D ) f(\D)^{-1-a;n}
\lambda_{-1}^{\alpha}(x)}
$$  .$$  
\end{prop}

\proof{We need consider only the continuous case. Let $\seq{q}{a}{\alpha}$ be
the graded sequence defined by:  
$$ q_{a}^{\alpha}(x) = \rn{a}! f'(\D ) f(\D)^{-1-a}
\lambda_{-1}^{\alpha}(x).$$  
It  suffices to verify that $\seq{q}{a}{\alpha}$ is the
$n\th $ basic graded sequence for $f(\D )$. $\seq{q}{a}{\alpha}$ is indeed is a
{\em graded sequence} since  
$\deg \left(f'(\D )f(\D)^{-1-a;n}\right)=-1-a$ and all the
operators in question are continuous. Then, if
$a\neq 0$, we have
\begin{eqnarray*}
\action{\alpha} {q_{a}^{\alpha}(x)}
&=& \action{\alpha} {\rn{a}!f'(\D )f(\D)^{-1-a;n} \lambda
_{-1}^{\alpha}(x)}\\* 
&=& \frac{\rn{a}!}{-a} \action{\alpha} {\left(f(\D )^{-a;n}\right) '
\lambda _{-1}^{\alpha}(x)}
\end{eqnarray*}
Now, $\action{\alpha} {\left(f(\D )^{-a}\right)'
\lambda_{-1}^{\alpha}(x)}= [\D ^{-1}] \left(f(\D
)^{-a}\right)'$. 

The crucial observation in any theorem related to Lagrange inversion is that
for all Artinian operators $g(\D )$, the coefficient
$[\D ^{-1}]g'(\D )$ is zero. Hence, $\action{\alpha}
{q_{a}^{\alpha}(x)}=0$. Thus, Property~\ref{prop2} 
of Definition~\ref{defBas} is satisfied.

Now, consider the  case. $a=0$, and $\alpha \neq  (0)$. Then
we have: 
$$ \action{\alpha} {q_{0}^{\alpha}(x)}=
\action{\alpha} {f'(\D )f(\D )^{-1}\lambda _{-1}^{\alpha}(x)}.$$ 
Say that $c_{b}=[\D ^{b}f(\D )]$.
The coefficient
 $[\D^{-1}] f(\D )^{-1} = c_{1}^{-1}$, and the coefficient 
$[\D ^{0}] f'(\D )=c_{1}$. Neither operator has any terms of
lower degree. Hence, the coefficient  $[\D ^{-1}] (f'(\D)
f(\D)^{-1}) =  1$, and there are no terms of lower degree.
Hence,  
$\action{\alpha} {q_{0}^{\alpha}(x)}=1.$ 
Thus, Property~\ref{prop1} of Definition~\ref{defBas} is satisfied.

Finally,
\begin{eqnarray*}
f(\D )^{b;n}q_{a}^{\alpha}(x)
&=& \rn{a}!f'(\D )f(\D )^{-a}\lambda _{-1}^{\alpha}(x)\\*
&=& \rn{n}q_{a-1}^{\alpha}(x)
\end{eqnarray*}
so Property~\ref{prop3} is also satisfied. Therefore,
$\seq{p}{a}{\alpha}=\seq{q}{a}{\alpha}$.}

Proposition~\ref{TranForm} is unusual
in that it
 does not readily generalize to ${\cal I}^{(0)}$ because
for $\alpha =(0)$,   $\lambda _{-1}^{\alpha}(x)$ equals 0.
Moreover, for $a\geq 0$, $f(\D )^{-1-a}$ (resp. $f(\D )^{-1-a};n$) has negative
degree, and thus is not a member of $\Lambda$. 

Nonetheless, its consequence---the transfer formula---still
holds for ${\cal I}^{(0)}$: 

\begin{thm}\label{rot}
\bold{Transfer Formula}
Let $f(\D )$ be a {\em differential} operator of degree one,
and let $\seq{p}{a}{\alpha}$ be its ($n\th $) associated 
graded sequence. Then for all $a$ and $\alpha $,
\dc{p_{a}^{\alpha}(x)}{f'(\D ) \left(\frac{\D
}{f(\D)}\right)^{a+1} \lambda_{a}^{\alpha}(x)}{p_{a}^{\alpha}(x)}{f'(\D )
\left(\frac{\D }{f(\D)}\right)^{a+1;n} \lambda_{a}^{\alpha}(x) }
In particular, its residual series is given by
$$ p_{-1}^{(1)}(x)=f'(\D ) x^{-1}. $$
\end{thm}

\proo{{\bf (When $\alpha \neq  (0)$)} The conclusion is immediate
from the preceding proposition. 

{\bf (When $\alpha =(0)$)} We need only consider the continuous case. By regularity,
we have the following string of equalities for all $a$:
\begin{eqnarray*}
p_{a}^{(0)}(x) &=& E_{(0),(1)}p_{a}^{(1)}(x)\\*
&=& E_{(0),(1)} f'(\D )\left(\frac{\D }{f(\D )}\right)^{a+1;n}
\lambda_{a}^{(1)}(x)\\ 
&=& f'(\D )\left(\frac{\D }{f(\D )}\right)^{a+1;n}
E_{(0),(1)} \lambda_{a}^{(1)}(x)\\*
&=&f'(\D )\left(\frac{\D }{f(\D )}\right)^{a+1;n}
\lambda_{a}^{(0)}(x). \Box   
\end{eqnarray*}}

The following variants of the transfer formula are
often useful:

\begin{cor}
Let  $\seq{p}{a}{\alpha}$ and
$\seq{q}{a}{\alpha}$ be the ($n\th $ and $m\th $) associated
sequences for the delta operators $f(\D ) $ and $g(\D )$
respectively, then for all $a$ and $\alpha $,
\dc{p_{a}^{\alpha}(x)}{\frac{ f'(\D )g'(\D )^{-1} g(\D)^{a+1}} {f(\D)^{a+1}}
q_{a}^{\alpha}(x)}{p_{a}^{\alpha}(x)}{\frac{ f'(\D )g'(\D )^{-1} g(\D)^{a+1;m}}
{f(\D)^{a+1;n}} q_{a}^{\alpha}(x).\Box }
\end{cor}

\begin{prop}
\label{another}
In the notation of Proposition~\ref{TranForm}, 
\dc{p_{a}^{\alpha}(x)}{g(\D )^{-a} \lambda_{a}^{\alpha}(x) -
\left(g(\D )^{-a} \right)' \lambda_{a-1}^{\alpha}(x)}{p_{a}^{\alpha}(x)}{g(\D
)^{-a;n} \lambda_{a}^{\alpha}(x) - 
\left(g(\D )^{-a;n} \right)' \lambda_{a-1}^{\alpha}(x)}
for $a \neq 0$,
where $g(\D )=f(\D )/\D $.
\end{prop}

\proo{We need only consider the continuous case. By Proposition~\ref{TranForm},
$$p_{a}^{\alpha}(x)=f'(\D  )g(\D )^{-1-a;n}\lambda _{a}^{\alpha}(x).$$
However,
\begin{eqnarray*}
f'(\D )g(\D )^{-1-a}&=& \left( \D g(\D )\right)'g(\D )^{-a-1;n}\\*
&=&\left(\D 'g(\D ) +\D g'(\D )\right)g(\D )^{-a-1;n}\\
&=&g(\D )^{-a} + \D g'(\D )g(\D )^{-a-1;n}\\*
&=&g(\D )^{-a}+\left( g(\D )^{-a;n}\right)'\D /a
\end{eqnarray*}
so that
$$ f'(\D ) g(\D )^{-1-a;n}\lambda _{a}^{\alpha}(x)=
 g(\D )^{-a;n}\lambda _{a}^{\alpha}(x)
 +\left( g(\D )^{-a;n}\right)' \lambda _{a-1}^{\alpha}(x).\Box$$}

Note that  Proposition~\ref{another} does not in general hold for $a=0$.

\begin{cor}\label{yetanother}
In the notation of Proposition~\ref{another}, 
\dc{p_{a}^{\alpha}(x)}{\sigma g(\D ) ^{-a} \lambda
_{a-1}^{\alpha}(x)}{p_{a}^{\alpha}(x)}{\sigma g(\D ) ^{-a;n} \lambda
_{a-1}^{\alpha}(x)}  
for $a \neq 0,1 $, where $\sigma$ is the standard Roman shift.
\end{cor}

\proo{We need only consider the continuous case. By Proposition~\ref{another},
 \begin{eqnarray*}
p_{a}^{\alpha}(x)&=&g(\D )^{-a;n} \lambda_{a}^{\alpha}(x)
-\left(g(\D)^{-a;n}\right)' \lambda_{a-1}^{\alpha}(x)\\*
&=&g(\D )^{-a;n} \lambda_{a}^{\alpha}(x) -g(\D)^{-a;n}
\sigma \lambda_{a-1}^{\alpha}(x) + \sigma g(\D)^{-a;n}
\lambda_{a-1}^{\alpha}(x).\Box  
\end{eqnarray*}}

We conclude with an important remark about graded sequences of
formal power series of logarithmic type.
Let $\seq{p}{a}{\alpha}$ be a Roman graded sequence with coefficients 
$ d_{ab}^{\alpha \beta }=[\lambda_{b}^{\beta}(x)]p_{a}^{\alpha}(x).$
By regularity,
 \begin{enumerate}
\item For $\alpha ,\beta \neq  (0)$ and for all $a$ and $b$,
$d_{ab}^{\alpha\alpha }=d_{ab}^{\beta \beta }$, 
\item For $a$ and $b$ not a negative integer,
$d_{ab}^{(0)(0)}=d_{ab}^{(0)(0)},$ and 
\item For  $\alpha \neq \beta $, and for all $a$ and $b$, 
$ d_{ab}^{\alpha \beta }=0. $
\end{enumerate}
In view of this, we see that in computations with Roman
graded sequences it suffices for most purposes to compute in the
subspace ${\cal I}^{(1)}$. In other words, even in
computations with polynomials it is preferable to deal with
logarithms first! A quick survey of the examples \cref{eggs} demonstrates
the utility of ${\cal I}^{(1)}$ as compared the algebra of polynomials.

\ssect{Composition of Formal Series}\label{cofs}

Let $f(\D )$ be a delta operator. All of the formula for
composition 
and inversions of series (in particular, the various
versions of the LaGrange
inversion formula)\index{LaGrange Inversion} are
consequences of the following theorem. See also \cite{ch1} and \cite{ch5}
for other results concerning the composition of series.
\begin{prop}\label{source}
\begin{enumerate}
\item If $f(\D )$ is  the delta operator with ($n\th $) associated
graded sequence $\seq{p}{a}{\alpha}$, 
then for every Artinian operator $g(\D )$ we have the
following  convergent series:
\dclab{cfseq}{g(f^{-1})}{\sum _{b}\frac{\action{\alpha} {g(\D
)p_{b}^{\alpha}(x)}} 
{\rn{b}!}\D ^{b}}{g(f^{(-1;n)})}{\sum _{b}\frac{\action{\alpha} {g(\D
)p_{b}^{\alpha}(x)}} {\rn{b}!}\D ^{b}}
where $\alpha \neq  (0)$.
\item \eref{cfseq} holds for $\alpha =(0)$ as well whenever $g(\D
)$ is a differential operator.
\end{enumerate}
\end{prop}

\proo{We need only consider the continuous case. By the expansion theorem
(Theorem~\ref{GET}),
$$ g(\D )=\sum_{b} \frac{\action{\alpha} {g(\D)
p_{b}^{\alpha}(x)}} {\rn{b}!} f(\D) ^{b;n} .$$
Hence, substituting $f^{(-1;n)}(\D )$ for $\D $ (by the characterization of
composition in \cite{ch1}),
\begin{eqnarray*}
g(f^{(-1;n)})
&=&\sum_{b} \frac{\action{\alpha}{g(\D) p_{b}^{\alpha}(x)}}
{\rn{b}!} f(f^{(-1;n)}))^{b;n}\\*
&=&\sum _{b} \frac{\action{\alpha}{g(\D)p_{b}^{\alpha}(x)}}
{\rn{b}!} \D^{b}\Box 
\end{eqnarray*}}

Theorem~\ref{source} has many versions and many corollaries.
We first deduce from the Expansion Theorem the following
convergent expansion
\dcthree{\sum _{n}\frac{\action{\alpha}
{g(\D)p_{n}^{\alpha}(x)}}{\rn{n}!}\D ^{n} }{g(f^{(-1)})}{\sum _{n}
\frac{\action{\alpha}
{g(f^{(-1)})\lambda_{n}^{\alpha}(x)}} {\rn{n}!}\D ^{n}}{\sum
_{b}\frac{\action{\alpha} 
{g(\D)p_{b}^{\alpha}(x)}}{\rn{b}!}\D ^{b} }{g(f^{(-1;n)})}{\sum _{b}
\frac{\action{\alpha}
{g(f^{(-1;n)})\lambda_{b}^{\alpha}(x)}} {\rn{b}!}\D ^{b}}
Hence, 
\dc{\action{\alpha} {g(\D )p_{n}^{\alpha}(x)}}{\action{\alpha}
{g(f^{(-1)})\lambda _{n}^{\alpha}(x)}}{\action{\alpha} {g(\D
)p_{b}^{\alpha}(x)}}{\action{\alpha} {g(f^{(-1;n)})\lambda _{b}^{\alpha}(x)}}.
An application of the Transfer Formula (Theorem~\ref{rot}) gives (in, for
example, the continuous case)
\begin{eqnarray*}
\rn{a}! \action{(1)} {g(\D )f'(\D )f(\D )^{-1-a;n} \lambda _{a}^{(1)}(x)}
&=&\action{(1)} {g(\D )p_{a}^{(1)}(x)}\\*
&=&\action{(1)} {g(f^{(-1;n)})\lambda _{a}^{(1)}(x)}.
\end{eqnarray*}
By the spanning argument (Proposition~\ref{spanning}), we have the following
corollary. 

\begin{cor}
Let $f(\D )$ be a delta operator, and let $g(\D )$ be any
Artinian  operator. Then we have the following convergent sum:
\dc{g(f^{(-1)})}{\sum _{k}\action{(1)} {g(\D) f'(\D )
f(\D)^{-1-k} \lambda _{-1}^{(1)}(x)}\D ^{k}}{g(f^{(-1;n)})}{\sum
_{a}\action{(1)} {g(\D) f'(\D ) 
f(\D)^{-1-a;n} \lambda _{-1}^{(1)}(x)}\D ^{a}.\Box}
\end{cor}

By taking $g(\D )=\D ^{a}$, we obtain powers of $f^{(-1;n)}(\D )$ (resp.
$f^{(-1)}(\D )$).

\begin{cor}
If $f(\D )$ is a delta operator with ($n\th $) associated
graded sequence $\seq{p}{a}{\alpha}$, then we have the
following convergent expansion
\dc{ f^{(-1)}(\D )^{a}}{\sum _{k} \frac{\action{(1)} {\D^{k}
p_{k}^{(1)}(x)}}{\rn{k}!} \D ^{k}}{ f^{(-1;n)}(\D )^{a}}{\sum _{b}
\frac{\action{(1)} {\D^{a} p_{b}^{(1)}(x)}}{\rn{b}!} \D ^{b} .\Box}
\end{cor}

\sect{Examples}\label{eggs}
\ssect{Roman Graded Sequences}\index{Roman Graded
Sequences} 

We prefix some general considerations about the computation
of the Roman graded sequence $\seq{p}{a}{\alpha}$ ($n$-)associated with a
delta operator $f(\D )$. The crucial step is the computation
of the residual series $p_{-1}^{(1)}(x)$. This is given by
the simple formula
$$ p_{-1}^{(1)}(x)=f'(\D )(1/x). $$
Once the residual series is known, any of the series
$p_{a}^{(1)}(x)$ can be obtained from the residual series by
applying a suitable power of $f(\D )$; that is, by the
Transfer Formula (Theorem~\ref{rot}),
\dc{p_{k}^{\alpha}(x)}{f'(\D ) \left(\frac{f(\D)
}{\D}\right)^{-k-1} x^{k}}{p_{a}^{\alpha}(x)}{f'(\D )
\left(\frac{f(\D) }{\D }\right)^{-a-1;n} x^{a}}
when $a$ is a negative integer, and
\dc{p_{k}^{\alpha}(x)}{f'(\D ) \left(\frac{f(\D)
}{\D}\right)^{-k-1} x^{k}\left(\log x - 1-\frac{1}{2}-\dots -\frac{1}{k}
\right)}{p_{a}^{\alpha}(x)}{f'(\D ) 
\left(\frac{f(\D) }{\D }\right)^{-a-1;n} x^{a}\left(\log
x-\frac{s(-a,1)}{\rn{-a}!} \right)} 
otherwise.

If we denote the coefficients of $p_{a}^{(1)}(x)$ by
$c_{ab}=[\lambda _{b}^{(1)}(x)]p_{a}^{(1)}(x)$,
then it follows from regularity that 
$$ p_{a}^{\alpha}(x)=\sum_{b}c_{ab} \lambda _{b}^{\alpha}(x)$$ 
for all  $\alpha $. 

Note that for $\alpha =(0)$ we obtain
\dc{p_{k}^{(0)}(x)}{\sum_{j=0}^{b}c_{kj}x^{j} }{p_{aa}^{(0)}(x)}{\sum _{b\in
\SB{R}-\SB{N}^{-}}c_{ab}x^{b}.}
In the discrete case, this is the sequence of polynomials of binomial type
associated with the delta operator $f(\D )$; thus we see
that even in the case of polynomials, it may be speedier to
compute via the logarithmic graded sequence.
\fig{Examples of Roman Graded Sequences}{{rcl|c|rcl|c}
&\thin{Delta}&&Associated&&\thin{Inverse}&&Conjugate\\
&\thin{Operator}&&Graded Sequence&&\thin{Operator}&&Graded Sequence\\
\hline 
&$\D $& &$\lambda _{n}^{\alpha}(x)$&&\thin{$\D $}&&$\lambda _{n}^{\alpha}(x)$\\
$\Delta $&$=$&$E-\I $&$(x)^{\alpha}_{n}$&& \thin{$\log (\I +\D )$}&&
$\phi_{n}^{\alpha}(x)$\\  
$\del$&$=$&$ \I -E^{-1}$&$\lr_{n}^{\alpha}$&& \thin{$-\log (\I -\D )$}&&
$q_{n}^{\alpha}(x)$\\ 
$A_{a}$&$=$&$\D E^{z}$&$A_{n}^{\alpha}(x)$&&&&$\mu_{n}^{\alpha}(x)$\\
&\thin{$E^{z}(E-1)$}&&$G_{n}^{\alpha}(x)$&&&&$g_{n}^{\alpha}(x)$\\
$K$&$=$&$\frac{\smallD }/(\smallD -\smallI)$& $L_{n}^{\alpha}(x)$& $K$&$=$&
$\smallD /(\smallD-\smallI)$& 
$L_{n}^{\alpha}(x)$ }  

\sssect{Lower Factorial}\label{psisec}

Other than the Harmonic logarithms (which we have already
discussed at great length), our first example of a Roman
graded sequence is the logarithmic lower factorial graded
sequence. 

 \begin{defn}\label{FDO}
\bold{Forward Difference Operator} Define the {\em forward
difference operator} $\Delta =E-\I =e^{\smallD }-\I $. Let
$(x)_{a}^{\alpha}$ denote its standard
associated graded sequence; it is
called the {\em logarithmic lower factorial graded sequence}\index{Lower
Factorial}. Let 
$\seq{\phi}{a}{\alpha}$ denote its standard conjugate graded
sequence; it is called the 
{\em logarithmic exponential graded
sequence}.\index{Exponential Graded Sequence}
\end{defn}

\fig{Lower Factorials $(x)_{n}^{\alpha}$}{{rcl|rcl}
$(x)_{2}^{(0)}$&=&$x(x-1)$&
$(x)_{2}^{(1)}$&=&$\lambda _{2}^{(1)}(x)-\lambda _{1}^{(1)}(x)+
B_{3,3}/3x -B_{4,3}/12x^{2}+\cdots $\\
$(x)_{1}^{(0)}$&=&$x$&
$(x)_{1}^{(1)}$&=&$x\log
(x)-x+ B_{2,2}/2x -B_{3,2}/6x^{2}+\cdots $\\
$(x)_{0}^{(0)}$&=&$1$&
$(x)_{0}^{(1)}$&=&$\log(x+1)
+B_{1}/(1+x)- B_{2}/2(1+x)^{2}
+B_{3}/3(1+x)^{3} -\cdots$\\
&&&$(x)_{-1}^{(1)}$&=&$1/(x+1)$\\
&&&$(x)_{-2}^{(1)}$&=&$1/(x+1)(x+2)$}
Now, $\Delta '=E$, so by Theorem~\ref{rot}, we calculate the
residual series
\begin{equation}\label{FD1}
 (x)_{-1}^{(1)}=Ex^{-1}=\frac{1}{x+1}.
\end{equation}
In general,
\begin{prop}\label{FD2}
For $n$ a negative integer, 
$$ (x)_{n}^{(1)}=\frac{1}{(x+1)\cdots (x-n)}=(x)_{n}. $$ 
\end{prop}

\proo{We proceed by induction. The case $n=1$
amounts to \eref{FD1} which we verified above. Now,
suppose proposition holds for $n$ then
\begin{eqnarray*}
 (x)_{-n-1}^{(1)} &=&\rn{-n}^{-1}\Delta (x)_{-n}^{(1)}\\*
&=& \frac{1}{n}\left( \frac{1}{(x+1)\cdots (x+n)} -
\frac{1}{(x+2)\cdots (x+n+1)} \right) \\*
&=&\frac{1}{(x+1)\cdots (x+n+1)}.\Box 
\end{eqnarray*}}

Similarly, we have the continuous analog of \cite[p.\
133]{UC}: 

\begin{prop}\label{FD3}
For $a$ not a negative integer, 
$(x)_{a}^{(0)}=(x)_{a}.$
\end{prop}

\proof{Proposition~\ref{assthm}} 

Thus, 
\begin{cor}
For all $a$, $\tilde{(x)}_{a}=(x)_{a}.\Box$
\end{cor}

Now, we determine  the series $(x)_{0}^{(1)}$.
By Theorem~\ref{rot}, we have
$$ (x)_{0}^{(1)} =E\frac{\D }{\Delta }\log x,$$
or equivalently
$$ \int_{x}^{x+1}(t)_{0}^{(1)}dt = \log (x+1). $$
Moreover, since by the Euler-MacLaurin formula
\begin{equation}\label{FD4}
\frac{\D }{\Delta }=\sum_{k\geq 0}B_{k}\D ^{k}/k!,
\end{equation}
we have
\begin{eqnarray}
(x)_{0}^{(1)}&=&\sum_{k\geq 0}\frac{B_{k}}{k!}\D ^{k}\log
(x+1)\nonumber \\*
&=& \log (x+1) +
\frac{B_{1}}{1+x}-\frac{B_{2}}{2(1+x)^{2}}
+\frac{B_{3}}{3(1+x)^{3}}-\cdots  \label{FD6}
\end{eqnarray}
where the $B_{k}$ are the Bernoulli numbers.
Hence, we discover that $(x)_{0}^{(1)}=\psi(x+1)$ coincides with
the classical $\psi$-function (the logarithmic derivative of
the gamma function)
introduced by Gauss. Similarly, one finds that
$(x)_{1}^{(1)}$ and $(x)_{2}^{(1)}$ coincide with the
digamma and trigamma functions.

The classical expansion
\dc{\left({\D }/{\Delta }\right)^{n} }{ \sum_{k\geq 0}
{B_{kn}} \D ^{k}/k! }{\left({\D }/{\Delta
}\right)^{a;0} }{ \sum_{k\geq 0} {B_{ka}} \D ^{k}/k!}
defines the graded sequence $B_{ab}$ of {\em Bernoulli numbers}\index{Bernoulli
Numbers} of {\em order $b$} and {\em degree $a$}. In terms of these
higher order Bernoulli 
numbers, we obtain for all $a$:
\begin{eqnarray*}
(x)_{a}^{\alpha} &=& \left({\D }/{\Delta
}\right)^{1+a;0} E\lambda _{a}^{\alpha}(x)\\*
&=&\sum_{k\geq 0}B_{k,a+1}\D ^{k}/k!
\lambda_{a}^{\alpha}(x+1)\\* 
&=& \sum_{k\geq 0} B_{k,a+1} \romcoeff{a}{k}
\lambda_{a-k}^{\alpha}(x+1)   
\end{eqnarray*}

\discrete{We now  give an explicit calculation of $(x)_{n}^{(2)}$
for $n$ a negative integer. To begin with, by
Theorem~\ref{rot}, we can calculate the residual series of
order (2).
\begin{eqnarray*}
(x)_{-1}^{(2)}&=& E\left( 2x^{-1}\log x\right)\\*
&=&\frac{2\log (x+1)}{x+1}.
\end{eqnarray*}
Continuing in this way,
\begin{eqnarray*}
(x)_{-2}^{(2)}&=&-\Delta (x)_{-1}^{(2)}\\*
&=&(x)_{-1}^{(2)}-(x+1)_{-1}^{(2)}\\
&=& 2\left(\frac{\log (x+1)}{x+1}-\frac{\log
(x+2)}{x+2}\right)\\*
&=&2 \left( \frac{\log (x+1)-(x+1)\log
\left(\frac{x+2}{x+1}\right)} {(x+1)(x+2)} \right),
\end{eqnarray*}
and 
\begin{eqnarray*}
(x)_{-3}^{(2)}&=&-\frac{1}{2}\Delta (x)_{-2}^{(2)}\\*
&=&\frac{\log (x+1)- (x+1)\log\left(\frac{x+2}{x+1}\right)}
{(x+1)(x+2)} - \frac{\log (x+2) -(x+2)\log
\left(\frac{x+3}{x+2}\right)} {(x+2)(x+3)}\\*
&=&\frac{2\log (x+1)}{(x+1)(x+2) (x+3)}+\frac{\log \left(
\frac{x+3}{x+1}\right)}{(x+3)}. 
\end{eqnarray*}
In general for $n$ a positive integer,
$$(x)_{-n}^{(2)} = \frac{2\log (x+1)}{(x+1)\cdots (x+n)}+\frac{2\log \left(
\frac{x+n}{x+1}\right)}{(n-1)!(x+n)}.$$

Similarly, we calculate $(x)_{n}^{(0,1)}$. The residual
series of order (0,1) is:
\begin{eqnarray*}
(x)_{-1}^{(0,1)} &=& E(1/x\log x)\\*
&=& 1/(x+1)\log (x+1).
\end{eqnarray*}
Next,
\begin{eqnarray*}
(x)_{-2}^{(0,1)} &=& -\Delta (x)_{-1}^{(0,1)}\\*
&=& \frac{(x+1)\log (\frac{x+2}{x+1}+\log
(x+2)}{(x+1)(x+2)\log (x+1)\log (x+2)},
\end{eqnarray*}
and so on.}

The identity 
\begin{equation}\label{FDI}
(x+a)_{0}^{(1)} = \sum _{k\geq 0} \romcoeff{0}{k}
(a)_{k} (x)_{-k}^{(1)}
\end{equation}
gives a classical identity satisfied by the $\psi$-function,
that is:
$$\psi(x+a+1)=\psi(x+1)+\sum _{k\geq
0}\frac{(-1)^{k+1}a(a-1)\cdots (a-k+1)}{k(x+1)(x+2)\cdots (x+k)}.$$
Similar identities can be obtained for the digamma and
trigamma functions.

The logarithmic Taylor's theorem (Theorem~\ref{GLTF}) gives the
following generalization of Newton's expansion:

\begin{prop}\label{Newton}
Every  logarithmic series $p(x)$ can be
uniquely expanded as a convergent  series 
$$ p(x)=\sum_{a,\alpha} {d_{a}^{\alpha}}
(x)_{a}^{\alpha} /{\rn{a}!}$$ 
where the coefficients $d_{a}^{\alpha}$ are given by
$$ d_{a}^{\alpha}=\action{\alpha} {\Delta ^{a;0}p(x)}. \Box $$
\end{prop}

For example, 
\begin{equation}\label{askey}
\frac{1}{x}=\sum_{k\geq 0} {1}/{k!(x+1)\cdots (x+k+1)}. 
\end{equation}

\digress{We {\em  digress} to indicate the meaning of such equalities.
Formally, we merely mean that when both sides of the equality
are expanded in terms of harmonic logarithms $\lambda
_{a}^{\alpha}(x)$ the resulting coefficients are identical.
However, because of such results as Theorem~\ref{Roman2},
we are allowed to make computations in the real or complex
numbers, and thus obtain asymptotic expansions. In 
\eref{askey},  the right side
and left side are both approximately $0.01$ for $x=100$;  in fact, the
error is about $0.0000015$ when you compute 20 or more terms
of the summation. Similarly for $x=69$, when one computes the first 14
terms of the summation, one finds that the left side is 0.14488
and the right side is 0.14493.}

\begin{prob}
How are the Stirling numbers $s(a,b)$ which are the coefficients of the formal
power series $(y)_{a}$  related to the coefficients of the logarithmic series
$(x)_{a}^{\alpha }$?
\end{prob}

\sssect{Upper Factorial}

 \begin{defn}
\bold{Backward Difference Operator}
 Define the {\em backward
difference operator} $\del =\I- E^{-1} =\I  -e^{-\smallD }$. Let
$\lr _{a}^{\alpha}$ denote its standard
associated graded sequence; it is
called the {\em logarithmic upper factorial graded sequence}.\index{Upper
Factorial} 
\end{defn}

\fig{Upper Factorials $\lr_{n}^{\alpha}$}{{rcl|rcl}
$\lr _{2}^{(0)}$&=&$x(x+1)$&
$\lr _{2}^{(1)}$&=&$\lambda _{2}^{(1)}(x)+\lambda _{1}^{(1)}(x)-
B_{3,3}/3x -B_{4,3}/12x^{2}-\cdots $\\
$\lr _{1}^{(0)}$&=&$x$&
$\lr _{1}^{(1)}$&=&$-x\log(x) +x -B_{2,2}/2x
-B_{3,2}/6x^{2} -\cdots $\\
$\lr _{0}^{(0)}$&=&$1$&
$\lr _{0}^{(1)}$&=&$\log (x-1) -
B_{1}/(x-1) +B_{2}/2(x-1)^{2} -B_{3}/3(x-1)^{3} +\cdots$\\ 
&&&$\lr _{-1}^{(1)}$&=&$1/(x-1)$\\
&&&$\lr _{-2}^{(1)}$&=&$1/(x-1)(x-2)$}

As before, the residual series is given by
$$ \lr _{-1}^{(1)}=E^{-1}x^{-1}=\frac{1}{x-1},$$
and for $n$ a negative integer
$$\lr _{n}^{(1)}={1}/{(x-1)\cdots (x+n)}. $$
Similarly, for $a$ not a negative integer,
we have:
$$\lr _{a}^{(0)}= \Gamma (x+a)/\Gamma (x).$$
Thus, for all $a$,
$$ \tilde{\langle x \rangle}_{a} = {\Gamma (x+a)}/{\Gamma (x)} $$
For $a=0$ we have:

\begin{prop}
$\lr _{0}^{(1)}= \log (x-1) -
{B_{1}}/{x-1}+ {B_{2}}/{2(x-1)^{2})} -{B_{3}}/{3(x-1)^{3}}+\cdots . $
\end{prop}

\proo{By Theorem~\ref{rot}, we have
$$ \lr _{0}^{(1)} =E^{-1}\frac{\D }{\del }\log x.$$
From the Euler-MacLaurin formula, namely from
\begin{eqnarray*}
\frac{\D }{\del }&=&\frac{-\D }{e^{-\smallD }-\I }\\*
&=& \frac{\D }{\Delta }(-\D ; 0)\\
&=&\sum_{k\geq 0}\frac{B_{k}}{k!}(-\D) ^{k}\\*
&=&\sum_{k\geq 0}(-1)^{k}\frac{B_{k}}{k!}\D ^{k}.
\end{eqnarray*}
we infer
\begin{eqnarray*}
\lr _{0}^{(1)}&=&\sum_{k\geq 0}(-1)^{k}\frac{B_{k}}{k!}\D ^{k}\log
(x+1)\\*
&=& \log (x-1) -
\frac{B_{1}}{x-1}+\frac{B_{2}}{2(x-1)^{2}}
-\frac{B_{3}}{3(x-1)^{3}}+\cdots . \Box 
\end{eqnarray*}}

We have, in terms of the Bernoulli numbers of higher order:
$$ \left(\frac{\D }{\del}\right)^{a;0}= \sum_{k\geq 0} (-1)^{k}
{B_{ka}} \D^{k}/{k!} .$$
Hence, for all $a$,
\begin{eqnarray*}
\lr _{a}^{\alpha}&=&
\left(\frac{\D }{\del}\right)^{1+a;0} E^{-1} \lambda _{a}^{\alpha}(x)\\*
&=&\sum_{k\geq 0} (-1)^{k}
{B_{k,a+1}} \D^{k} \lambda_{a}^{\alpha}(x-1)/k!\\*
&=& \sum_{k\geq 0} (-1)^{k} B_{k,a+1} \romcoeff{a}{k}
\lambda_{a-k}^{\alpha}(x-1) . 
\end{eqnarray*}

\discrete{As with the upper factorial graded sequence, we compute $\lr
_{n}^{(2)}$ for $n$ a negative integer as follows starting with the residual
series of order (2).
\begin{eqnarray*}
\lr _{-1}^{(2)}&=& E^{-1}\left( 2x^{-1}\log x\right)\\*
&=&\frac{2\log (x-1)}{x-1}\\
\lr _{-2}^{(2)}&=&-\del \lr _{-1}^{(2)}\\
&=& 2\left(\frac{\log (x-2)}{x-2}-\frac{\log
(x-1)}{x-1}\right)\\*
&=&2\frac{\log (x-2)-(x-2)\log \left(\frac{x-1}{x-2}\right)}{(x-1)(x-2)},
\end{eqnarray*}
and in general
$$\lr _{-n}^{(2)}=\frac{2\log (x-n)}{(x-1)\cdots
(x-n)}+\frac{2\log \left(
\frac{x-1}{x-n}\right)}{(n-1)!(x-n)}.$$}

\sssect{Abel}

The logarithmic extension of the Abel polynomials turns out
to be surprisingly pleasing.

\begin{defn}
\bold{Abel Operator}
Define the {\em Abel operator} $A_{z}=\D E^{z}$. 
Its standard associated graded sequence $\seq{A}{a}{\alpha}$
is called the {\em logarithmic Abel graded
sequence}.\index{Abel Sequence}
Its standard conjugate graded sequence $\seq{\mu }{a}{\alpha}$
is called the {\em logarithmic 
inverse-Abel graded sequence}.\index{Inverse-Abel}
\end{defn}

\fig{Abel Graded Sequence $A_{n}^{\alpha}(x)$}{{rcl|rcl}
$A_{2}^{(0)}(x)$&=&$x(x-2z)$&
$A_{2}^{(1)}(x)$&=&$\sigma_{}\lambda _{1}^{(1)}(x-2z)$\\
$A_{1}^{(0)}(x)$&=&$x$&
$A_{1}^{(1)}(x)$&=&$\sigma_{}\log (x-z)$\\
$A_{0}^{(0)}(x)$&=&$1$&
$A_{0}^{(1)}(x)$&=&$\log(x)-{z}/{x}$\\
&&&$A_{-1}^{(1)}(x)$&=&$x(x+z)^{-2}$\\
&&&$A_{-2}^{(1)}(x)$&=&$x(x+2z)^{-3}$}

By Corollary~\ref{yetanother}, for $a\neq
0,1$, we obtain
\begin{eqnarray*}
A_{a}^{\alpha}(x)&=& \sigma E^{-az}\lambda _{a-1}^{\alpha}(x)\\*
&=&\sigma_{}\sum_{k\geq 0}\romcoeff{a-1}{k} (-az)^{k}
\lambda_{a-k-1}^{\alpha}(x)\\*  
&=&\sum_{\scriptstyle k\geq 0 \atop \scriptstyle k\neq a }
\romcoeff{a-1}{k} (-az)^{k} \lambda _{a-k}^{\alpha}(x).
\end{eqnarray*}
In particular, when $a$ is not a negative integer, we obtain
a simple generalization of the classical Abel polynomials.
\dclab{A1}{A_{n}^{(0)}(x)}{x(x-nz)^{n-1}}{A_{a}^{(0)}(x)}{x(x-az)^{a-1;0},}
and for $a$ a negative integer, we  have
\begin{equation}\label{A2}
A_{a}^{(1)}(x)=x(x-az)^{a-1}.
\end{equation}
Thus, its residual series is
$$ A_{-1}^{(1)}=\frac{x}{(x+z)^{2}}, $$
and its principal sequence is 
\dc{\tilde{A}(x)}{x/(x-az)^{a;0}}{\tilde{A}(x)}{x/(x-az)^{a}.}
 $A_{0}^{(1)}(x)$ can be calculated via Theorem~\ref{rot}. It is expressed
 very simply.
\begin{eqnarray}
A_{0}^{(1)}(x)&=&A'_{z}E^{-z}\log x\nonumber \\*
&=&(\I  +z\D )\log x\nonumber \\*
&=&\log x+{z}/{x},\label{A3}
\end{eqnarray}
so 
$$ A_{0}^{\alpha }=\lambda_{0}^{\alpha}(x)+z\lambda_{-1}^{\alpha}(x). $$

From the logarithmic binomial identity, we have the sum
$$A_{0}^{(1)}(x+b)=\sum_{k\geq 0} \romcoeff{0}{k}
A_{k}^{(0)}(b)A_{-k}^{(1)}(x)$$
over all nonnegative integers $k$.
Thus, we infer the remarkable identity
\begin{equation}\label{A4}
\frac{a}{x+b}+\log (x+b)=\frac{a}{x}+\log x+\sum
_{k\geq 1}\frac{(-1)^{k+1}b(b-ak)^{k-1}x}{k(x+ak)^{k+1}}
\end{equation}
For example, we can substitute here the
values $a=1$, $b=2$, and 
$x=5$. If we compute the first 12 terms of the series, the
left hand and right hand sides of \eref{A4} are both approximately 2.0887673.

In general, by Theorem~\ref{rot}
\begin{eqnarray*}
A_{a}^{\alpha}(x)&=& E^{-az}(\I  +z\D )\lambda _{a}^{\alpha}(x)\nonumber \\*
&=& \lambda_{a}^{\alpha}(x-az)
+z\rn{a}\lambda_{a-1}^{\alpha}(x-az).
\end{eqnarray*}
For example,
$$A_{1}^{(1)}(x)=x\log (x-z)+z-x.$$
Again, by Theorem~\ref{GLTF}, every formal power series of
logarithmic type can be expanded in terms of Abel series
$$ p(x)=\sum_{a,\alpha } \frac{d_{a}^{\alpha}} {\rn{a}!}
A_{a}^{\alpha}(x) $$ 
where $$d_{a}^{\alpha}=\action{\alpha} {E^{az}\D ^{a}p(x)}.$$ 
For example, (See \eref{RemId})
\begin{eqnarray}\label{Neid}
 \log x& =& 
\sum_{k\leq 0}(ka)^{-k}\romcoeff{0}{k}A_{k}^{(1)}(x)\\*
&=&A_{0}^{(1)}(x) +z
A_{-1}^{(1)}(x)-2z^{2}A_{-2}^{(1)}(x) +
9z^{3}A_{-3}^{(1)}(x)-\cdots \nonumber \\*
&=&\log x-\frac{z}{x}-\sum_{k>0}
 \left(\frac{k}{x+kz}\right) ^{k+1}z^{k}x. \nonumber
\end{eqnarray}
That is,
$$ x^{-2}= \sum_{k>0}\left(\frac{ka}{x+ka}\right) ^{k+1}.$$

\sssect{Gould}

\begin{defn}\bld{Logarithmic Gould Graded Sequence}\index{Gould Sequence}
We define the {\em logarithmic Gould graded sequence} $\seq{G}{a}{\alpha}$
to be the standard graded sequence associated with the delta
operator $E^{z}\Delta =E^{z+1}-E^{z}$. 
\end{defn}

\fig{Gould Graded Sequence $G_{n}^{\alpha}(x)$}{{rcl}
$G_{2}^{(0)}(x)$&=&$x(x-2a-1)$\\
$G_{1}^{(0)}(x)$&=&$x$\\
$G_{0}^{(0)}(x)$&=&$1$\\
\hline 
$G_{2}^{(1)}(x)$&=&$\lr _{2}^{(1)}-a\lr _{1}^{(1)}+
 \sum_{ k\geq 2} \frac{(-1)^{n+k}{na\choose k}{\rn{n-k-1}!}}
{{\rn{n-1}!(x-1)\cdots (x-n+k)}}$\\[0.1in]
$G_{1}^{(1)}(x)$&=&$x\log
(x)-x+ \sum_{k\geq 2} \left[ B_{k}^{(2)}\romcoeff{2}{k}x^{1-k}
+\frac{\left({a\choose k }+a(k-1){a-1\choose
k-1}\right)(-1)^{k}} {(k-1)!(x-1)\cdots (x-k+1)}\right]$\\[0.1in]
$G_{0}^{(1)}(x)$&=& $\log(x+1)+{B_{1}}/({x+1})- {B_{2}}/{2!(1+x)^{2}}
+\cdots$\\[-0.05in]
&& $+{a}/({x+1}) -{a}/{(x+1)(x+2)}
+{a}/{2!(x+1)(x+2)(x+3)} +\cdots $\\[0.1in]
$G_{-1}^{(1)}(x)$&=&${x}/{(x+a)(x+a+1)}$\\[0.1in]
$G_{-2}^{(1)}(x)$&=&${x}/{(x+2a)(x+2a+1)(x+2a+2)}$}

The Pincherle derivative of $E^{z}\Delta $ is
$(z+1)E^{z+1}-aE^{z}$, so the residual series is given by
\begin{eqnarray*}
 G_{-1}^{(1)}(x)
& =&\left((z+1)E^{z+1} -aE^{z}\right)\left(\frac{1}{x} \right)\\* 
&=& \frac{z+1}{x+z+1}-\frac{z}{x+z}\\*
&=& \frac{x}{(x+z)(x+z+1)}.
\end{eqnarray*}
Since Roman graded sequences are basic, we have
\begin{eqnarray*}
G_{-2}^{(1)}(x)&=& -E^{z}\Delta G_{-1}^{(1)}(x)\\*
&=&E^{z}\left(\frac{x}{(x+z)(x+z+1)}-\frac{x+1}{(x+z+1)(x+z+2)}
\right)\\
&=&\frac{x}{(x+2z)(x+2z+1)(x+2z+2)}\\*
&=&\frac{x}{x+2z}(x+2z)_{-2}.
\end{eqnarray*}
Similarly, by induction we have for $n$ positive
$$ G_{-n}^{(1)}(x)={x}(x+nz-1)_{-n-1}, $$
and, by induction for $n$ nonnegative
$$ G_{n}^{(0)}(x)=x(x-nz-1)_{n-1}. $$
In general for all $a$, 
$$ \tilde{G}_{a}(x)=x(x-az-1)_{a-1}. $$

See \S\ref{UFGsec} for an explicit computation of $G_{a}^{\alpha }(x)$.

\sssect{Laguerre}\label{LagSec}
Our final example of a Roman graded sequence is the logarithmic
Laguerre graded sequence.

\begin{defn}\label{LagOp}
\bold{Laguerre Operator}
Define the {\em Laguerre operator} to be 
$$K={\D }/(\D -\I ).$$
Define the {\em logarithmic Laguerre graded
sequence}\index{Laguerre Sequence} to be its
($0\th $) associated graded sequence denoted
$\seq{L}{a}{\alpha}$. 
\end{defn}
\fig{Laguerre Graded Sequence
$L_{n}^{\alpha}(x)$}{{rcl|rcl} 
$L_{2}^{(0)}(x)$&=&$x^{2}-2x$&
$L_{2}^{(1)}(x)$&=&$\lambda _{2}^{(1)}(x)-2\lambda _{1}^{(1)}(x)$\\
$L_{1}^{(0)}(x)$&=&$-x$&
$L_{1}^{(1)}(x)$&=&$-\lambda _{1}^{(1)}(x)$\\
$L_{0}^{(0)}(x)$&=&$1$&
$L_{0}^{(1)}(x)$&=&$\log(x)-\sum_{k\geq 0}(-1)^{k}k!x^{-k-1}$\\
&&&$L_{-1}^{(1)}(x)$&=&$\sum_{k\geq 1}(-1)^{k}k!x^{-k}$\\
&&&$L_{-2}^{(1)}(x)$&=&$-\sum_{k\geq 2}(-1)^{k}(k-1){k!}x^{-k}/2$} 

\cont{Note that the Laguerre operator $K$ and the Laguerre graded
sequence $L_{a}^{\alpha }(x)$
are not standard, since the leading of $K$ is -1. However, $-K$ and $K(-\D ;0)$
are both standard operators we  discuss later (\S\ref{UF-LF} and
\S\ref{LF-UF}), and if $r_{a}^{\alpha }(x)$ and $p_{a}^{\alpha }(x)$ are their
standard associated sequences, then
\begin{equation}\label{star100}
L_{a}^{\alpha }(x)=\psi_{-1/2}r_{a}^{\alpha }(x)
\end{equation}
where $\psi_{a}$ is as defined in Example~\ref{psi}, and
\begin{equation}\label{dot100}
L_{a}^{\alpha }(x)=(-1)^{\pi ia}p_{a}^{\alpha }(x).
\end{equation}}

Now, the Pincherle derivative of the Laguerre operator is
given by $K'=-(\D -\I  )^{-2}$, so by Theorem~\ref{rot}, 
\begin{eqnarray*}
L_{a}^{\alpha}(x)&=& - (\D - \I  )^{-2}(\D - \I  )^{a+1;0}
\lambda_{a}^{\alpha}(x)\\* 
&=& - (\D -\I  )^{a-1;0} \lambda _{a}^{\alpha}(x).
\end{eqnarray*}
Hence, for all $a$,
\begin{eqnarray*}
L_{a}^{\alpha}(x)&=& e^{\pi ia}
\left( \sum_{k\geq 0}  (-1)^{k}
{a-1\choose k} \D ^{k} \right) \lambda_{a}^{\alpha}(x)\\*
&=& e^{\pi ia} \sum_{k\geq 0} (-1)^{k} {a-1\choose k} \D^{k}
\lambda_{a}^{\alpha}(x)\\* 
&=&e^{\pi ia} \sum_{k\geq 0} (-1)^{k} {a-1\choose k}
\frac{\rn{a}!}{\rn{a-k}!} \lambda_{a-k}^{\alpha}(x).
\end{eqnarray*}
Thus, $$ L_{1}^{\alpha}(x)=-\lambda _{1}^{\alpha}(x). $$
Similarly, we have derived the following expansion
$$ L_{0}^{(1)}(x)=\log
(x)+ \frac{1}{x}- \frac{1}{x^{2}}+ \frac{2}{x^{3}}-
\frac{6}{x^{4}}+ \cdots.$$

\begin{description}
\item[Discrete] Thus, $L_{n}^{\alpha }(x)$ does not contain any terms of
negative degree when $n$ is a positive integer. This is true only for Roman
sequences associated with  ${\D }/(a\D +b\I)$ where $a$ and $b$ are a nonzero
complex numbers. Note that the series of degree 0 contains no negative terms
only if the operator in question is $a\D $ where $a$ is a nonzero complex
number. 
\item[Continuous] In the continuous case, $L_{a}^{\alpha }$ contains terms of
negative degree for all $a$. The only Roman sequences which do not always
contain terms of negative degree are those associated with delta operators of
the form $z\D $ where $z$ is a nonzero complex number.
\end{description}
 
From Theorem~\ref{esig}, we derive a logarithmic extension
of the classical 
Rodrigues' formula for Laguerre polynomials.
$$ L_{a}^{\alpha}(x)=-e^{\sigma_{}}\D
^{a-1}e^{-\sigma_{}}\lambda _{a}^{\alpha}(x). $$

\ssect{Connection Constants}
Given two graded sequences $\seq{p}{a}{\alpha}$ and $\seq{q}{a}{\alpha}$; we
would like to express one in terms of the other.
$$ p_{a}^{\alpha}(x)= \sum_{b} d_{ab} q_{b}^{\alpha}(x). $$ 
The coefficients\index{Coefficients} $d_{ab}$ are called the {\em connection
constants}\index{Connection Constants} from the graded sequence
$\seq{q}{a}{\alpha}$ to the graded sequence
$\seq{p}{a}{\alpha}$, and are denoted $[q_{b}^{\alpha
}(x)]p_{a}^{\alpha }(x)$.

If $\seq{p}{a}{\alpha}$ and $\seq{q}{a}{\alpha}$
are the standard Roman graded sequences associated with the
delta operators $f(\D )$ and $g(\D ) $ respectively, then by
Proposition~\ref{Umbral}
$$ r_{a}^{\alpha}(x)=\sum_{b} c_{ab} \lambda_{b}^{\alpha}(x) $$
is the standard Roman graded sequence associated with $g(f^{(-1)};0)$.
Thus, to determine the connection constants it  suffices merely
to calculate $r_{a}^{\alpha}(x)$. This easy device for the
computation of connection constants is the most effective
application of the present theory.

\sssect{Upper Factorial to Lower Factorial}\index{Connection
Constants,Upper Factorial to Lower Factorial}\index{Upper
Factorial}\index{Lower Factorial}\label{UF-LF}
To express $(x)_{a}^{\alpha}$ in terms of $\lr _{a}^{\alpha}$ we
first calculate
\begin{eqnarray*}
\Delta (\del^{(-1)};0)&=& e^{-\log (\I -\D)}-\I \\*
&=&\frac{\I }{\I -\D }-\I \\
&=&\frac{\D }{\I  -\D }\\*
&=&-K(\D )
\end{eqnarray*}
where $K(\D )$ is the Laguerre operator
(Definition~\ref{LagOp}). Thus, $r_{a}^{\alpha}(x)$ is the standard graded
sequence related to the logarithmic Laguerre graded sequence defined by
\eref{star100},
$$ r_{a}^{\alpha}(x)= \Psi_{1/2}L_{a}^{\alpha}(x) $$
where $\Psi_{1/2}$ is as defined in Example~\ref{psi}
$$ \psi_{1/2}\lambda_{a}^{\alpha}(x)=e^{\pi ia}\lambda_{a}^{\alpha}(x). $$
In a sense, $\psi_{1/2}$ is the logarithmic analog of
substitution of $-x$ for $x$.

We now apply the results from \S\ref{LagSec}. 
\fig{Lower Factorial in Terms of
Upper Factorial}{{rcl|rcl}
$(x)_{2}^{(0)}$&=&$\lr _{2}^{(0)}+2\lr _{1}^{(0)}$&
$(x)_{2}^{(1)}$&=&$\lr _{2}^{(1)}(x)+2\lr _{1}^{(1)}$\\
$(x)_{1}^{(0)}$&=&$\lr _{1}^{(0)}$&
$(x)_{1}^{(1)}$&=&$\lr _{1}^{(1)}$\\
$(x)_{0}^{(0)}$&=&$\lr _{0}^{(0)}$&
$(x)_{0}^{(1)}$&=&$\lr_{0}^{(1)}
-\sum_{k\geq 0} k!\lr _{-k-1}^{(1)}$\\
&&&$(x)_{-1}^{(1)}$&=&
$-\sum_{k\geq 1}k!\lr  _{-k}^{(1)}$\\
&&&$(x)_{-2}^{(1)}$&=&$ \sum_{k\geq 2}(k-1)k! \lr_{-k}^{(1)}/2$}
Thus, for all $a$ and $\alpha $
$$ (x)_{a}^{\alpha } = \sum_{k\geq 0} {a-1\choose k}
\frac{\rn{a}!}{\rn{a-k}!} \lr_{a-k}^{\alpha}(x). $$

From Theorem~\ref{esig}, we have the identity
\dc{(x)_{a}^{\alpha}}{-e^{\sigma_{\del}}\del^{a-1} e^{-\sigma_{\del}}\lr
_{a}^{\alpha}. }{(x)_{a}^{\alpha}}{-e^{\sigma_{\del}}\del^{a-1;0}
e^{-\sigma_{\del}}\lr _{a}^{\alpha}.}  

\sssect{Lower Factorial to Upper Factorial}\index{Connection
Constants,Lower Factorial to Upper Factorial}\index{Upper
Factorial}\index{Lower Factorial}\label{LF-UF}

Conversely, to compute $\lr_{a}^{\alpha}$ in terms of
$(x)_{a}^{\alpha}$, we need the delta operator 
$$ \del (\Delta ^{(-1;0)})=\frac{\D }{\D +\I }=K(-\D ). $$
Hence, the connection constants from
$(x)_{a}^{\alpha}$ to $\lr_{a}^{\alpha}$
are given by the coefficients of
$e^{-\pi ia} L_{a}^{\alpha}(x)$ (\eref{dot100}). Hence, for all $a$ and
$\alpha$, 
$$ \lr_{a}^{\alpha } = \sum_{k\geq 0} (-1)^{k} {a-1\choose k}
\frac{\rn{a}!}{\rn{a-k}!} \lambda_{a-k}^{\alpha}(x). $$

\fig{Upper Factorial in Terms of
Lower Factorial}{{rcl|rcl}
$\lr_{2}^{(0)}$&=&$(x) _{2}^{(0)}-2(x) _{1}^{(0)}$&
$\lr_{2}^{(1)}$&=&$(x) _{2}^{(1)}\lr-2(x) _{1}^{(1)}$\\
$\lr_{1}^{(0)}$&=&$(x) _{1}^{(0)}$&
$\lr_{1}^{(1)}$&=&$(x) _{1}^{(1)}$\\
$\lr_{0}^{(0)}$&=&$(x) _{0}^{(0)}$&
$\lr_{0}^{(1)}$&=&$(x)_{0}^{(1)}
-\sum_{k\geq 0}(-1)^{k}k!\lr _{-k-1}^{(1)}$\\
&&&$\lr_{-1}^{(1)}$&=&
$\sum_{k\geq 1} (-1)^{k+1}k!(x)_{-k}^{(1)}$\\
&&&$\lr_{-2}^{(1)}$&=&$ \sum_{k\geq 2} (-1)^{k}
(k-1)k!(x)_{-k}^{(1)}/2$}  

\sssect{Laguerre to Harmonic}\index{Connection
Constants,Laguerre to Harmonic}\index{Laguerre
Sequence}\index{Harmonic Logarithm}
If $\seq{p}{a}{\alpha}$ is the standard Roman graded
sequence  associated with the delta operator $f(\D )$,
then finding the connection constants from $\seq{p}{a}{\alpha}$
to $\seq{\lambda}{a}{\alpha}$ is tantamount to finding the
standard graded sequence associated with $f^{(-1;0)}(\D )$; that is, the
standard conjugate graded sequence for $f(\D )$.

Note that since $x/(x-1)$ is the 0-compositional inverse
of itself, the logarithmic Laguerre graded sequence is
self-conjugate, so we have for all $a$ and $\alpha $
$$ \lambda_{a}^{\alpha}(x) = \sum_{k\geq 0} (-1)^{k}
{a-1\choose k} \frac{\rn{a}!}{\rn{a-k}!}
L_{a-k}^{\alpha}(x)$$
and from Theorem~\ref{esig}, 
$$ \lambda _{a}^{\alpha}(x)= L_{a}^{\alpha}({\bf L})=
-e^{\sigma_{\SB{L}}}K^{a-1;0} e^{-\sigma_{\SB{L}}}
L_{a}^{\alpha}(x). $$

\sssect{Lower Factorial to Harmonic}\index{Connection
Constants,Lower Factorial to Harmonic}\index{Harmonic Logarithm}
\index{Lower Factorial}
Similarly, the connection constants from
$(x)_{a}^{\alpha}$ to $\seq{\lambda}{a}{\alpha}$ are
given by $\seq{\phi}{a}{\alpha}$---the logarithmic exponential
graded sequence---which we  now compute.

The relevant delta operator is 
$$\log (\I +\D
)=\sum_{k>0}(-1)^{k+1} {\D ^{k}}/{k}.$$
Thus, by Theorem~\ref{rot}, we can calculate the residual series.
\begin{eqnarray*}
\phi _{-1}^{\alpha}(x) &=& \left(\frac{\I }{\I  +\D} \right)
\lambda_{-1}^{\alpha}(x)\\*  
&=&\sum_{k\geq 0}(-1)^{k}\D ^{k}\lambda _{-1}^{\alpha}(x)\\
&=& \sum_{k\geq 0} \frac{(-1)^{k}}{\rn{-k-1}!} \lambda_{-1-k}^{\alpha}(x)\\*
&=& \sum_{k\geq 0}k!\lambda _{-1-k}^{\alpha}(x).
\end{eqnarray*}
Finally, by Theorems~\ref{rodrigues} and \ref{esig}, for $a\neq -1$,
we obtain the 
recursion formula
$$ \phi_{a+1}^{\alpha}(x) = \sigma_{}(\I  + \D)
\phi_{n}^{\alpha}(x) =\sigma_{}e^{-\sigma_{}}\D
e^{\sigma_{}} \phi_{n}^{\alpha}(x). 
$$

\sssect{Abel to Harmonic}\index{Connection
Constants,Abel to Harmonic}\index{Harmonic
Logarithm}\index{Abel Sequence}

The compositional inverse  of the Abel
operator is not easily calculated. Nevertheless, we may calculate the
logarithmic inverse-Abel graded sequence\index{Inverse-Abel}, using
the theory of conjugate graded sequences. 

Thus, by Definition~\ref{defConj}
$$ \mu_{a}^{\alpha}(x) = \sum_{b} \frac{\action{\alpha}
{\D^{b}E^{bz} \lambda _{a}^{\alpha}(x)}}{\rn{b}!}
\lambda_{b}^{\alpha}(x). $$ 
Now,
\begin{eqnarray*}
\D ^{b}E^{bz} \lambda_{a}^{\alpha}(x) &=& \D ^{b}
\sum_{j\geq 0} \romcoeff{a}{j} (bz)^{j}
\lambda_{a-j}^{\alpha}(x)\\* 
&=&\sum_{j\geq 0} (bz)^{j} \frac{\rn{a}!}{j!\rn{a-b-j}!}
\lambda_{a-b-j}^{\alpha}(x)  .
\end{eqnarray*}
Thus,
\begin{eqnarray*}
\mu_{a}^{\alpha}(x) &=& \sum_{k\geq 0} {((a-k)z)^{k}
\frac{\rn{a}!}{k!0!}} \lambda _{a-k}^{\alpha}(x) /{\rn{a-k}!}\\*
&=&\sum_{k\geq 0}((a-k)z)^{k} \romcoeff{a}{k} \lambda_{a-k}^{\alpha}(x).
\end{eqnarray*}
Hence, the remarkable identity
\begin{equation}\label{RemId}
\lambda _{a}^{\alpha}(x)=\sum_{k\geq 0}((a-k)z)^{k}
\romcoeff{a}{k} A_{a-k}^{\alpha}(x) .
\end{equation}

\sssect{Upper Factorial to Gould}\label{UFGsec}\index{Connection
Constants,Upper Factorial to Gould}\index{Upper
Factorial}\index{Gould Sequence}
Assume that Gould parameter  $z$ is a real number which we denote by $t$.
The relevant operator is $f(\D )=-\D (\I - \D)^{-t;0}$ whose
associated graded sequence is:
\begin{eqnarray}
r_{a}^{\alpha}(x)&=& f'(\D )\left(\frac{f(\D )}{\D }\right)^{-a-1;0}
\lambda _{a}^{\alpha}(x)\nonumber \\*
&=&e^{2\pi ia}((t-1)\D +\I )(\I -\D )^{at-1;0}
\lambda_{a}^{\alpha}(x). \label{UF2G} 
\end{eqnarray}
As expected, for $t=1$ we rederive the Laguerre graded sequence $L_{a}^{\alpha
}(x)$, and for $t=0$ we get a variant on the harmonic graded sequence
$\psi_{1/2}\lambda_{a}^{\alpha}(x)$.

For $a\neq 0,1$, instead of \eref{UF2G}, we may use
Corollary~\ref{yetanother}.
\begin{eqnarray*}
r_{0}^{\alpha}(x) &=& e^{\pi ia} \sigma (\I -\D )^{at}
\lambda_{a-1}^{\alpha}(x)\\* 
&=&e^{\pi ia}\sigma_{}\left(\sum_{k\geq 0}(-1)^{k}
{at\choose k}\D ^{k}\right) \lambda_{a-1}^{\alpha}(x)\\*
&=&\sum_{\scriptstyle k\geq 0 \atop \scriptstyle k\neq a}
e^{\pi i(a+k)} {at\choose k }\frac{\rn{a-1}!}{\rn{a-k-1}!}
\lambda _{a-k}^{\alpha}(x). 
\end{eqnarray*}
Thus, for $a\neq 0,1$, we obtain the following remarkable
identity relating the Gould graded sequence to the lower
factorial graded sequence.
$$ G_{a}^{\alpha}(x)= \sum_{\scriptstyle k\geq 0 \atop
\scriptstyle k\neq n} e^{\pi i (a+k)} {at\choose k}
\frac{\rn{a-1}!}{\rn{a-k-1}!} \lr_{a-k}^{\alpha}.$$

For $a=0$, \eref{UF2G} reduces to
\begin{eqnarray*}
r_{0}^{\alpha}(x)&=& ((t-1)\D +\I )(\I  -\D )^{-1}\lambda _{0}^{\alpha}(x)\\*
&=& (1+t\D +t\D ^{2}+t\D ^{3}+\cdots )\lambda _{0}^{\alpha}(x)\\*
&=&\lambda _{0}^{\alpha}(x) + t\sum_{k>0} (-1)^{k+1} (k-1)!
\lambda_{-k}^{\alpha}(x). 
\end{eqnarray*}
Thus, the   Gould series of degree zero are given by the
elegant expression 
$$ G_{0}^{\alpha}(x)=\lr _{0}^{\alpha} + t\sum_{k>0} (-1)^{k+1}
(k-1)! \lr_{-k}^{\alpha}. $$ 
Thus,
\begin{eqnarray*}
{ G_{0}^{(1)}(x) } &= & \log
(x+1)+\frac{B_{1}}{x+1}-\frac{2!B_{2}}{(1+x)^{2}}+\cdots
\phil
+\frac{t}{x+1} -\frac{t}{(x+1)(x+2)}
+\frac{t}{2!(x+1)(x+2)(x+3)} +\cdots. 
\end{eqnarray*}
Finally, from Proposition~\ref{another}, we obtain
\begin{eqnarray*}
r_{1}^{\alpha}(x)&=& -(\I -\D )^{t}\lambda _{1}^{\alpha}(x)-
t(\I -\D )^{t-1}\lambda _{0}^{\alpha}(x)\\* 
&=&\lambda _{1}^{\alpha}(x)+
\sum_{k\geq 2} \left({t \choose k } +t{t-1\choose k-1}\right)
\lambda _{1-k}^{\alpha}(x) / \rn{1-k}!\\*
&=&\lambda _{1}^{\alpha}(x)+
\sum_{k\geq 2} \left({t\choose k } +t{t-1\choose k-1}\right)
(-1)^{k}(k-2)! \lambda _{1-k}^{\alpha}(x),
\end{eqnarray*}
so
\begin{eqnarray*}
 G_{1}^{\alpha}(x)&=&\lr _{1}^{\alpha}+
 \sum_{k\geq 2}\left({t\choose k } +t{t-1\choose k-1}\right)
(-1)^{k}(k-2)! \lr _{1-k}^{\alpha}\\*
G_{1}^{(1)}(x)&=&x\log
(x)-x+\sum_{k\geq 2}B_{k}^{(2)}\romcoeff{2}{k}x^{1-k}+
 \sum_{k\geq 2}
\frac{(-1)^{k}\left({t\choose k }+t{t-1\choose
k-1}\right) (k-2)!}{(x-1)\cdots (x-k+1)}\\*
\end{eqnarray*}

\sssect{Lower Factorial to Gould}\index{Connection
Constants,Lower Factorial to Gould}\index{Lower
Factorial}\index{Gould Sequence}
Similarly, here we are interested in the operator $f(\D )=\D
(\I + \D)^{t}$ where $t$ is real. For $t=0$, we rederive the harmonic
logarithms $\lambda_{a}^{\alpha}(x)$, and for $t=-1$, we get a variant of the
Laguerre graded sequence $(-1)^{\pi ia}L_{a}^{\alpha }(x)$. (See
\eref{dot100} and \S\ref{LF-UF}.)

The standard associated graded sequence for $f(\D )$ is
\begin{eqnarray*}
r_{a}^{\alpha}(x)&=& f'(\D ) \left(\frac{f(\D )}{\D
}\right)^{-a-1;0} \lambda _{a}^{\alpha}(x)\\
&=&((t+1)\D +\I )(\I +\D )^{-at-1;0} \lambda_{a}^{\alpha}(x). 
\end{eqnarray*}
For $a\neq 0,1$,
\begin{eqnarray*}
r_{a}^{\alpha}(x) &=& \sigma (\I +\D )^{-at;0}
\lambda_{a-1}^{\alpha}(x)\\ 
&=&\sum_{\scriptstyle k\geq 0 \atop \scriptstyle k\neq a}
{-at\choose k} \frac{\rn{a-1}!}{\rn{a-k-1}!}\lambda _{a}^{\alpha}(x).
\end{eqnarray*}
Thus, for $a\neq 0,1$,
$$ G_{a}^{\alpha}(x) = \sum_{\scriptstyle k\geq 0 \atop
\scriptstyle k\neq n} {-at\choose k}
\frac{\rn{a-1}!}{\rn{a-k-1}!} \lr_{a-k}^{\alpha}.$$ 
Now, for $a=0$, 
\begin{eqnarray*}
r_{0}^{\alpha}(x)&=& ((t+1)\D +\I )(\I  +\D )^{-1}\lambda _{0}^{\alpha}(x)\\
&=& (1+t\D -t\D ^{2}+t\D ^{3}-\cdots )\lambda _{0}^{\alpha}(x)\\
&=&\lambda _{0}^{\alpha}(x)+t\sum_{k>0} (k-1)!
\lambda _{-k}^{\alpha}(x),
\end{eqnarray*}
so that again
$$ G_{0}^{\alpha}(x)=\lr
_{0}^{\alpha}(x)+t\sum_{k>0} (k-1)! \lr _{-k}^{\alpha}. $$
For $n=1$,
\begin{eqnarray*}
r_{1}^{\alpha}(x)&=& (\I +\D )^{-t}\lambda _{1}^{\alpha}(x)-
t(\I +\D )^{-t-1}\lambda _{0}^{\alpha}(x)\\ 
&=&\lambda _{1}^{\alpha}(x)+\sum_{k\geq 2}\left({-t\choose k
}+t(k-1){-t-1\choose k-1 
}\right) (-1)^{k}k! \lambda _{-1-k}^{\alpha}(x),
\end{eqnarray*}
hence another remarkable expansion for a residual series:
$$ G_{-1}^{\alpha}(x)=\sum_{k\geq 0} \left({t\choose k }+t{t-1\choose k-1
}\right) (-1)^{k} k! \lr _{-1-k}^{\alpha }.$$

As noted in \cite{ch2}, nearly any sequence may be used as a ``factorial'' on
which to base on umbral calculus. In particular, we could have chosen the
Gaussian coefficients as our factorials. If so, then would have derived the
$q$-analog of our theory.

{\begin{prob}
What are the $q$-analogs of these examples?
\end{prob}}

We hope the preceding examples display the utility of the
theory of formal power series of logarithmic type.

\end{document}